\documentclass[12pt,a4paper]{amsart}

\usepackage{amssymb}
\usepackage{amsthm}
\usepackage{tikz-cd}
\usepackage{enumerate}

\theoremstyle{plain}

\theoremstyle{remark}

\DeclareFontFamily{U}{cyr}{}
\DeclareFontShape{U}{cyr}{m}{n}{<5> wncyr5 <6> wncyr6 <7> wncyr7 <8> wncyr8 <9> wncyr9 <10->
wncyr10}{}
\DeclareMathAlphabet{\mathcyr}{U}{cyr}{m}{n}

\input cyracc.def

\newcommand{\ts}{\vspace{\baselineskip}\noindent{\bf Proof.}$\;\;$}
\newcommand{\ZZ}{{\mathbb Z}}
\newcommand{\QQ}{{\mathbb Q}}
\newcommand{\RR}{{\mathbb R}}
\newcommand{\CC}{{\mathbb C}}
\newcommand{\FF}{{\mathbb F}}

\newcommand{\PP}{{\mathbb P}}

\newcommand{\cN}{{\mathcal N}}

\newcommand{\cT}{{\mathcal T}}

\newcommand{\im}{\mbox{im}}

\newcommand{\End}{\mbox{End}}

\newcommand{\half}{\mbox{$\frac{1}{2}$}}

\newcommand{\Pic}{\operatorname{Pic}}
\newcommand{\Ext}{\operatorname{Ext}}
\newcommand{\Fix}{\operatorname{Fix}}
\newcommand{\MBM}{\operatorname{MBM}}
\newcommand{\SPE}{\operatorname{SPE}}
\newcommand{\Kthree}{\operatorname{K3}}

\newcommand{\bes}{\begin{equation*}}
	\newcommand{\ees}{\end{equation*}}

\begin{document}
\title[]{Hyperk\"ahler sixfolds, abelian fourfolds of Weil type and a Hodge class}

\author{Bert van Geemen}
\address{Dipartimento di Matematica, Universit\`a di Milano, Via Cesare Saldini 50, Milano, Italy}
\email{lambertus.vangeemen@unimi.it}

\author{Antonio Rapagnetta}
\address{Dipartimento di Matematica,  Universit\`a di Roma Tor Vergata, Via della Ricerca Scientifica 1, Roma, Italy}
\email{rapagnet@mat.uniroma2.it}

\begin{abstract}
There are now several proofs of the Hodge conjecture for the general abelian fourfold of Weil type with trivial discriminant. This paper provides another one.
The abelian fourfolds under consideration allow a map to a hyperk\"ahler sixfold of K3$^{[3]}$ type.
The pull-back of the second Chern class of the tangent bundle of the sixfold
is an algebraic class in codimension two that is not an intersection of divisor classes and the main result follows.

After recalling the basic facts on abelian fourfolds of Weil type
we establish the existence of the map using results on the birational geometry of these hyperk\"ahler manifolds and deformation theory.
\end{abstract}

\maketitle

\section*{Introduction}

Recently there has been important progress on the Hodge conjecture. 
This conjecture states that, for a projective complex manifold $X$, the subspace of Hodge classes
$$
H^{p,p}(X,\QQ)\,:=\,H^{2p}(X,\QQ)\cap H^{p,p}(X)
$$ 
is spanned by the classes of codimension $p$ subvarieties of $X$.
In particular E.\ Markman (\cite{Markman25}) proved that this conjecture holds for all abelian varieties of dimension at most five.

In an earlier paper, Markman proved a result on  abelian varieties of Weil type of dimension four.
These come in four dimensional families which are determined, up to isogeny, by an imaginary quadratic field $K$ which lies in the endomorphism algebra of all members of the family, and a discriminant.
There is a 2--dimensional subspace of  Hodge-Weil classes $HW(A)=HW(A,K)\subset H^{2,2}(A,\QQ)$
in the space of Hodge classes of such a fourfold $A$.
The field $K$ acts on $HW(A)$ and if one non-trivial class in $HW(A)$ is algebraic, then all classes in $HW(A)$ are algebraic.

We offer a new, simpler and more explicit, proof of Markman's theorem below in this paper.

\subsection{Theorem}\label{thm:main} (Markman \cite{Markman23})
For any abelian fourfold of Weil type $A$ with discriminant one,
the subspace $HW(A)\subset H^{2,2}(A,\QQ)$ of Hodge-Weil classes is spanned
by classes of algebraic cycles.

\

A simpler proof of Theorem \ref{thm:main}  was already given by S.\ Floccari and L.\ Fu in \cite{Floccari_F}.
It uses a map from 
fourfolds of Weil type to (singular) OG6 hyperk\"ahler varieties as well as 
other recent results on the Hodge conjecture.

Our proof is very much inspired by theirs. 
We also provide a rather direct, explicit, description of the algebraic cycles involved.
Moreover, we do not need the OG6 varieties but instead
we work with the better known K3$^{[3]}$ type hyperk\"ahler manifolds.
The second part of the paper is devoted to the proof of the existence of the manifolds we need.
Actually, these K3$^{[3]}$ type manifolds are natural double covers of the singular OG6 type varieties used \cite{Floccari_F} (see \cite{Mongardi_RS18}),
so this aspect of our proof is really very similar to the one in \cite{Floccari_F}.

\

In order to state the main, technical, result that we need, we recall some definitions.
The Kummer variety $K(A)$ of an abelian variety $A$ is the blow-up of the quotient
variety, the singular Kummer variety, $K^s(A):=A/\{\pm1\}$ in its singular locus, 
which consists of the images of the $2$--torsion points.
$$
A\,\longrightarrow\,K^s(A) \,=\,A/\pm 1_A\,\stackrel{}{\longleftarrow}\, K(A)~.
$$
A hyperk\"ahler manifold of K3$^{[3]}$ type is a deformation
of the Hilbert scheme of degree three cycles of dimension zero on a K3 surface.

Our simplified proof depends on the following theorem. 
The manifold $X$ `replaces' the singular OG6 variety in \cite[Theorem 1.3]{Floccari_F}.

\subsection{Theorem} \label{thm:WeilK33}
Let $B'$ be a very general abelian fourfold of Weil type, for an
imaginary quadratic field $K$, with trivial discriminant.

Then $B'$ is isogeneous to an abelian fourfold $B$ such that the Kummer variety
$K(B)$ of $B$ is a submanifold of a six-dimensional projective hyperk\"ahler manifold $X$ of K3$^{[3]}$ type whose transcendental Hodge structure $T_X\subset H^2(X,\QQ)$ has dimension six.

\

Our proof is based on the observation that any K3$^{[3]}$ type manifold $X$,
even a non-projective one, has a natural Hodge class in $H^4(X,\QQ)$.
It is the second Chern class of the tangent bundle $\cT_X$.
The following theorem thus implies Theorem \ref{thm:main},
we prove it, assuming Theorem \ref{thm:WeilK33}, in Section \ref{sec:prf03}.

\subsection{Theorem}\label{thm:c2TX}
Let $B'$ be a very general abelian fourfold of Weil type with imaginary quadratic field $K$ 
and discriminant one. Let
$$
B'\,\longrightarrow\,B\,\dasharrow\,K(B)\,\hookrightarrow\,X
$$
be the composition of the maps in Theorem \ref{thm:WeilK33}. Then the pull-back of the second Chern class $c_2(\cT_X)$ of the tangent bundle of $X$ to $B'$ is a Hodge-Weil class
which is also algebraic. Therefore all Hodge-Weil classes in $HW(B',K)$ are algebraic.

\

In more detail, to show that that $c_2(\cT_X)$ pulls-back to a Hodge-Weil class, we use its well-known relation with the class $q_X^\vee\in H^4(X,\QQ)$ defined by the Beauville-Bogomolov-Fujiki (BBF) form $q_X$ on $H^2(X,\QQ)$. This form provides a polarization on the Hodge structure of
$H^2(X,\QQ)$. The pull-back of the transcendental Hodge structure $T_X\subset H^2(X,\QQ)$
is a K3 type sub Hodge structure $T_1\subset H^2(B',\QQ)$ and $q_X$ induces a polarization
$q_1$ on $T_1$.
In Proposition \ref{prop:polHW} we observe that $q_1$ defines a Hodge-Weil class
$q_1^\vee\in H^4(B',\QQ)$. 
As $q_X^\vee$ is algebraic and pulls back to $q_1^\vee$
we conclude that the Hodge-Weil classes on $B'$ are algebraic.

\

Let $B'_0$ be the complement of the $2$--torsion points in the abelian variety $B'$.
The (trivial) tangent bundle of $B'_0$
is a subbundle of the pull-back of $\cT_X$ to $B'_0$. The quotient bundle $\cN_{0}$
is a rank two bundle whose second Chern class is thus a Hodge-Weil class. After tensoring $\cT_X$ or $\cN_0$ by a sufficiently ample line bundle, this rank two bundle will have a global section whose zero locus is a 2--cycle on $B'$ that, up to a complete intersection
of divisors, defines a Hodge-Weil class.
This observation might be useful to find explicit surfaces in $B'$
with a cohomology class that is not a complete intersection of divisor classes.

\

To construct a hyperk\"ahler manifold as in Theorem~\ref{thm:WeilK33}, 
we start with a specific example: the Hilbert scheme $S^{[3]}$ of length three subschemes on the smooth Kummer surface $S=K(A)$ associated with a general Jacobian $A$ of a genus two curve.

The surface $S$ is the minimal resolution of a nodal quartic surface $S'\cong K^s(A)$,
which is isomorphic to the quotient $A/\{\pm 1\}$ in $\mathbb{P}^3$. The linear
system of plane sections on $S'$ induces, via pullback, a linear system
$|L|$ of genus three curves on $S$.
The Hilbert scheme $S^{[3]}$ is birational to the Beauville-Mukai system $M\rightarrow |L|$,
it is a Lagrangian fibration whose fiber over a curve $C\in |L|$ 
parametrize rank one torsion free sheaves of degree $3$ on $C$.

The variety $M$ is equipped with a birational symplectic involution $\iota_M$ (which
was already used to compute the Hodge numbers of $OG6$ in \cite{Mongardi_RS18}). 
In fact, let $R$ be the square root of the line bundle associated with
the exceptional divisor of the desingularization $S \rightarrow S'$ (that is, the line
bundle defining the ramified double cover from the blow-up of $A$ at its
$2$--torsion points to $S$).

The restriction of $R$ to any smooth curve $C$ in $|L|$ is a point of order two in the Jacobian of $C$. 
It follows that tensoring by $R$ induces the translation by this point of order two on the fiber over $C$
of the Beauville-Mukai system and thus defines a rational involution $\iota_{M}$ on $M$.

We show that a $4$--dimensional component of the fixed locus of $\iota_M$
is birational to the Kummer fourfold $K(A \times A)$ of $A \times A$, and that,
after replacing $M$ by a suitable birational model if necessary, we may assume
that this component is isomorphic to $K(A \times A)$.

The smooth pair $(K(A \times A), M)$ deforms, at least locally, as a smooth
pair over the locus where the rational involution $\iota_M$ deforms. By the
Global Torelli Theorem, this locus coincides with the Hodge locus of the
anti-invariant lattice of $\iota_M$.

By ergodic  results of Ratner, or alternatively by the local Torelli Theorem and density results on the Hodge loci in the period domain, we deduce that the rank seven invariant sublattice of $\iota_M$ in $H^2(M,\ZZ)$ contains sublattices of rank six of any determinant $d>0$. 
If the period of a deformation of $M$ lies in the complexification of such a sublattice, 
then we show that the corresponding deformation of $K(A\times A)$ is the Kummer variety of an abelian
fourfold of Weil type with field $K=\QQ(\sqrt{-d})$ and discriminant one.
That one obtains a variety of Weil type can be deduced using Kuga-Satake varieties.
The triviality of the discriminant follows from the results
of Lombardo \cite{Lombardo} on the second cohomology group of an abelian fourfold of Weil type. 
We give an alternative proofs for some of these results.

\

\subsection*{Structure of the paper.}
In Section~\ref{s:HWclasses} we recall the definition and the basic properties of Weil Hodge classes on Weil abelian varieties as well as polarizations of weight-two Hodge structures.
This allows us to prove our main result, Theorem~\ref{thm:c2TX},
assuming Theorem \ref{thm:WeilK33}, in Section \ref{sec:prf03}.

In the remainder of the paper, we explain how to realize these Kummer fourfolds as subvarieties of such hyperk\"ahler manifolds. More precisely, they arise as components of the fixed locus of suitable birational involutions.

In Section~\ref{s:BWlat} we study an appropriate variant $BW$ of the Barnes-Wall lattice $BW_{16}$, which will appear as the anti-invariant lattice of the symplectic birational involutions considered above.
See also  the recent preprint \cite{Floccari_F26} for this lattice and applications related to this paper.

In Section~\ref{s:MBMclasess} we recall the notion of MBM classes on hyperk\"ahler manifolds and their role in governing the birational geometry. We then classify the MBM classes in $BW$, viewed as a sublattice of the 
second cohomology group of a K3$^{[3]}$ manifold (see Proposition~\ref{prop:walldiv}).

In Section~\ref{s:Bimero} we study non-projective hyperk\"ahler manifolds of K3$^{[3]}$ type with Picard lattice isometric to $BW$. We show that such manifolds contain $256$ projective $3$-spaces which can be contracted simultaneously (see Proposition~\ref{256P3I}), and that they are endowed with a symplectic birational involution, which becomes regular after the contraction (see Corollary~\ref{cor:indet}).

In Section~\ref{s:BimeroKummer} we construct,
via suitable birational modifications of the geometrical example described above,
a projective hyperk\"ahler manifold $M$ of K3$^{[3]}$ type endowed with a symplectic birational involution whose anti-invariant lattice is isometric to $BW$. We show that one component of the fixed locus of this involution is bimeromorphic to a  Kummer fourfold (see Proposition~\ref{prop:vfix}), and that $M$ admits a contraction of $256$ projective $3$--spaces to a variety $\widehat{M}$ on which the involution becomes regular (see Corollaries~\ref{256P3M} and~\ref{cor:involution-properties}).

In Section~\ref{s:compfix} we analyze the component of the fixed locus of the involution on $M$ which is bimeromorphic to a Kummer fourfold, and prove that it is in fact isomorphic to it (see Corollary~\ref{cor:Z}).

Finally, in Section~\ref{s:defth} we use deformation-theoretic arguments to show that, up to isogeny, the Kummer fourfold associated with a very general Weil abelian fourfold of discriminant one and arbitrary imaginary quadratic field is contained in a hyperk\"ahler manifold of K3$^{[3]}$ type (see Proposition~\ref{prop:defo} and Section \ref{prf:ratner}).


\section{Hodge-Weil classes} \label{s:HWclasses}

\subsection{Abelian fourfolds of Weil type}
Let $(B,K)$ be a four dimensional abelian variety of Weil type, here
$K\subset End(B)_\QQ:=End(B)\otimes_\ZZ\QQ$ is an imaginary quadratic field
such that the action, by pull-back, of any $x\in K$ on $H^{1,0}(B)$ has eigenvalues $x$ and $\bar{x}$, the complex conjugate of $x$, each with multiplicity two.
Combining this with the Hodge decomposition leads to an eigenspace decomposition
$$
H^1(B,\CC)\,=\,H^{1,0}(B)_\sigma\;\oplus\;H^{1,0}(B)_{\bar{\sigma}}\;\oplus\;
H^{0,1}(B)_\sigma\;\oplus\;H^{0,1}(B)_{\bar{\sigma}}~,
$$
in which each summand has dimension two and $\sigma,\bar{\sigma}$ are the two
embeddings of $K$ into $\CC$. We often identify $K$ with $\sigma(K)\subset\CC$
and then $\sigma(x)=x$, $\bar{\sigma}(x)=\bar{x}$.


\

\subsection{Proposition}
A very general abelian fourfold $(B,K)$ of Weil type has a rank one N\'eron Severi group:
$$
NS(B)_\QQ\,=\,H^{1,1}(B,\QQ)\,:=\,H^2(B,\QQ)\cap H^{1,1}(B)\,=\,\QQ\omega
$$
and the eigenvalue of $x\in K$ on a basis $\omega$ is $x\bar{x}$.
The Hodge classes in codimension two are a $3$--dimensional $\QQ$--vector space
$$
H^{2,2}(B,\QQ)\,=\,\QQ \omega^2\,\oplus\,HW(B,K)~. 
$$
The eigenvalue of $x$ on $\QQ\omega^2$ is  $(x\bar{x})^2$ and the eigenvalues
on the $2$--dimensional space of Hodge-Weil cycles are $x^4$, $\bar{x}^4$.
In particular, if a non-trivial class $c\in HW(B,K)$ is algebraic, then $HW(B,K)$
is spanned by the algebraic classes $c$ and $x^*c$ for a general $x\in K$.

\ts See \cite{Weil}, \cite[Thm 6.12]{vG_hcav}. We recall some aspects of the proof.
Notice that since $NS(B)$ has rank one,
one must have $x^*\omega=x\bar{x}\omega$. 

The eigenspace decomposition
of $H^1$ induces decompositions of all $H^q=\wedge^qH^1$, in particular
$$
H^4(B,K)\,=\,\wedge^4H^1(B,K)\,=\,
\bigoplus_{k=0}^4
\wedge^kH^1(B,K)_\sigma\otimes \wedge^{4-k}H^1(B,K)_{\bar{\sigma}}~,
$$
is the decomposition into $K$--eigenspaces with eigenvalue $x^k\bar{x}^{4-k}$ on the $k$--th summand.
There are two $1$--dimensional $K$--subspaces, the summands with $k=0,4$,
these are conjugates. This allows us to define a $2$--dimensional $\QQ$--subspace:
$$
HW(B,K)\,\subset\,H^4(B,\QQ),\qquad HW(B,K)\otimes_\QQ K\,=\,
\wedge^4H^1(B,K)_\sigma\oplus  \wedge^4 H^1(B,K)_{\bar{\sigma}}~.
$$
The eigenvalues of $x\in K$ on $HW(B,K)$ are $x^4,\bar{x}^4$ and any $K$-eigenvector
in $H^4(B,\CC)$ with eigenvalue $x^4$ or $\bar{x}^4$ lies in $HW(B,K)\otimes_\QQ\CC$.
We do observe that also $HW(B,K)\subset H^{2,2}(B)$:
$$
HW(B,K)\otimes\CC\,=\,
\big(\wedge^2H^{1,0}(B)_\sigma\big)\otimes (\wedge^2H^{0,1}(B)_\sigma\big)\;\oplus\;
\big(\wedge^2H^{1,0}(B)_{\bar{\sigma}}\big)\otimes (\wedge^2H^{0,1}(B)_{\bar{\sigma}}\big)~,
$$
since each summand has Hodge type $(2,2)$ and $K$--eigenvalue $x^4$, $\overline{x}^4$, respectively.
Thus $HW(B,K)$ lies in the space of Hodge classes $H^{2,2}(B,\QQ)$.
\qed

\

\subsection{The discriminant} \label{sec:defdis}
Let $(B,K)$ be a general abelian fourfold of Weil type
and let $\omega\in NS_\QQ(B)\subset \wedge^2H^1(B,\QQ)$ be the class of an ample divisor.
Then $\omega$ defines an alternating form on the dual vector space
$H_1(B,\QQ)=H^1(B,\QQ)^\vee$. Considering $H_1(B,\QQ)$ as a
4-dimensional $K=\QQ(\sqrt{-d})$--vector space,
one defines a non-degenerate Hermitian form $H$ on it (so
$H(v,w)=\overline{H(w,v)}$  and $H$ is $K$-linear in the second factor):
$$
H:\;H_1(B,\QQ)\,\times\,H_1(B,\QQ)\,\longrightarrow\,K,\qquad
H(v,w)\,=\,\omega(v,\sqrt{-d}_*w)\,+\,\sqrt{-d}\omega(v,w)~.
$$
This hermitian form has signature $(2,2)$.

The discriminant of $(B,K,\omega)$ is the determinant of the Hermitian matrix that defines $H$
w.r.t.\ some $K$-basis of $H_1(B,\QQ)$. Changing the basis changes this determinant by
$Nm(s):=s\bar{s}$ for some $s\in K^\times$,
hence the discriminant is well-defined in the group
$\QQ^\times/Nm_{K/\QQ}(K^\times)$ (\cite[Lemma 5.2]{vG_hcav}). 

Two non-degenerate hermitian $K$-vector spaces $(V_1,H_1),(V_2,H_2)$ of the same dimension and the same signature are isomorphic if and only if they have same discriminant.

\subsection{Complete families}
A complete family of $2n$--dimensional abelian varieties of Weil type is determined by a $2n$--dimensional non-degenerate hermitian $K$-vector space $(V,H)$, of signature $(n,n)$, and a free abelian subgroup
$\Gamma\subset V$ of rank $4n$, so that $\Gamma\otimes_\ZZ\QQ\cong V$ and $\omega:=\im(H)$
is integer valued on $\Gamma$.
The complex structures $J$ on $V_\RR=V\otimes_\QQ\RR$ which commute with $K$ and such that
the Riemann conditions are satisfied for $\omega$, so $\omega(Jv,Jw)=\omega(v,w)$ and $\omega(v,Jv)>0$ for all $v,w\in V_\RR$, $v\neq 0$,
define the members of the family; $B_J:=V_\RR/\Gamma$ is an abelian variety of Weil type, polarized by $\omega$ (cf.\ \cite{Weil}, \cite{vG_hcav}).

A polarized abelian variety of Weil type $(B,K,\omega)$ determines such a complete family by choosing 
$V=H_1(B,\QQ)$, $\Gamma=H_1(B,\ZZ)$ and the Hermitian form $H$ is defined by $\omega$ as in Section \ref{sec:defdis}.

\

\noindent
The $2n$--dimensional abelian varieties of Weil type considered in this paper, as well as those in \cite{Markman23} and \cite{Markman25}, are obtained as deformations of a self product $A\times A$ where $A$ is a general $n$--dimensional abelian variety. The following lemma shows that the discriminant of these deformations is $(-1)^n$.

\subsection{Lemma}\label{lem:trivdisc}
Let $B$ be an $2n$-dimensional abelian variety of Weil type with field $K$ and polarization $\omega$. Assume that the complete family determined by $(B,K,\omega)$ contains an abelian variety $B'$ which is isogeneous to $A\times A$, where $A$ is a general abelian variety. Then the discriminant of $(B,K,\omega)$ is $(-1)^n$.

\ts 
Since $A$ is general, we may assume that $A$ is simple and has Picard number one.
The isogeny induces an isomorphism $H_1(B',\QQ)\cong H_1(A,\QQ)^2$.
Let $K=\QQ(\sqrt{-d})$ and let $f:=(\sqrt{-d})_*\in \End(H_1(B',\QQ))$. Then, for any copy of $H_1(A,\QQ)\subset H_1(B,\QQ)$, $f(H_1(A,\QQ))\cap H_1(A,\QQ)=0$ since otherwise either $f\in \End(A)_\QQ$,
which contradicts $\End(A)_\QQ=\QQ$, or $A$ is not simple. Thus, for any copy of $H_1(A,\QQ)\subset H_1(B',\QQ)$, we find
$$
H_1(B',\QQ)\,=\,H_1(A,\QQ)\,\oplus\,f(H_1(A,\QQ)),\qquad f:\,(a,f(b))\,\longmapsto\,(-db,f(a))~,
$$
where we used that $f^2=-d$. Let $p:B\rightarrow A$ be the map induced by the projection on the first factor.

The rank of the N\'eron Severi group of $A\times A$ is three, and it contains
$\omega_A:=p^*\omega$ and $f^*\omega_A$, where $\omega$ is an ample generator of $NS(A)$. Notice that
$f^*(f^*\omega_A)=[-d]^*\omega_A=d^2\omega_A$ since $[-d]^*$ acts as $d^2$ on
$H^2(A,\QQ)=\wedge^2H^1(A,\QQ)$.
We define $\omega_K:=d\omega_A+f^*\omega_A$, then
$$
f^*\omega_K\,=\,f^*(d\omega_A\,+\,f^*\omega_A)\,=\,d^2\omega_A\,+\,df^*\omega_A\,=\,d\omega_K~,
$$
hence $\omega_K\in NS(A)_\QQ$ is an eigenvector for the $\sqrt{-d}$--action with eigenvalue $d$. Similarly
$d\omega_A-f^*\omega_A$ has eigenvalue $-d$. On the other hand, the eigenvalues of $x\in K$ on $H^2(B',\QQ)$ are $x^2,\bar{x}^2$ and $x\bar{x}$.
Hence $x=\sqrt{-d}$ has eigenvalues $-d$ and $d$,
and the multiplicity of the eigenvalue $-d$ on $NS(A)_\QQ$ must be even.
We conclude that $\omega_K$ is a polarization of Weil type on $(B',K)$ and it is the unique one
up to scalar multiple.
Therefore the Hermitian form $H$ which determines the complete family is defined as in 
Section \ref{sec:defdis} with $\omega=\omega_K$, up to a positive scalar multiple.

Now we show that there is an $H$--isotropic subspace of $K$-dimension $n$ in $(H_1(B',\QQ),K)\cong K^{2n}$.
Since $\omega_A$ is alternating on $H_1(A,\QQ)$, there is a subspace $W_A\subset H_1(A,\QQ)$ of $\QQ$-dimension $n$ such that $\omega_{A|W_A\times W_A}=0$. Let $W:=W_A+f(W_A)\subset H_1(B',\QQ)$, then $W$ is $K$--vector space with $K$--dimension $n$ and, as $\omega_K=d\omega_A+f^*\omega_A$, it is easy to check that $H_{|W\times W}=0$.

Let $W'$ be a complementary subspace to $W$, so that $V=W\oplus W'$.
Choose a $K$--basis of $(H_1(B',\QQ),K)$
such that $e_1,\ldots,e_n$ generate $W$ and the $e_{n+1},\ldots,e_{2n}$ generates $W'$.
The Gram matrix $M$ of $H$ then has a $2\times 2$ block form with matrices $A,\ldots,D$ and
$A=0$, $C={}^t\overline{B}$. Therefore $\det(M)=(-1)^n\det(B)\overline{\det(B)}\in (-1)^n\cdot Nm_{K/\QQ}(K^\times)$, hence the discriminant of $(V,H)$ is $(-1)^n$.

For completeness sake we observe that if
$(V,H)$ be a $2n$-dimensional non-degenerate hermitian $K$-vector space of signature $(n,n)$, then
there is a $K$-basis of $V$ such that, for some $a\in\QQ_{>0}$,
$$
H(z,w)\,=\,az_1\bar{w}_1\,+\,z_2\bar{w_2}+\ldots+\,z_n\bar{w}_n\,-\,\big(z_{n+1}\bar{w}_{n+1}\,+\,\ldots\, 
+\,z_{2n}\bar{w}_{2n}\big)~.
$$
If the discriminant of $H$ is $(-1)^n$, we may assume that $a=1$  and then
$$
W\,:=\,\langle\ldots, e_j+e_{n+j},\ldots\rangle_{j=1,\ldots,n}
$$
is a maximal isotropic subspace. Thus the discriminant of $H$ is $(-1)^n$ if and only if there exists an $n$--dimensional $K$-subspace $W\subset V$ such that $H_{|W\times W}=0$.
\qed

\

\

\noindent
The following lemma will be used to show that certain families of abelian fourfolds are of Weil type in Section \ref{prf:ratner}. Notice that the proof avoids the (explicit) use of Kuga-Satake varieties.

The assumption that these families are deformations of the self-product of an abelian surface
implies by Lemma \ref{lem:trivdisc} that the discriminant of these Weil type deformations is trivial,
but this also follows, without this assumption, from the results of Lombardo
(see \cite[Corollary 3.6]{Lombardo}, cf.\ Theorem \ref{thm:lomb} below).

\subsection{Lemma}\label{lem:noKS}
Let $B$ be an abelian fourfold and assume that there exists a polarized sub Hodge structure $T\subset H^2(B,\QQ)$
with Hodge numbers $\dim T^{2,0}=1$ and $\dim T^{1,1}=4$. If $T$ is general, then $B$ is an abelian variety of Weil type.

In case the complete family of fourfolds of Weil type determined by $B$ has a member isogeneous to
$A\times A$, for a general abelian surface $A$,
the imaginary quadratic field is $K=\QQ(\sqrt{-d})$, with $d=\det(q_T)$,
here $q_T$ is the polarization on $T$ viewed as a quadratic form, and the discriminant is one.

\ts 
Let $G$ be the special Mumford-Tate group of $T$, it is an algebraic subgroup of $SL(T)$.
 The fact that $T$ is general and polarized implies that $G_\CC\cong SO(6)$. 
In a suitable basis of $T_\CC$, the homomorphism $h$ defining the Hodge structure on $T$ induces $h_\CC:\CC^\times\rightarrow SO(6)$, $h_\CC(t)=diag(t^2,1,1,1,1,t^{-2})$ with eigenvalue $t^{p-q}$ on $T^{p,q}$. 

As $T\subset \wedge^2H^1(B,\QQ)$, there is a homomorphism of the special Mumford-Tate group $SMT(B)$ of $B$, that is, of $H^1(B,\QQ)$, onto $G$. The (reductive) Lie algebra $smt(B)_\CC$ of $SMT(B)_\CC$ thus  has a summand which is the Lie algebra $so(6)$ of $SO(6)$.

The homomorphism $h$ lifts to a homomorphism $\tilde{h}$ and the eigenvalues of $\tilde{h}_\CC(t)$ on $H^1(B,\CC)$ must be $t,t^{-1}$. 
The Lie algebras $so(6)$ and $sl(4)$ (of $SL(4)$) are isomorphic. Considering the (co-)weight lattice of $sl(4)$, the only possible lift is given, on the complex points, by 
$\tilde{h}:\CC^\times\rightarrow SL(4)$, $\tilde{h}(t)=(t,t,t^{-1},t^{-1})$. 
The standard representation $\CC^4$ of $SL(4)$ (or its dual $(\CC^4)^\vee$) 
is thus a summand of the $SMT(B)_\CC$ representation on $H^1(B,\CC)$.
The representation of $SL(4)$ on $T_\CC$ is the one on $\wedge^2\CC^4$.

As $B$ is an abelian variety, $H^1(B,\QQ)$ is polarized by a non-degenerate alternating form $E$ and $SMT(B)\subset Sp(E)$.
However, there is no $SL(4)$--invariant bilinear form on $\CC^4$. One finds that to obtain an 
$SL(4)$--invariant anti-symmetric form, the one induced by a polarization of $B$, one must have 
$$
H^1(B,\CC)\,\cong\,\CC^4\,\oplus\,(\CC^4)^\vee
$$ 
as $SL(4)$--representations.
The representation of $SMT(B)$ on $H^1(B,\QQ)$ is irreducible as 
none of the subrepresentations of $SL(4)$ is obtained from a polarized sub Hodge structure of $H^1(B,\QQ)$. Thus $B$ is a simple abelian fourfold.

The algebraic group $Sp(E)_\CC\cong Sp(8)$ has rank 4. As $SL(4)$ has rank 3, there can be at most one copy of $sl(4)$ in $smt(B)_\CC$. The other summand, if any, is a rank one reductive Lie algebra. In particular, there is a double cover $\tilde{G}\rightarrow G$ and $\tilde{G}\subset SMT(B)$ is a normal subgroup defined over $\QQ$. The homomorphism $\tilde{h}$, defining the Hodge structure on $H^1(B)$,
has values in $\tilde{G}$, hence by minimality of the special Mumford-Tate group, $\tilde{G}=SMT(B)$, it is a $\QQ$--form of $SL(4)$.

The endomorphism algebra $\End(B)_\QQ$ of $B$ consists of those endomorphisms of  $H^1(B,\QQ)$ that commute with the action of $SMT(B)$, so $\End(B)_\CC$ consists
of the $\CC$--linear maps commuting with the $SL(4)$ action, which by Schur's lemma are the maps that are scalar multiplications on both irreducible components $\CC^4$ and $(\CC^4)^\vee$, hence
$\dim \End(B)_\QQ=2$. 
The two $SL(4)$--subrepresentations of $H^1(B,\CC)$ are thus the eigenspaces of the action of 
$\End(B)_\QQ$. 
The eigenvalues of $\tilde{h}(t)$, $\tilde{h}^\vee(t)=\tilde{h}(t)^{-1}$ on the eigenspaces
are $t,t^{-1}$ with multiplicity two, hence the intersections of the two eigenspaces with $H^{1,0}$ 
both have dimension 2.

Since  $B$ is simple, $K:=\End(B)_\QQ$ must be a quadratic field extension of
$\QQ$. If it were a totally real field, the Picard number of $B$ would be two (see [Prop.\ 5.5.7]\cite{LB})
but the Picard number is the dimension of the space of $SL(4)$--invariants in 
$$
H^2(B,\CC)\,=\,\wedge^2(\CC^4)\,\oplus \,\wedge^2(\CC^4)^\vee\,\oplus\,\CC^4\otimes
(\CC^4)^\vee~.
$$
The first two summands are irreducible (and $SL(4)$ acts through its quotient $SO(6)$ on them)
whereas the last summand is $\End(\CC^4)=\CC \mbox{id}_{\CC^4}\oplus sl(4)_\CC$, and the representation of $SL(4)$ on $sl(4)_\CC$ is the adjoint representation which is irreducible.
Hence the subspace of invariants is one dimensional and is spanned by $\mbox{id}_{\CC^4}$
and the Picard number of $B$ is one. Therefore $\End(B)_\QQ$ is an imaginary quadratic field.

We already observed that the $K$--eigenspaces in $H^1(B,\CC)$ are the $SL(4)$--subrepresentations
and that each of these intersects $H^{1,0}(B)$ in a two dimensional subspace.
So $H^{1,0}(B)$ is the direct sum of two $2$--dimensional $K$--eigenspaces, and
hence $B$ is of Weil type. 
That the discriminant is one follows from Lemma \ref{lem:trivdisc} (or \cite[Corollary 3.6]{Lombardo}).

From \cite[Thm 3.8.3]{Lombardo}, for $B$ general of Weil type with field $K=\QQ(\sqrt{-d})$ and trivial discriminant, for any polarized rank six Hodge substructure $T\subset H^2(B,\QQ)$,
the quadratic space $(T,q_T)$ is isometric to
$U\oplus U\oplus \langle-2\rangle\oplus\langle-2d\rangle$, which has determinant $4d$ and obviously
$\QQ(\sqrt{-4d})=\QQ(\sqrt{-d})$.
\qed

\

\subsection{The second cohomology group of an abelian fourfold of Weil type} \label{TH^2}
In Lemma \ref{lem:noKS} we showed that an abelian fourfold $B$ having a sub Hodge structure with Hodge numbers $(1,4,1)$ in $H^2(B,\QQ)$ must be of Weil type, with discriminant one. 
Lombardo studied the Hodge structure of a general abelian fourfold of Weil type in detail.
We collect some of the results of Lombardo in the following theorem.

\subsection{Theorem}(Lombardo \cite{Lombardo})\label{thm:lomb}
Let $(B,K)$ be a very general abelian fourfold of Weil type with
discriminant $a\in \QQ^\times /Nm_{K/\QQ}(K^\times)$.
\begin{enumerate}
\item[a)] The action of $K$ on $H^2(B,\QQ)$ has, over $\QQ$, the eigenspace decomposition
$$
H^2(B,\QQ)\,=\,H^2(B,\QQ)_0\,\oplus\,T,\qquad \dim H^2(B,\QQ)_0\,=\,16,\quad \dim T\,=\,12~,
$$
where the eigenvalues of $x\in K$ on $H^2(B,\QQ)_0$ are $x\bar{x}$ and they are
$x^2,\bar{x}^2$ on $T$.
\item[b)] One has $\dim\big(H^2(B)_0\big)^{2,0}=4$, so $\dim \big(H^2(B)_0\big)^{1,1}=8$. 
This sub Hodge structure is a direct sum of $\QQ \omega$,
where $\omega\in H^{1,1}(B,\QQ)$ is the first Chern class of an ample line bundle on $B$, and a simple $15$-dimensional complement.
\item[c)] One has $\dim T^{2,0}=2$, so $\dim T^{1,1}=8$.
The algebra $End_{Hod}(T)$ of Hodge endomorphisms of $T$ (the $\QQ$-linear endomorphisms mapping $T^{p,q}$ into $T^{p,q}$ for all $p,q$) is a quaternion algebra and it is a skew field except if the discriminant
$a$ of $(B,K)$ is one, in that case $End_{Hod}(T)\cong M_2(\QQ)$.
\item[d)] In case the discriminant $a=1$, there exists a K3 type Hodge structure $T_1$ such that
$$
T\,\cong\, T_1^{\oplus 2},
\qquad \dim T_1^{2,0}\,=\,1,\quad \dim T_1^{1,1}\,=\,4,\quad End_\QQ(T_1)\,=\,\QQ~,
$$
and $T_1$ is uniquely determined by $(B,K)$.
\end{enumerate}

\

\ts
We recall some aspects of the proof from \cite{Lombardo}, using the notation introduced in the proof of Lemma \ref{lem:noKS}. 
For any $K$ and any discriminant, the $smt(B)_\CC$--representation on $T_\CC$
is the direct sum of $\wedge^2V$ and $\wedge^2V^\vee$,
where $V=\CC^4$ is the standard representation of $SL(4)$. 
As $\dim V=4$, there is a
perfect pairing $\wedge^2V\times\wedge^2V\rightarrow\wedge^4V\cong\CC$ and thus $\wedge^2V\cong \wedge^2V^\vee$.
Therefore $T_\CC\cong \big(\wedge^2V)^{\oplus 2}$ and $End(T)_\CC$ is isomorphic to the matrix algebra 
$M_2(\CC)$.

The subtle point is to see whether there exist, over $\QQ$, non-trivial $smt(B)$--representations in $T$;
these subrepresentations correspond to sub Hodge structures.
A careful analysis, using a natural Hodge endomorphism of $T$, shows that this happens if and only if the discriminant is trivial, \cite[Cor.\ 3.6]{Lombardo}.
The existence of a K3 type sub Hodge structure in that case can also be deduced from \cite[Section 6]{vG23}.
\qed

\

\noindent 
We recall the definition of a polarization of a weight two Hodge structure and that  
these polarizations define Hodge classes. 
In the case of fourfolds of Weil type with discriminant one, we show in Proposition \ref{prop:polHW} that
the polarizations of the rank six sub Hodge structures in $H^2$ define Hodge-Weil classes.

\subsection{Polarized Hodge structures and Hodge classes}
A (rational) Hodge structure on a (finite dimensional) $\QQ$-vector space $V$
is of K3 type if $V$ has a (Hodge) decomposition
$$
V_\CC\,:=\,V\otimes_\QQ\CC\,=\,V^{2,0}\,\oplus\,V^{1,1}\,\oplus\,V^{0,2},\qquad
\overline{V^{p,q}}\,=\,V^{q,p},\quad\dim V^{2,0}\,=\,1
$$
and a polarization, which is a non-degenerate quadratic form $q_V:V\rightarrow\QQ$
such that
$$
V^{p,q}\,\perp\,V^{p'q'}\quad\mbox{if}\quad (p+p',q+q')\neq (2,2)~,
$$
where perpendicular refers to the bilinear form $b_{V,\CC}$ associated to $q_{V,\CC}$.
In particular $q_{V,\CC}$ is zero on $V^{2,0}$, $V^{0,2}$ and is non-degenerate on
the subspaces $V^{2,0}\oplus V^{0,2}$ and $V^{1,1}$.
We will not need positivity properties of $q_V$ on $V_\RR$ in this section.

Thus $q_V\in Sym^2 V^\vee$, the second symmetric power of the dual $V^\vee$ of $V$.
The bilinear form $b_V$ induces an isomorphism $V\rightarrow V^\vee$ and its inverse gives
an isomorphism $Sym^2V^\vee\rightarrow Sym^2V$. We denote by $q_V^\vee\in Sym^2V$ the image of $q_V$, we call it the polarization class (of $(V,q_V)$).
The vector space $Sym^2V$ has a natural Hodge structure
with $(Sym^2V)^{r,s}$ being the direct sum of the
$V^{p,q}V^{p',q'}$ with $r=p+p'$, $s=q+q'$.
The following lemma is well-known and a variant holds more generally for any polarized  Hodge structure.

\subsection{Lemma}\label{lem:pol22} Let $(V,q)$ be a K3 type polarized Hodge structure.
Then the polarization class $q_V^\vee\in Sym^2V$
is a rational Hodge class, that is,
$$
q_V^\vee \,\in\, (Sym^2V)\,\cap\,(Sym^2V)^{2,2}~.
$$
Assume that the special Mumford-Tate group of $V$ is
the orthogonal group $SO(V,q)$. Then any Hodge class in $Sym^2V$ is a scalar multiple of $q^\vee$:
$$
(Sym^2V)\,\cap\,(Sym^2V)^{2,2}\,=\,\QQ q_V^\vee~.
$$

\ts We only need to show that $q_V^\vee\in (Sym^2V)^{2,2}$.
Choosing a $\CC$-basis $e_j$ of $V_\CC$ such that each $e_j$ lies in one of the $V^{p,q}$,
let $e_j^*:=b_{V,\CC}(e_j,\cdot)\in V^\vee_\CC$. 
Then $q_\CC$ is a linear combination of the $e_j^* e_k^*$ with $e_j,e_k\in V^{1,1}$
or $e_j \in V^{2,0}$ and $e_k \in V^{0,2}$. Therefore $q_{V,\CC}^\vee$ is the same linear combination of the $e_je_k$, all of which have Hodge type $(2,2)$, so $q_V^\vee\in (Sym^2V)^{2,2}$.

The sub Hodge structures of $Sym^2 V$ correspond to the subrepresentations of
the special Mumford-Tate group, or equivalently, those of its Lie algebra. The Lie algebra representation of $so(V,q)_\CC$ on $Sym^2 V_\CC$ is the direct sum of two irreducible components, one of which is one dimensional, cf.\ \cite[Excercise 19.21]{Fulton_H}.
\qed

\

\noindent
The following proposition relates any copy of $T_1$ in $T\subset H^2(B,\QQ)$, for a Weil type fourfold $B$ with discriminant one, to the Hodge-Weil classes.

\

\subsection{Proposition} \label{prop:polHW}
Let $B$ be a very general abelian fourfold of Weil type with discriminant one and
let $T_1\hookrightarrow H^2(B,\QQ)$ be a(ny) sub Hodge structure of K3 type. 
Then the image of the cup product map
$$
Sym^2(T_1)\,\longrightarrow\,H^4(B,\QQ)
$$
contains an exceptional Hodge class. More precisely, let $q_1^\vee\in Sym^2(T_1)$ denote the Hodge class defined by the polarization $q_1$ on $T_1$, then the image of $q_1^\vee$ is in $HW(B,K)$, so it is a
Hodge-Weil class.

\

\noindent
{\bf Proof.}$\quad$Recall that $HW(B,K)$ is the subspace of $H^4(B,\QQ)$ on which $x\in K$ acts
with eigenvalues $x^4,\bar{x}^4$.
The map $Sym^2H^2(B,\QQ)\rightarrow H^4(B,\QQ)$ is surjective and
$$
Sym^2H^2(B,\QQ)\,=\, Sym^2H^2(B,\QQ)_0\,\oplus\,H^2(B,\QQ)_0\otimes T\,\oplus\,Sym^2T~.
$$
The eigenvalues on the first and second summand are $(x\bar{x})^2$ and $(x\bar{x})x^2,(x\bar{x})\bar{x}^2$ respectively.
Since the eigenvalues on $HW(B,K)$ are $x^4,\bar{x}^4$,
this implies that the intersection of the image of the first two summands
with $HW(B,K)$ is zero. Therefore $HW(B,K)\subset \mbox{im}(Sym^2T)$, more precisely, $HW(B,K)$ lies in the image of the space of Hodge classes in $Sym^2T$.

As $End_{Hod}(T_1)=\QQ$, if $x \notin\QQ$, we have $xT_1\neq T_1$.
Since $T_1$ is simple, the intersection between $xT_1$ and $T_1$ must be zero, so $T=T_1\oplus xT_1$.
Therefore
$$
Sym^2T\,=\,Sym^2(T_1)\;\oplus\;Sym^2(xT_1)\; \oplus\; T_1\otimes(xT_1)~.
$$
By assumption, the special Mumford-Tate group of $T_1$ is $SO(T_1,q_1)$,
where $q_1$ is the polarization on $T_1$. This implies by Lemma \ref{lem:pol22}
that the Hodge classes in the first two summands are spanned by their polarizations.
Multiplication by $x$ gives an isomorphism $T_1\rightarrow xT_1$, hence
as representations of the special Mumford-Tate group of $B$ we get
$$
T_1\otimes (xT_1)\,\cong\,T_1\otimes T_1\,\cong\,Sym^2(T_1)\,\oplus\,\wedge^2T_1~.
$$
The representation on $\wedge^2T_1$ is isomorphic to the adjoint representation of the Lie
algebra of $SO(T_1,q_1)$, which is irreducible, hence there are no Hodge classes in that summand. The space of Hodge classes in $Sym^2(T_1)$ is again 1-dimensional by Lemma
\ref{lem:pol22}.
Thus $Sym^2T$ has a three dimensional space of Hodge classes.

The subspace $HW(B,K)\subset \mbox{im}(Sym^2T)$ is mapped into itself under the action of $x^*$ for any $x\in K$.
Assuming that $\mbox{im}(Sym^2(T_1))\cap HW(B,K)=0$
then implies that also $\mbox{im}(Sym^2(xT_1))\cap HW(B,K)=0$.
But now the image of the space of Hodge classes in $Sym^2T$ intersects $HW(B,K)$ in a subspace of dimension at most one, which contradicts that the 2-dimensional
$HW(B,K)\subset \mbox{im}(Sym^2T)$.
Therefore $\mbox{im}(Sym^2(T_1))\cap HW(B,K)\neq 0$, and the same holds for any choice of sub Hodge structure $T_1'\subset T$ with $T_1'\cong T_1$.
Since any Hodge class in $Sym^2(T_1)$ lies in $\QQ q_1^\vee$,
it follows that $q_1^\vee$ maps to a non-zero class in $HW(B,K)$.
\qed

\

\subsection{Remark} For completeness sake we observe that the action of $K$ does not preserve the summands
of $Sym^2(T)$ given in the proof of Proposition \ref{prop:polHW}.
The action preserves $T$ and thus also the 3-dimensional space of Hodge classes in $Sym^2(T)$. Since $HW(B,K)\subset\mbox{im}(Sym^2(T))$, both $x^4$ and $\bar{x}^4$
are eigenvalues.
Since the dimension of this space is three, the remaining eigenvalue must be $x^2\bar{x}^2$.

\subsection{A Hodge class on a K3$^{[3]}$ type manifold}\label{sec:hodgeclassK33}
We refer to \cite[\S 2]{Rizzo} for a summary with references of the results collected here.
Let $X$ be a HK manifold of dimension six deformation equivalent to the Hilbert cube of
a K3 surface. The Beauville-Bogomolov-Fujiki (BBF) form $q_X$ is a quadratic form on
$H^2(X,\ZZ)$ which defines a polarization of this K3 type Hodge structure of dimension $23$.
There is an isomorphism of Hodge structures
$$
H^4(X,\ZZ)\,\cong\,Sym^2H^2(X,\ZZ)\,\oplus\,H^2(X,\ZZ)(-1)
$$
where the second component has a Tate twist so that it is a weight four Hodge structure.

The class $q^\vee_X\in Sym^2H^2(X,\ZZ)$ defined by the BBF form is,
up to a non-zero scalar multiple, the second Chern class of $X$:
$$
q^\vee_X\,=\,\lambda c_2(X),\qquad \lambda\,\in\,\QQ,\quad \lambda\,\neq\, 0~.
$$
In fact, following O'Grady \cite[\S 2]{O'Grady}, by deforming $X$ to a very general non-projective manifold,
we may assume that there are no Hodge classes in $H^2(X,\QQ)$.
Then the special Mumford-Tate group of $H^2(X,\QQ)$ is the special orthogonal group $SO(q_X)$ of the BBF form $q_X$.
Therefore the only non-trivial Hodge classes in $Sym^2H^2(X,\QQ)$ are the scalar multiples of $q_X^\vee$.
Also the second Chern class $c_2(\cT_X)\in H^4(X,\QQ)$ is a non-trivial Hodge class
(in fact $c_2(\cT_X)^3 =36800\neq 0$), and it is thus proportional to $q_X^\vee$.

\subsection{Proof of Theorem \ref{thm:c2TX}}\label{sec:prf03}
We now prove Theorem \ref{thm:c2TX}, which implies Theorem \ref{thm:main},
assuming Theorem \ref{thm:WeilK33}.
Let $B'$ be a general abelian fourfold of Weil type with an isogeny $B'\rightarrow B$
such that  $K(B)\subset X$ for a projective HK manifold $X$ of K3$^{[3]}$ type.
Then the pull-back of a non-zero holomorphic 2-form
on $X$ to $B$ must be a non-zero holomorphic 2-form on $B$
(else the tangent space in a point $x$ to $K(B_0)$
would be isotropic in $T_xX$ which is impossible for dimension reasons).
Hence we obtain an injective map from the transcendental lattice $T_X$ of $X$
into $H^2(B,\ZZ)$.
Tensoring with $\QQ$, we see that the six dimensional image of $T_{X,\QQ}$ must lie in the sub Hodge structure $T$ defined in Theorem \ref{thm:lomb} since $T_{X,\QQ}$ is of K3 type and the
only sub Hodge structures of K3 type in $H^2(B,\QQ)$ lie in $T$.
We will denote the image by $T_1$.

This induces an isomorphism of rational Hodge structures $Sym^2(T_X)\rightarrow Sym^2(T_1)$.
The Hodge class $q^\vee_X$ maps to a non-zero element in $Sym^2(T_1)$
which must be a scalar multiple of $q_1^\vee\in Sym^2(T_1)$.
By Proposition \ref{prop:polHW}, $q_1^\vee$ maps to a non-zero
Hodge-Weil class in $HW(B,K)\subset H^4(B,\QQ)$.
Thus $q^\vee_X$ pulls back to a Hodge-Weil class.
Since $X$ is projective, $c_2(\cT_X)$ is algebraic and since
$q^\vee_{X}$ and $c_2(X)$ are proportional, also $q^\vee_X$ is algebraic.
Hence the pull-back of that algebraic cycle represents
a Hodge-Weil class, which is thus also algebraic.

As observed earlier, this implies, using the $K$--action, that
then all Hodge-Weil classes in $B$ are algebraic. Pulling back along the isogeny
$B'\rightarrow B$ we conclude that the Hodge-Weil classes of $B'$ are also algebraic.
\qed

\

\subsection{Generalizations} An important aspect of the proof is the fact that the
rational map $B\rightarrow X$ allows one to pull-back $H^2(X)$ to $H^2(B)$, which then defines
an inclusion of Hodge structures $T_{X,\QQ}\hookrightarrow H^2(B,\QQ)$.
In case the discriminant of a general abelian fourfold of Weil type $B$ is non-trivial, there is no sub Hodge structure of K3 type in $H^2(B,\QQ)$, see Theorem \ref{thm:lomb},
hence in this case there cannot exist such a map to a hyperk\"ahler manifold.

A general abelian variety $B$ of Weil type of dimension $2n$ for $n>2$ does not have Hodge structures of K3 type in $H^2(B,\QQ)$, so it is not clear how to generalize the results from this paper. 
In fact, the Special Mumford-Tate group of $B$ over $\CC$ is $SL(2n)$
and the representation on $H^1(B,\CC)$ is the direct sum of the standard representation
and its dual. This implies that there are three simple sub Hodge structures in $H^2(B,\QQ)$,
these have $h^{2,0}=2\binom{n}{2},n^2,0$ respectively.

It is maybe interesting to observe that if the discriminant of $B$ is trivial,
then $H^n(B,\QQ)$ does have a sub Hodge structure $T_1$ of CY type, so with $T_1^{n,0}=1$,
see \cite{Cacciatori_F} for the case $n=3$.


\

\section{The Barnes-Wall lattice and an involution}\label{s:BWlat}

\subsection{The Kummer lattice}\label{sec:Kumlat}
Let $A$ be an abelian surface. The quotient $K^s(A)=A/\langle -1_A\rangle$ is a surface with sixteen nodes, the images of the 2--torsion points,
whose desingularization is a K3 surface $K(A)$.
The 2--torsion points $A[2]$ form a vector space over the field of two elements $\FF_2$
of dimension four, $A[2]\cong\FF_2^4$.
Let $N_x$ be the class of the exceptional curve over the node defined by $x\in A[2]$.
These exceptional curves generate a non-primitive sublattice of rank $16$ of
$H^2(K(A),\ZZ)$.
The Kummer lattice $Kum$ is smallest primitive sublattice of $H^2(K(A),\ZZ)$
that contains all the sixteen $N_x$. More explicitly, in \cite{Nikulin} it is shown that
$Kum$ is generated by the (classes of) the
$$
N_x,\qquad N_x^2\,=\,-2,\quad N_x\cdot N_y\,=\,0\quad(x\neq y)~,
$$
and the classes
$$
\half\sum_{x\in W} N_x
$$
where $W\subset A[2]$ is either a subgroup of order eight or its complement.
There are $15$ such subgroups and thus there are $30$ possible $W$'s.
Notice that there is a well-defined homomorphism $Kum\rightarrow\ZZ$ sending
$\sum a_xN_x$ to $\sum a_x$, the sum is over all $x\in A[2]$ with $a_x\in\half\ZZ$,
but of course only the elements in $Kum$ are considered.
Reducing modulo two we get a surjective homomorphism
$Kum\rightarrow\ZZ/2\ZZ$. We learned the following definition from a lecture of Mukai,
see also \cite[\S 4.10]{Conway_S} for other definitions.

\subsection{Definition}\label{def:BWlat}
The Barnes-Wall lattice $BW$ is the index two sublattice of the Kummer lattice that is the kernel of the surjective homomorphism
$$
BW\,=\,BW_{16}(-1)\,:=\ker\Big({Kum}\,\longrightarrow\,\ZZ/2\ZZ,\quad \sum a_xN_x\,\longmapsto\,\sum a_x\mod 2\Big)~.
$$
In particular, the $N_x$ are not in $BW$ and in fact $BW$ does not have classes $c$ with
$c^2=-2$. The lattice $BW$ is negative definite, even, of rank sixteen.

\subsection{Properties of $BW$}\label{sec:propBW}
We list some of the properties of the Barnes-Wall lattice.
\begin{enumerate}
\item[a)]
The lattice $BW$ has discriminant group $A_{BW}:=BW^*/BW\cong\FF_2^8$, its discriminant quadratic
form can be identified with a quadratic form $q:A_{BW}\rightarrow\FF_2$.
One has
$$
q(x+y)\,=\,q(x)\,+\,q(y)\,+\,b(x,y)~,
$$
where the bilinear alternating form $b$ (i.e.\ $b(x,x)=0$ for all $x$)
is non-degenerate, so $(A_{BW},b)$ is a symplectic vector space over $\FF_2$.
On a suitable symplectic basis $e_i$ of $A_{BW}$ (so $b(e_i,e_j)=0$ except if $j=i\pm4$
and $b(e_i,e_{i\pm 4})=1$) the discriminant
form $q$ is given by
$q(\sum x_ie_i)=\sum_{i=1}^4x_ix_{i+4}$, in particular it has Arf invariant
$\sum q(e_i)q(e_{i+4})=0$ (cf.\ \cite[\S 5.1.2]{Dolgachev}).

\item[b)] The orthogonal group $O(BW)$ of $BW$ has a normal subgroup $H_8$ of order $2^9$
and the quotient is the orthogonal group $O(A_{BW},q)\cong O^+_8(\FF_2)$ of the discriminant group.
The action of $O(BW)$ on the discriminant group is through this quotient.
The group $H_8$ is generated by $N_x\mapsto N_{x+y}$ and by
$N_x\mapsto (-1)^{b(x,y)}N_x$ for any fixed $y\in A[2]$, the center of $H_x$ is $\pm I$, where $I$ is the identity on $BW$.

\item[c)] The classes $N_x\in Kum$ have scalar products $N_x\cdot w\in\ZZ$ for all
$w\in BW$. Since they are not in $BW$, they define non-trivial classes
$\bar{N}_x \in A_{BW}$ with $q(N_x)=0$. Actually 32 classes $\pm2N_x\in BW$ form a single $H_8$ orbit and they
define a unique class $\bar{N}_x=(\pm2N_x)/2\in A_{BW}$, in fact $\pm N_x\pm N_y\in BW$.
The $O(BW)$-orbit of $2N_x$ has cardinality $4320=32\cdot 135$
and the elements in this orbit map onto the 135 non-zero
elements $a$ in $A_{BW}$ with $q(a)=0$.

\item[d)] The lattice $BW$ is generated by the classes $\pm N_x\pm N_y$ and the
$\half\sum_{x\in W} \pm N_x$ with an even number of minus signs. These are all the vectors of length $-4$ in $BW$, 
there are $16\cdot15\cdot 2+30\cdot 2^7=480+3840=4320$ classes of length $-4$.
\end{enumerate}

\subsection{The $BW$ and the Mukai lattice}\label{sec:BWM}
The discriminant group $A_{BW}$ of $BW$ is isomorphic to the discriminant lattice
of the even, rank eight, lattice $U(2)^{\oplus 4}$ of signature $(4,4)$,
here $U(2)$ is defined by the quadratic form on $\ZZ^2$ given by $(x,y)\mapsto 4xy$.
Let $\phi:A_{BW}\rightarrow A_{U(2)^{\oplus 4}}$ be such an isomorphism. Then we glue
$BW$ to $U(2)^{\oplus 4}$ by adding to $BW\oplus U(2)^{\oplus 4}$ the eight
vectors $\half(\tilde{e}_i\oplus \tilde{f}_i)$, $i=1,\ldots,8$, where the
$\tilde{e}_i\in BW$ map to a basis $e_i$ of $A_{BW}$ and the $\tilde{f}_i\in U(2)^{\oplus 4}$ map to $\phi(e_i)\in A_{U(2)^{\oplus 4}}$.

In this way we obtain an overlattice $\widetilde{H}$ of $BW\oplus U(2)^{\oplus 4}$
which has trivial discriminant group and thus $\widetilde{H}$ is even, unimodular, of signature $(4,20)$. This overlattice is therefore unique up to isometry
(\cite[Cor.\ 1.13.3]{NikulinL}):
$$
BW\oplus U(2)^{\oplus 4}\;\hookrightarrow\;
\widetilde{H}\;\cong\; U^{\oplus 4}\oplus E_8(-1)^{\oplus 2}~.
$$

The lattice $\widetilde{H}$ is known as the Mukai lattice,
it is isomorphic to $\oplus_i H^i(S,\ZZ)$ for a K3 surface $S$ equipped
with the quadratic form $(r,x,s)^2:=-2rs+x^2\in H^4(S,\ZZ)=\ZZ$,
the latter identification is due to the fact that $S$ is a compact complex surface, hence it is an oriented 4-manifold.

The embedding $BW\oplus U(2)^{\oplus 4}\subset \widetilde{H}$ is unique up to isomorphism,
in fact, such an embedding induces an isomorphism $\phi'$ between the discriminant groups
and $\phi^{-1}\circ\phi'$ then induces isomorphisms on the two discriminant groups. As any such isomorphism is induced by an isomorphism of $BW$ and $U(2)^{\oplus 4}$
respectively, the uniqueness follows.

\subsection{The cohomology of K3$^{[3]}$ type manifolds}\label{cohK33}
Up to isometry, there is a unique vector $v\in \widetilde{H}$ with $v^2=4$,
see \cite[Thm.\ 1.14.4]{NikulinL} applied to the sublattice $\ZZ v$ of $\widetilde{H}$.
Thus we may assume that $v=(1,2)\in U$, where $U$ is the quadratic form on $\ZZ^2$
with $(x,y)^2=2xy$, and it is, say, the last copy of $U$ in $U^{\oplus 4}
\subset \widetilde{H}$ under the isomorphism from Section \ref{sec:BWM}.
Then $v^\perp\cap U$ is spanned by $(1,-2)$ and
$$
v^\perp\,\cong \,E_8(-1)^2\,\oplus \,U^3\,\oplus\,\langle -4\rangle~,
$$
notice that $v^\perp$ is thus isometric to $H^2(X,\ZZ)$, with the BBF-form, for any manifold $X$ of K3$^{[3]}$ type, we denote this lattice also by $\Lambda_{\Kthree^{[3]}}$.

The vector $w=(1,1)\in U(2)\subset U(2)^{\oplus 4}\oplus BW$ also has $w^2=4$. 
This implies that there is a copy of the Barnes-Wall lattice in contained in $v^\perp$ 
and its perpendicular in $v^\perp$ is the lattice 
$$
BW^{\perp_{v^\perp}}\;\cong\; U(2)^{\oplus 3}\,\oplus\,\,\langle -4\rangle~.
$$

\subsection{The involution on $v^\perp$}\label{sub:invret}
The involution on $BW\oplus U(2)^{\oplus 4}$ that is $-1$ on $BW$ and
$+1$ on $U(2)^{\oplus 4}$ extends to an involution on $\widetilde{H}$.
The vector $w$ in Section \ref{cohK33} is invariant under this involution and thus
the involution on $\tilde{H}$ restricts to an involution $\iota_v$ on $v^\perp$, the invariant and anti-invariant lattices are:
$$
(v^\perp)^{\iota_v=1}\,=\,BW,\qquad (v^\perp)^{\iota_v=-1}\,=\,U(2)^{\oplus 3}\,\oplus\,\,\langle -4\rangle~.
$$

\subsection{Remark}
We observe that a manifold of K3$^{[3]}$ type with such an involution cannot be
isomorphic to a Hilbert cube $S^{[3]}$ for a K3 surface $S$ in such a way that the involution fixes the divisor of non-reduced subschemes.

In fact, assuming that this happens, there is an identification
$$
v^\perp\,=\, H^2(S,\QQ)\oplus\ZZ\eta\,\cong\, H^2(S^{[3]},\ZZ)~,
$$
where $\eta\in v^\perp$ would be a primitive class with $\eta^2=-4$
and divisibility $4$, and moreover, $\eta$ would be invariant under the involution.
Hence
$\eta\in U(2)^{\oplus 4}\cap v^\perp$. However, since $\eta\in U(2)^{\oplus 4}$
is primitive, there is a $b\in BW$ such that $(b+\eta)/2\in \widetilde{H}$ and then in fact
$(b+\eta)/2\in v^\perp$. But now $\eta(b+\eta)/2=\eta^2/2=2$, in contradiction with the
fact that the divisibility of $\eta$ in $v^\perp$ is $4$.

\section{MBM classes and birational geometry}\label{s:MBMclasess}

\subsection{MBM classes and wall divisors}\label{sec:walldivisors}
The birational geometry of a hyperk\"ahler manifold $X$ of K3$^{[3]}$ type is,
as for any hyperk\"ahler manifold, governed by MBM classes (or wall divisors) in
$H^2(X,\ZZ)$ of type $(1,1)$. 
For completeness, we briefly recall their definition and basic properties; 
these classes were introduced independently by Amerik and Verbitsky (see \cite{AVimrn})
and Mongardi (see \cite{Mongardi15}).

\subsection{Definition} \label{def:mbm}\begin{enumerate}\item Let $X$ be a hyperk\"ahler manifold. A class
\[
\alpha \in H^2(X,\mathbb{Z})
\]
with negative BBF square is called a Monodromy Birationally Minimal (MBM) class  or wall divisor  if for every deformation
$(X',\alpha')$ of $(X,\alpha)$ through families of hyperk\"ahler manifolds
such that $\alpha'$ is of type $(1,1)$ and $\mathrm{Pic}(X') \simeq \mathbb{Z}$,
the hyperk\"ahler manifold $X'$ contains a rational curve.
\item  An MBM class $\alpha$ on a hyperk\"ahler manifold $X$ is called primitive positive (ppMBM) class if $\alpha\in H^2(X,\mathbb{Z})$ is indivisible, $\alpha\in \Pic(X)$ and $B_X(\alpha, \omega)>0$ for some (any) K\"ahler class $\omega$.
\end{enumerate}

Following \cite[Definition~5.6]{MaSPE}, in Definition \ref{def:mbm}, a deformation $(X',\alpha')$ of $(X,\alpha)$ is a parallel transport operator, induced by a smooth family over a connected base, sending $\alpha$ to $\alpha'$.

We denote by $\operatorname{MBM}(X)$ and $\operatorname{ppMBM}(X)$ the sets of MBM classes and ppMBM classes of the hyperk\"ahler manifold $X$, respectively.

\subsection{Remark}\label{rem:mbmdef}
By \cite[Theorem~5.16]{AVimrn}, in order to prove that a class $\alpha \in H^{2}(X,\mathbb{Z})$ with negative BBF square is MBM, it suffices to exhibit a single deformation $(X',\alpha')$ of $(X,\alpha)$ through a family of hyperk\"ahler manifolds such that $\alpha' \in \Pic(X') \cong \mathbb{Z}$ and $X'$ contains a rational curve. Furthermore, by \cite[Theorem~5.11]{AVimrn}, the existence of a rational curve on $X'$ is equivalent to the existence of an effective curve on $X'$.

By definition, given two hyperk\"ahler varieties $X_1$ and $X_2$, any parallel transport operator
\[
t:H^{2}(X_1,\mathbb{Z})\rightarrow H^{2}(X_2,\mathbb{Z})
\]
arising from a family of hyperk\"ahler varieties induces a bijection between the sets of MBM classes $\operatorname{MBM}(X_1)$ and $\operatorname{MBM}(X_2)$. On the other hand, if $\alpha_1$ is a ppMBM class of $X_1$, then, even when $\alpha_2:=t(\alpha_1)$ is of type $(1,1)$, one can only assert that one of the two classes $\pm \alpha_2$ is ppMBM on $X_2$. Indeed, for the parallel transport operator induced by a Mukai flop, the class $t(\alpha_1)$ is the negative of a ppMBM class.

Finally, assume moreover that $\alpha_1$ and $\alpha_2$ are both ppMBM classes. By \cite[Lemma~5.17]{MaSPE}, there also exists a family of hyperk\"ahler varieties realizing a deformation from $(X_1,\alpha_1)$ to $(X_2,\alpha_2)$ such that the class corresponding to $\alpha_1$ remains of type $(1,1)$ throughout the family.

\subsection{MBM classes, rational curves and contractions}\label{sub:mbmratcon}

To clarify the relationship between MBM classes on the hyperkähler variety $X$ and rational curves on its deformations, we recall that the BBF form on $X$ induces a natural  isomorphism

\begin{equation}\label{j}
j_X : H_{2}(X,\mathbb{Q}) \longrightarrow H^{2}(X,\mathbb{Q}).
\end{equation}
This isomorphism is defined by requiring that $j_X(\gamma)$ is the unique class
in $H^{2}(X,\mathbb{Q})$ such that
\[
\int_{\gamma} \beta = B_X\bigl(j_X(\gamma), \beta\bigr)
\]
for every $\beta \in H^{2}(X,\mathbb{Q})$.

We first consider the case $\Pic(X)\cong \mathbb{Z}$. Let $\alpha \in \Pic(X)$ be an MBM class. Then, for every curve class $\gamma \in H_2(X,\mathbb{Q})$, the class $j_X(\gamma)$ is proportional to $\alpha$. Moreover, if $\alpha$ is primitive, then, by definition, $\alpha$ is a ppMBM class if and only if there exist a positive rational number $k$ and an effective curve class $\gamma$ such that
\[
j_X(\gamma)=k\alpha.
\]

Moreover, if the second Betti number of $X$ is greater than $5$, then by \cite[Corollary~5.9]{BLjems} (see also \cite[Theorem~4.6]{AVsel}) there exists a contraction
\[
X \,\longrightarrow\, \widehat{X}
\]
onto a singular K\"ahler variety $\widehat{X}$, which contracts precisely the curves whose classes are mapped by $j_X$ to multiples of $\alpha$. As $\Pic(X)\cong \mathbb{Z}$, every curve class has this property, and hence all curves on $X$ are contracted.

Finally, let $X_1$ and $X_2$ be hyperk\"ahler manifolds with
\[
\Pic(X_1)\cong \Pic(X_2)\cong \mathbb{Z},
\]
and let $\alpha_i\in H^2(X_i,\mathbb{Z})$, $i=1,2$, be MBM Hodge classes. Assume that there exists a parallel transport operator arising from a family of hyperk\"ahler manifolds and mapping $\alpha_1$ to $\alpha_2$. Let
\[
X_i\,\longrightarrow\, \widehat{X}_i
\]
denote the contraction associated with $\alpha_i$. Then the singular varieties $\widehat{X}_1$ and $\widehat{X}_2$ are locally trivially deformation equivalent.

This follows directly from \cite[Corollary~5.9]{BLjems} when both $\alpha_1$ and $\alpha_2$ are both ppMBM classes. In the general case, the statement follows again from \cite[Corollary~5.9]{BLjems} together with the fact that, since $\alpha_i$ generates $\operatorname{Pic}(X_i)$, the manifold $X_i$ admits at most two hyperk\"ahler birational models, and the contractions associated with these models have the same image.

The only case needed in the present paper is that of manifolds of $\mathrm{K3}^{[3]}$ type. In this situation, the local trivial deformation equivalence of $\widehat{X}_1$ and $\widehat{X}_2$ is proved in \cite[Corollary~4.7]{Amerik_V}.

The discussion above in the Picard rank one case has the following immediate consequences in general.

\begin{enumerate}
\item[(a)] Since effective curves specialize to effective curves (possibly reducible), the ppMBM classes of a hyperk\"ahler manifold $X$ can equivalently be defined as the MBM classes $\alpha\in H^2(X,\mathbb{Z})$ for which there exist an effective curve class $\gamma$ on $X$ and a positive rational number $k$ such that
\[
j_X(\gamma)=k\alpha.
\]

\item[(b)] Let $(X,\alpha)$ and $(X',\alpha')$ be as in the definition of an MBM class, and assume moreover that the second Betti number of $X$ is greater than $5$. Then the locally trivial deformation type of the contraction $\widehat{X}'$ of $X'$ depends only on the monodromy orbit of $\alpha$, namely its orbit under the action of parallel transport operators from $H^2(X,\mathbb{Z})$ to $H^2(X,\mathbb{Z})$ arising from families of hyperk\"ahler manifolds.
\end{enumerate}

Since a contraction of a hyperk\"ahler manifold is divisorial if and only if
the singular locus of the contracted variety has codimension $2$, item b) allows
 to define an important subclass of MBM classes (see \cite{MaSPE} and \cite{Msurv}).

\subsection{Definition}
An MBM class $\alpha$ of a hyperk\"ahler manifold $X$    is  a Stably Prime
Exceptional (SPE) class if the associated contraction (on deformations with cyclic Picard group) is divisorial, $\alpha\in \Pic(X)$ and $B_X(\alpha,\omega)>0$ for some (any) K\"ahler class $\omega$.

We denote by $\SPE(X)$ the set of SPE classes of $X$-

\subsection{MBM classes and cones}
To illustrate how MBM classes govern the birational geometry of a
hyperk\"ahler manifold $X$, we introduce the standard notation.

We denote by $\mathcal{K}_X$ the K\"ahler cone of $X$, and by
\[
\mathcal{C}_X := \{ \beta \in H^{1,1}(X,\mathbb{R}) \mid
B(\beta,\beta) > 0 \text{ and } B(\beta,\omega) > 0 \text{ for some }
\omega \in \mathcal{K}_X \}
\]
the positive cone of $X$.

We also define the birational K\"ahler cone of $X$ as
\[
\mathcal{BK}_X := \bigcup_{Y \text{ hyperk\"ahler},\atop \text{bimeromorphic to } X}
\mathcal{K}_Y,
\]
where the union is taken inside $\mathcal{C}_X$ via the pullback
identification of $H^2(Y,\mathbb{R})$ with $H^2(X,\mathbb{R})$.

Finally, we denote by $\mathcal{FE}_X$ the fundamental exceptional chamber,
that is, the interior in $\mathcal{C}_X$ of the closure of $\mathcal{BK}_X$.

The following result, due to  Amerik-Verbitsky (see \cite[Theorem 1.19]{AVimrn}) and Markman (see  \cite[Proposition 5.6 and Theorem 6.17]{Msurv}), describes the
K\"ahler cone and the fundamental exceptional chamber of a hyperk\"ahler
manifold in terms of MBM classes and SPE classes, respectively.

\subsection{Theorem} \label{thm:coni}
Let $X$ be a hyperk\"ahler manifold, and let $\MBM(X)$ and $\SPE(X)$ denote
the sets of MBM classes and SPE classes of $X$, respectively.

\begin{enumerate}
\item 
The K\"ahler cone $\mathcal{K}_X$ of $X$ is the
connected component of
\[
\mathcal{C}_X \setminus 
\bigcup_{\alpha \in \MBM(X) \cap \Pic(X)} \alpha^{\perp}
\]
containing a K\"ahler class.

\item The fundamental exceptional chamber $\mathcal{FE}_X$ is the connected
component of
\[
\mathcal{C}_X \setminus 
\bigcup_{\alpha \in \SPE(X)} \alpha^{\perp}
\]
containing a K\"ahler class.
\end{enumerate}

Since any integral multiple of an MBM class is again an MBM class and defines the same orthogonal hyperplane, the set of hyperplanes $\alpha^{\perp}$, with $\alpha \in \operatorname{MBM}(X)\cap\operatorname{Pic}(X)$, is in natural bijection with $\operatorname{ppMBM}(X)$. Hence, by Theorem~\ref{thm:coni}(1),  the K\"ahler cone of $X$ is also the
connected component of
\[
\mathcal{C}_X \setminus 
\bigcup_{\alpha \in \operatorname{ppMBM}(X)} \alpha^{\perp}
\]
containing a K\"ahler class.

\subsection{MBM classes and monodromy for $K3^{[3]}$ manifolds}
As already observed in Remark~\ref{rem:mbmdef}, the set $\operatorname{MBM}(X)$ of MBM classes of a hyperk\"ahler manifold $X$ is invariant under the action of parallel transport operators from $H^2(X,\mathbb{Z})$ to itself arising from families of hyperk\"ahler manifolds. In other words, the monodromy group of $X$ acts on $\operatorname{MBM}(X)$ and, in particular, on the subset of primitive elements of $\operatorname{MBM}(X)$. Hence the set of primitive elements of $\operatorname{MBM}(X)$ is a union of monodromy orbits.

In the case of interest for us, namely hyperk\"ahler manifolds of
K3$^{[3]}$ type, the monodromy group coincides with the full group of
orientation-preserving isometries. Hence, the orbits of primitive classes
under the action of the monodromy group are determined by their square and
divisibility with respect to the  
BBF form on $H^2(X,\mathbb{Z})$.

The orbits of primitive MBM classes (wall divisors) on a hyperk\"ahler manifold $X$ of K3$^{[3]}$ type
are listed in \cite[\S~2.2, Table]{Mongardi15} (see also
\cite[Table~H3, p.~309]{Hassett_T_intnum} and \cite[Cor.~4.7]{Amerik_V}).
They are the  classes $\alpha\in H^{2}(X,\mathbb{Z})$ satisfying
\[
(\alpha^2,div(\alpha)) = (-12,2),\; (-36,4),\; (-2,1),\; (-4,4),\; (-4,2).
\]
Moreover, in the first case the associated contraction on a hyperk\"ahler
manifold with cyclic Picard group contracts a projective space of dimension
$3$ to a point; in the second case, it contracts the fibers of a
$\mathbb{P}^2$--bundle over a K3 surface; and in the last three cases the
contractions are divisorial, i.e.\ these give rise to SPE classes.

\subsection{Proposition} \label{prop:walldiv}
Let $X$ be a hyperk\"ahler sixfold of K3$^{[3]}$ type with
$\Pic(X) = BW$, the Barnes-Wall lattice defined in Section \ref{def:BWlat}. 
Any ppMBM class $D \in \Pic(X)$ satisfies $D^2 = -12$ and has divisibility $2$; in particular, there are no  SPE classes on $X$.
There are $256$ such divisor classes. If $D,D'$ are distinct ppMBM classes of $X$, then
$b_X(D,D')\in\{-4,0,4\}$.

\ts
We first exclude the existence of classes $\alpha\in BW$ in
all cases listed in Section \ref{sec:walldivisors} except $(-12,2)$.
Since $\alpha^2\leq -4$ for all $\alpha\in BW$, the case $(-2,1)$ is excluded.
The discriminant group of $BW$ is $(\ZZ/2\ZZ)^8$, so the divisibility
$div_{BW}(x)=1,2$ for $x\in BW$ hence also $div_{H^2(X,\ZZ)}(x)=1,2$. This excludes
the cases $(-36,4)$ and $(-4,4)$.
It is easy to verify that the divisibility of a(ny) vector with $D^2=-4$ is 1,
this excludes the last two cases.

Now we consider the remaining case $(-12,2)$.
With Magma we found that there are three orbits of vectors of length $-12$ in $BW$.
Only one orbit consists of vectors that have divisibility two in $BW$,
it is characterized by its cardinality $61440=2^{12}\cdot3\cdot5$.

For a primitive MBM class $D$ one needs that $(D-v)/2$ is in the Mukai lattice, where
$v$ has $v^2=4$ and $v^\perp=\Lambda_{\Kthree^{[3]}}$.
As the class of $v/2\in A_{U(2)^{\oplus 4}}$ is non-trivial,
this class determines the class of $\bar{D}\in A_{BW}$ by the choice of the isomorphism $\phi$
between these two discriminant groups, we will now identify both groups with $A_{BW}$.

The class $\bar{D}\in A_{BW}$ has $q(\bar{D})=(D/2)^2=12/4=3\equiv 1\mod 2$.
The group $O(A_{BW})$ acts transitively on the set of the
$120$ elements $x\in A_{BW}$ with $q(x)=1$. There are thus
$61440/120=2^9$ elements in the orbit of the $(-12,2)$ classes mapping to $\bar{v}$,
these are one orbit of the subgroup $H_8$ of $O(BW)$.
As $D$ and $-D$ have the same image in $A_{BW}$, we conclude that there are
$2^8$ ppMBM classes.

The ppMBM classes can be described as follows. 
There are 16 subsets $S_1,\ldots,S_{16}\subset(\ZZ/2\ZZ)^4$, each with $6$ elements,
such that the $256$ MBM classes, with norm $-12$, are given by
$$
\pm N_{x_1}\pm N_{x_2}\ldots \pm N_{x_6},\qquad S_i\,=\,\{x_1,\ldots,x_6\}~.
$$
for certain choices of signs.
One choice of $x_i\in\FF_2^4$ is $S_1:=\{e_1,e_2,e_1+e_2,e_3,e_4,e_3+e_4\}$ where
the $e_i$ are the standard basis of $\FF_2^4$ and $D_1:=\sum N_x$ where we sum over $x\in S_1$.
Obviously $D_1\cdot (\pm N_x\pm N_y)\in2\ZZ$ and one checks that for any subspace $W$ of
codimension one also $D_1\cdot e_W\in 2\ZZ$, in fact the cardinality of $S\cap W$ is either
$2$ or $4$. Then also $D_1\cdot e_W\in 2\ZZ$ where now $W$ is the complement of codimension
one subspace.

There is one overall choice of sign. If $D,D'$ differ by an odd number of signs, then
$D-D'$ is twice the sum of an odd number of the $N_x$, so $D/2$ and $D'/2$ define distinct elements in $A_{BW}$ since the sum of an odd number of $N_x$ is not in $BW$. Since the primitive MBM classes, 
divided by two, map to a fixed $\bar{D}$, this implies that there are only sixteen choices for the sign. 
In fact, $16\cdot 16=256$.

Any two of these subsets $S_i$ have $2$ elements in common. 
This implies that if $D\ne D'$ then $b_X(D,D')=\pm 2\pm 2\in\{-4,0,4\}$.
We observe that the $16$ subsets $S_i$ correspond to the 16 tropes of $K^s(A)$
for a principally polarized, indecomposable, abelian surface $A$;
each trope is a hyperplane section of the 16-nodal quartic surface $K^s(A)\subset\PP^3$
which is a double conic and this conic contains six of the nodes.
\qed

\

\section{Bimeromorphic involutions on non-projective K3$^{[3]}$ type varieties}\label{s:Bimero}

In this section, we first consider irreducible holomorphic symplectic manifolds of
K3$^{[3]}$ type whose Picard lattice is isomorphic to the Barnes-Wall lattice
$BW$. We prove that every such non-projective manifold $X$ admits a symplectic bimeromorphic
involution. More precisely, we show that each of these manifolds possesses a
singular birational model $\widehat{X}$ on which the involution becomes regular; this model
is obtained by contracting $256$ projective $3$--spaces.

\subsection{Birational Geometry of K3$^{[3]}$ Varieties with Picard Lattice $BW$}
Before discussing the birational properties of K3$^{[3]}$ varieties with Picard
lattice $BW$, we begin with a simple but important consequence of Proposition \ref{prop:walldiv}: any two
projective spaces contained in such a variety and arising from distinct MBM 
classes must be disjoint. In order to state the result we need the natural isomorphism $j_X:H_2(X,\mathbb{Q})\rightarrow H^2(X,\mathbb{Q})$ induced by the BBF form introduced in (\ref{j}) in 
Section \ref{sub:mbmratcon}.

\subsection{Lemma}\label{lemmaintvuota}
Let $X$ be a hyperk\"ahler manifold of K3$^{[3]}$ type with
$\Pic(X) \simeq BW$, and let $P_1, P_2 \subset X$ be projective
spaces of dimension $3$. Let $\gamma_i \in H_2(X,\mathbb{Q})$ be the
class of a line in $P_i$, and set $\alpha_i := 2j_X(\gamma_i)$ for $i=1,2$.
Assume that $\alpha_1 \neq \alpha_2$ and that each $\alpha_i$ is a 
ppMBM class. Then the intersection product
\[
[P_1]\cdot [P_2] \in H^{6}(X,\mathbb{Z})
\]
is zero. In particular, $P_1 \cap P_2 = \varnothing$.

\ts By \cite[Theorem 1]{Harvey_HT}, the cohomology 
class of each projective space $P_i$ is:
\[
[P_i]=\frac{\alpha_i^{3}+c_2(X)\,\alpha_i}{48}~.
\]
Hence
\[
[P_1]\cdot[P_2]
=
\frac{1}{48^{2}}\Bigl(
\int_X \alpha_1^{3}\alpha_2^{3}
+
\int_X \alpha_1^{3}\alpha_2 \, c_{2}(X)
+
\int_X \alpha_1\alpha_2^{3} \, c_{2}(X)
+
\int_X \alpha_1\alpha_2 \, c_{2}(X)^{2}
\Bigr).
\]

By the Fujiki type formulas for K3$^{[3]}$'s (see page 279 of \cite{Harvey_HT}),
for any $\alpha\in H^{2}(X)$, writing $B$ for the BBF form:
\[
\int_X \alpha^{6} = 15\, B(\alpha,\alpha)^{3},\qquad
\int_X \alpha^{4}c_2(X) = 108\, B(\alpha,\alpha)^{2},\qquad
\int_X \alpha^{2}c_2(X)^{2} = 1200\, B(\alpha,\alpha).
\]

The polarizations of these formulas give:
\[
6!
\int_X \alpha_1\cdots\alpha_6
=
15\sum_{\sigma\in S_6}
B(\alpha_{\sigma(1)},\alpha_{\sigma(2)})
B(\alpha_{\sigma(3)},\alpha_{\sigma(4)})
B(\alpha_{\sigma(5)},\alpha_{\sigma(6)}),
\]
\[4!
\int_X \alpha_1\alpha_2\alpha_3\alpha_4\,c_2(X)
=
108
\sum_{\sigma\in S_4}
B(\alpha_{\sigma(1)},\alpha_{\sigma(2)})
B(\alpha_{\sigma(3)},\alpha_{\sigma(4)}),
\]
\[2
\int_X \alpha_1\alpha_2 c_2(X)^2 = 1200 \, \sum_{\sigma\in S_2}B(\alpha_{\sigma(1)},\alpha_{\sigma(2)}).
\]

From Proposition \ref{prop:walldiv}, we know:
\[
B(\alpha_1,\alpha_1)=B(\alpha_2,\alpha_2)=-12,\qquad
B(\alpha_1,\alpha_2)\in\{0,\pm 4\}.
\]

In case  $B(\alpha_1,\alpha_2)=0$, all the mixed Fujiki integrals vanish,
and therefore $[P_1]\cdot[P_2]=0$.

In the second case, where $B(\alpha_1,\alpha_2)=\pm 4$,
by evaluating the polarizations, we obtain:
\[
\int_X \alpha_1^{3}\alpha_2^{3}
=
\frac{15\cdot 3! 2^3}{6!}\Bigl(6 B(\alpha_1,\alpha_2)^3+ 9 B(\alpha_1,\alpha_1) B(\alpha_2,\alpha_2)B(\alpha_1,\alpha_2)\Bigr)
=\]
\[
\pm\Bigl(
6\cdot 4^{3}\;+
9(-12)^{2}\cdot 4 
\Bigr)
=
\pm 5568,
\]
\[
\int_X \alpha_1^{3}\alpha_2\,c_{2}(X)
=\frac{108\cdot 4!}{4!}B(\alpha_1,\alpha_1) B(\alpha_1,\alpha_2)=108\cdot (-12)\cdot 4=
\mp\, 5184= \int_X \alpha_1\alpha_2^{3}\,c_{2}(X) , 
\]
\[
\int_X \alpha_1\alpha_2\,c_{2}(X)^{2}
= 1200\; B(\alpha_1,\alpha_2)= 1200\cdot \pm 4=
\pm 4800.
\]

Summing the four contributions gives:
\[
\int_X \alpha_1^{3}\alpha_2^{3}
+
\int_X \alpha_1^{3}\alpha_2\,c_{2}(X)
+
\int_X \alpha_1\alpha_2^{3}\,c_{2}(X)
+
\int_X \alpha_1\alpha_2\,c_{2}(X)^{2}
=0.
\]
Thus $[P_1]\cdot[P_2]=0$ in all cases.

Hence $P_1 \cap P_2$ cannot be finite and non-empty.
If it were positive dimensional, it would contain a curve whose class is proportional
to both $\alpha_1$ and $\alpha_2$.
Since these are distinct ppMBM classes of square $-12$, this is impossible.
Therefore $P_1 \cap P_2=\varnothing$.
\qed

\

\noindent
The next proposition is devoted to deriving geometric properties of  a K3$^{[3]}$ type manifold
whose Picard lattice is isomorphic to the lattice $BW$. More precisely, we provide a geometric interpretation of the $256$ ppMBM classes of $X$ whose existence was established in Proposition~\ref{prop:walldiv}.

\subsection{Proposition}\label{256P3I}
Let $X$ be a $K3^{[3]}$ type variety whose Picard lattice is isometric to the Barnes-Wall lattice $BW$.
Then $X$ contains $256$ pairwise disjoint projective $3$--spaces.

Under the identification $j_X: H_2(X,\mathbb{Q}) \rightarrow H^2(X,\mathbb{Q})$
induced by the BBF form,
the doubles of the classes of lines in these  $\mathbb{P}^3$'s yield the  $256$ distinct  ppMBM classes in $\Pic(X)$, of square $-12$ and divisibility $2$.

\ts 
By Proposition \ref{prop:walldiv} there are no SPE classes in $\Pic(X)$, and by Theorem \ref{thm:coni} one has $\mathcal{FE}_X = \mathcal{C}_X$. It follows that
\[
\mathcal{BK}_X
=
\mathcal{C}_X \setminus
\bigcup_{\alpha \in \MBM(X)\cap \Pic(X)} \alpha^{\perp}= \mathcal{C}_X \setminus
\bigcup_{\alpha \in \operatorname{ppMBM(X)}} \alpha^{\perp}
\]
and its connected components are the K\"ahler cones of hyperk\"ahler manifolds bimeromorphic to $X$. In particular, for every $\alpha \in  \operatorname{ppMBM(X)}$, the hyperplane $\alpha^{\perp}$ is a codimension one face of the K\"ahler cone of some hyperk\"ahler manifold bimeromorphic to $X$.

 We first observe that, if  $\alpha^{\perp}$ is a
codimension one face of the K\"ahler cone of  $X$, then $X$
contains a projective $3$--space $\mathbb{P}$ such that twice the class of a line in
$\mathbb{P}$ coincides with $\alpha$ (up to sign), and every curve whose class
is proportional to $\alpha$ is contained in $\mathbb{P}$.

Indeed, by \cite[Corollary~5.9]{BLjems}, since $\alpha^{\perp}$ is an
 extremal face of the K\"ahler cone of  $X$,  there exists a contraction morphism
contracting precisely the curves whose classes are proportional to $\alpha$.
Moreover, by \cite[Proposition~5.8]{BLjems}, this contraction deforms along the
Hodge locus of $\alpha$ in such a way that the resulting family of singular
varieties is locally trivial.
By \cite[Corollary~5.7]{Amerik_V}, the contraction at a generic point of the Hodge
locus contracts exactly a projective space of dimension $3$ whose lines have class $\alpha/2$. 
By local triviality, the same description must hold for the contraction on $X$.

By the same argument, if $\alpha \in \operatorname{ppMBM(X)}$ and $\alpha^{\perp}$ is a common
codimension one face of the K\"ahler cones of two hyperk\"ahler manifolds $X_1$ and
$X_2$ that are bimeromorphic to $X$, then both $X_1$ and $X_2$ admit a contraction of
a $3$--dimensional projective space. Moreover, the curves contracted by these morphisms are
precisely the curves on $X_1$ or $X_2$, respectively, whose classes are
proportional to $\alpha$. 
This shows that crossing the wall $\alpha^\perp$ defined by $\alpha$ in order to pass from
$X_1$ to $X_2$ is equivalent to performing the  Mukai flop of a $3$--dimensional  projective space
$\mathbb{P}_1 \subset X_1$.

Since any hyperk\"ahler manifold $X'$ bimeromorphic to $X$ can be reached from any other by crossing finitely
 many codimension one faces of the K\"ahler cones, and since projective $3$--spaces
 corresponding to distinct MBM classes on a fixed model are pairwise disjoint by Lemma~\ref{lemmaintvuota}, 
 the existence of a projective $3$--space on $X'$ whose curves have class
 proportional to a given $\alpha \in \operatorname{ppMBM(X)}$ propagates, via a finite sequence of Mukai flops 
 along disjoint projective spaces, to  $X$. In particular, $X$ contains such a projective $3$--space.
 
 Since, as observed above, for every $\alpha \in \operatorname{ppMBM(X)}$ 
 the hyperplane $\alpha^{\perp}$ is a codimension one face of the K\"ahler cone of some hyperk\"ahler
 manifold $X'$ bimeromorphic to $X$, it follows that $X$ contains $256$ projective $3$--spaces whose lines represent all the $256$  ppMBM classes of $X$.
\qed

\

\subsection{Proposition}\label{prop:256P3II}
Let $X$ be a hyperk\"ahler manifold of K3$^{[3]}$ type whose Picard lattice is isometric to the Barnes-Wall lattice $BW$. Then there exists a normal K\"ahler variety $\widehat{X}$ with symplectic singularities and no complete curves such that the following holds.

For every hyperk\"ahler manifold $X'$ bimeromorphic to $X$, there exists a contraction
\[
c_{X'} \colon\, X' \,\longrightarrow\, \widehat{X}
\]
contracting precisely the $256$ projective $3$--spaces of Proposition~\ref{256P3I}. In particular $X'$ contains exactly $256$ copies of $\mathbb{P}^3$

\ts
As in the proof of Proposition~\ref{256P3I}, the connected components of
\[
\mathcal{C}_X \setminus \bigcup_{\alpha \in \operatorname{ppMBM}(X)}\alpha^{\perp}
\]
coincide with the K\"ahler cones of the hyperk\"ahler manifolds bimeromorphic to $X$.

Since the  ppMBM classes in $\mathrm{Pic}(X)$ generate the $\mathbb{Q}$ vector space $\mathrm{Pic}(X)\otimes_{\mathbb{Z}}\mathbb{Q}$ of dimension $16$,  the intersection
\[
\bigcap_{\alpha \in \operatorname{ppMBM(X)}} \alpha^{\perp}
\]
is a linear subspace of codimension $16$ which meets the positive cone $\mathcal{C}_X$. It follows that this subspace defines an extremal face of codimension $16$ of the K\"ahler cone of some hyperk\"ahler manifold bimeromorphic to $X$. Without loss of generality, we may assume that this manifold is $X$ itself.

By extremality, \cite[Corollary 5.9]{BLjems} 
(see also \cite[Theorem 4.6]{AVsel}) yields 
the  existence of a contraction
\[
c_X:X\longrightarrow \widehat{X}
\]
which contracts all curves whose classes are positive multiple of  ppMBM classes.  
Since, as before, the classes in $\operatorname{ppMBM}(X)$ generate $\Pic(X)$, the class of every curve on $X$ is a linear combination of the classes of these contracted curves. Hence the morphism
$c_X$ contracts all curves on $X$.

Since $\widehat{X}$ admits a symplectic resolution, it is a symplectic variety. By Proposition~\ref{256P3I}, $X$ contains $256$ projective $3$--spaces.  We now show  that the only irreducible complete curves in $X$ are those contained in these subvarieties; this ensures that $c_X$ contracts precisely them and that $\widehat{X}$ has no complete curves.

Since $\widehat{X}$ is a symplectic variety, it has canonical singularities. The rational connectedness conjecture of Shokurov, proved by Fujino in the complex analytic setting (see \cite{Fu}), then implies that the fibres of the contraction $c_X$ are rationally connected.

As observed above, every curve on $X$ is contracted by $c_X$. In particular, the existence of an irreducible complete curve not contained in any of the $256$ projective $3$--spaces would yield an irreducible rational curve on $X$ outside these subvarieties. 

By a standard argument, this produces an irreducible rational curve  of minimal degree (with respect to a K\"ahler class) whose deformations cover the same locus as the deformations of the original curve. 

Therefore, it suffices to show that every minimal rational curve on $X$ is contained in one of the $256$ projective $3$--spaces.

Let $C \subset X$ be such a curve. By \cite[Proposition 2.4]{Amerik_V}
(see also \cite[Definition~2.5]{Amerik_V}), the cohomology class of $C$
in $H^2(X,\mathbb{Z})$ is a primitive  MBM class and, since $C$ is effective, it is also a ppMBM class. Hence it is proportional to one of
the $256$ primitive ppMBM classes $\alpha_0 \in \Pic(X)$.

We now show that $C$ is contained in the $3$--dimensional projective space
$\mathbb{P}_{\alpha_0}$ corresponding to $\alpha_0$, that is, the one such that
twice the class of a line in $\mathbb{P}_{\alpha_0}$ equals $\alpha_0$. 

Since, by \cite[Theorem 1]{Harvey_HT}, the cohomology class of
$\mathbb{P}_{\alpha_0}$ is given by
\[
\frac{\alpha_0^{3} + c_2(X)\,\alpha_0}{48},
\]
it follows from \cite[Th\'eor\`eme 0.1]{Vlagr} that the pair
$(\mathbb{P}_{\alpha_0}, X)$ deforms over a Euclidean open subset $U$
of the Hodge locus of $\alpha_0$.

Since $C$ is a minimal rational curve and its class is proportional to
$\alpha_0$, the pair $(C, X)$ also deforms over $U$ by \cite[Corollary 4.8]{AVimrn}.
For every $t \in U$ such that the corresponding deformation $X_t$ of $X$
has cyclic Picard group, every curve on $X_t$ is contained in the unique
projective three-space in $X_t$. Therefore, the deformation of $C$ over $t$
is contained in the deformation of $\mathbb{P}_{\alpha_0}$ over $t$.

This shows that $C$ is a specialization of a family of curves contained in
the family of $\mathbb{P}^3$'s arising from the deformation of
$\mathbb{P}_{\alpha_0}$ over $U$. Since this family is proper, it follows
that $C$ is contained in $\mathbb{P}_{\alpha_0}$.

Finally, let $X'$ be a hyperk\"ahler manifold bimeromorphic to $X$. 
Since $\widehat{X}$ contains no rational curves, there exists a contraction
\[
c_{X'} \colon \, X'\, \longrightarrow\, \widehat{X}.
\]
It remains to determine its exceptional locus.

As $\widehat{X}$ has $256$ isolated singularities, their preimages in the symplectic resolution $X'$ are $256$ projective $3$--spaces (see \cite[Lemma 2.1]{WW}). Hence the exceptional locus of $c_{X'}$ consists precisely of these subvarieties, and the proof is complete.
\qed

\

\noindent
As a consequence of Proposition \ref{prop:256P3II},
bimeromorphic automorphisms of $\widehat{X}$ are biregular.

\subsection{Corollary}\label{cor:no-curves-Y}
Every meromorphic map from a normal compact complex space to $\widehat{X}$ is regular.  
In particular, every bimeromorphic map from $\widehat{X}$ to $\widehat{X}$ is an isomorphism.

\ts Given a meromorphic map 
\[
f : Y \dashrightarrow \widehat{X},
\]
by resolution of singularities (see, for instance, \cite{Wl}) there exists a
smooth complex space $\widetilde{Y}$ obtained from $Y$ by a finite sequence of
blow-ups along smooth centers, together with a lifting
\[
\widetilde{f} : \widetilde{Y} \longrightarrow \widehat{X}
\]
of $f$.  
In particular, the morphism
\[
\pi : \widetilde{Y} \longrightarrow Y
\]
has connected projective fibres, which are therefore covered by curves.

Since $\widehat{X}$ contains no curves, the fibres of $\pi$ are contracted by
$\widetilde{f}$.  
This means that the closure $\Gamma$ of the graph of $f$ in
$Y \times \widehat{X}$
is finite over $Y$.  

Since it is generically injective, by the normality of $Y$ the projection
$\Gamma \longrightarrow Y$
is an isomorphism.  
Therefore $f$ is regular.

\qed

\

\noindent
As a further consequence, we can identify the singularities of $\widehat{X}$ as
those arising from point contractions of Lagrangian projective spaces on hyperk\"ahler
manifolds of dimension six.

\subsection{Corollary}\label{cor:gsing}
The variety $\widehat{X}$ has precisely $256$ singular points.  
At each of these points, the analytic germ of $\widehat{X}$ is isomorphic to the cone 
over the incidence variety
\[
I := \{([x],[L]) \in \mathbb{P}^3\times \mathbb{P}^{3\vee} \mid L(x)=0\}
\subset \mathbb{P}^3 \times \mathbb{P}^{3\vee} \subset \mathbb{P}^{15}.
\]

\ts It is well known (see e.g.~\cite{Be}) that the singularity arising from
the contraction of a projective three-space in a six dimensional hyperk\"ahler
manifold is analytically isomorphic to the singularity at the origin of the
closure of the minimal nilpotent orbit in the Lie algebra $\mathfrak{sl}(4)$.

This minimal nilpotent cone consists of trace-free nilpotent matrices of rank~$1$,
and its closure is isomorphic to the affine cone over the incidence variety $I\subset \mathbb{P}^{15}$. 
\qed

\

\noindent
Thanks to the global Torelli theorem for K3$^{[3]}$ varieties, the previous 
proposition and corollaries allow us to deduce the existence of a bimeromorphic  
symplectic involution on $X$.  
Moreover, using the structure of the contraction $c_X : X \rightarrow \widehat{X}$ and the
description of the singularities of $\widehat{X}$, we can describe  the 
indeterminacy locus of this involution explicitly.  

\subsection{Corollary}\label{cor:indet}
Let $X$ be a K3$^{[3]}$ variety whose Picard lattice is isometric to $BW$, 
and let 
\[
c_X : X \longrightarrow \widehat{X}
\]
be the contraction constructed in Proposition~\ref{prop:256P3II}. Then:
\begin{enumerate}[(a)]
    \item $X$ admits a bimeromorphic symplectic involution 
    \(
        \iota_X : X \dashrightarrow X
    \),
    and $BW$ is the anti-invariant lattice of $\iota_X^*$.

    \item The bimeromorphic involution $\iota_X$ induces a regular involution 
    \(
        \iota_{\widehat{X}} : \widehat{X} \,\longrightarrow\, \widehat{X}
    \)
    which fixes each of the $256$ singular points of $\widehat{X}$.

    \item The involution $\iota_X$ is regular on the complement of the
$256$ projective three-spaces in $X$, and it does not extend to a
regular morphism along any of these $256$ $\mathbb{P}^3$'s.
   
\end{enumerate}

\ts 
(a) We recall, from subsection Section~\ref{sub:invret}, 
that there exists an involution on the lattice 
$H^{2}(X,\mathbb{Z})$ equipped with the BBF form,
whose anti-invariant lattice is precisely $\Pic(X)=BW$.  

The involution
above, being the identity on the orthogonal
complement of the negative definite lattice $BW$, preserves the orientation
of the positive cone.  By the global Torelli theorem for 
K3$^{[3]}$ varieties (see \cite[Corollary 5.7]{Msurv}), since the monodromy 
group of a K3$^{[3]}$ is maximal, such a lattice automorphism is induced by a 
bimeromorphic  automorphism of $X$ provided it preserves the closure of the birational K\"ahler 
cone.  
This condition is satisfied because, by Proposition \ref{prop:walldiv}, the set $\SPE(X)$ of stably prime exceptional classes on $X$ is empty, and therefore, by Theorem \ref{thm:coni}, 
its birational K\"ahler cone is an open dense subset of  the whole positive cone.

(b) We have already seen in Corollary~\ref{cor:no-curves-Y} that the bimeromorphic
involution $\iota_X$ becomes regular on $\widehat{X}$.
By construction, $\iota_X$ acts as $-1$ on  $\Pic(X)$, and therefore it
preserves the orthogonal hyperplane of each   ppMBM class in $\Pic(X)$.
If $\iota_{\widehat{X}}$ were to map one singular point of $\widehat{X}$ to a
different one, then the induced action on the lattice would have to send the
orthogonal hyperplane of one ppMBM class of $X$  to that of a different, hence  non proportional,  ppMBM class.

(c) Observe that the contraction \(c_X\) induces an isomorphism between the 
complement of the \(256\) projective three–spaces in \(X\) and the smooth 
locus of $\widehat{X}$.  
Since $\iota_{\widehat{X}}$ is an isomorphism, it follows that \(\iota_X\) restricts to an 
automorphism of the complement of these \(\mathbb{P}^3\)'s.  
Finally, $\iota_X$ cannot be extended to  a regular morphism  on any of the $256$ projective 
three-spaces.  
Indeed, if $\iota_X$ were regular  at every  point of such a $\mathbb{P}^3$, then the 
class of a line in that $\mathbb{P}^3$ would have to be mapped to the class of 
an effective curve.  
However, by construction $\iota_X$ acts as $-1$ on the corresponding MBM class, 
so the line class is sent to its negative, which cannot be effective.  
Thus $\iota_X$ has no regular extension  on any of the $256$ exceptional 
$\mathbb{P}^3$'s.
\qed

\


\section{Bimeromorphic involutions on Hilbert schemes of Kummer surfaces}\label{s:BimeroKummer}

\subsection{A birational involution on the  Hilbert scheme $S^{[3]}$ of a Kummer surface}
\label{sub:birinvs3}
The lattice-theoretic tools employed in the previous section allowed us to
extract strong birational information about the involution $\iota_X$ on a
hyperk\"ahler manifold $X$ of K3$^{[3]}$ type with $\Pic(X)\cong BW$.  
These methods are not sufficient to determine its fixed locus.
In order to show that a Kummer fourfold appears 
as a component of maximal dimensional of the fixed locus, 
we now specialize $X$ and turn to a geometric
construction of the involution based on the geometry of the Jacobian of a genus two curve.

\ 
  
We begin by fixing notation.
Let $A$ be the Jacobian of a very general genus two curve $C_0$, let
$S' =K^s(A)= A / \{\pm 1\}$ be the 
singular Kummer surface and let
$S \rightarrow S'$ denote the minimal resolution of $S'$, i.e.\ $S=K(A)$ is the Kummer K3
surface associated with the abelian surface $A$.

Recall that $S'$ admits a natural embedding as a quartic surface in
$\mathbb{P}^{3}$.
We shall say that three points of $S$ are \emph{collinear} if the images of
their projections in $S' \subset \mathbb{P}^3$ lie on a line.

Let $R$ denote the line bundle on $S$ whose square admits a section vanishing
precisely along the $16$ exceptional curves of $S$. 
The line bundle $R$ defines the
ramified double cover $\pi: \widetilde{A}\rightarrow S$, where $\widetilde{A}$
is the blow-up of $A$ along its $2$--torsion points, and $\pi$ ramifies exactly
along the union of the nodal curves of $S$, so that
\[
\pi_{*}\mathcal{O}_{\widetilde{A}} \simeq \mathcal{O}_{S} \oplus R^{\vee}.
\]
The restriction of $R$ to a smooth curve $C\subset S$ that doesn't intersect the 16 exceptional curves defines an unramified  double cover $\widetilde{C}\rightarrow C$ where $\widetilde{C}$ is the inverse image of $C$ in $\widetilde{A}$. Thus $R|_C$ is a line bundle of order two (and degree zero) in
$\Pic(C)$.

The birational symplectic involution we shall study acts on $S^{[3]}$.
A general point of $S^{[3]}$ is represented by a length three subscheme
consisting of three non-collinear points of $S$ (in the sense defined above).

For three non-collinear points $p_{1},p_{2},p_{3} \in S$, let
$\Pi \subset \mathbb{P}^{3}$ be the unique plane whose intersection with
$S' \subset \mathbb{P}^{3}$ contains the images of
$p_{1},p_{2},p_{3}$.
We denote by
\[
C \subset S
\]
the schematic preimage in $S$ of the plane section $\Pi \cap S'$ under the
resolution map $S \rightarrow S'$.

Let $I \subset \mathcal{O}_{C}$ denote the ideal sheaf of the subscheme
$\{p_{1},p_{2},p_{3}\}$ inside $C$.
We consider the twisted ideal sheaf
\[
I \otimes R|_{C}.
\]
For three general points on $S$, this sheaf is isomorphic to the ideal sheaf of
a unique Artinian subscheme of length three on $C \subset S$.
In particular, it defines a length three subscheme of $S$, that is, an element
of $S^{[3]}$.

The birational involution
\[
\iota_{S^{[3]}} : S^{[3]} \dashrightarrow S^{[3]}
\]
is defined by sending a general subscheme
\(
\xi = \{p_{1},p_{2},p_{3}\}
\)
to the length three Artinian subscheme whose ideal sheaf on $C$ is isomorphic to 
$I \otimes R|_{C}$.

We observe that $\iota_{S^{[3]}}$ is indeed an involution:

\[
(I \otimes R|_{C}) \otimes R|_{C}
\simeq I \otimes (R|_{C})^{\otimes 2}
\simeq I,
\]
and therefore $\iota_{S^{[3]}}^{2} = \mathrm{id}$ on the open subset where $\iota_{S^{[3]}}$ is defined.

\subsection{The Picard lattice of $S^{[3]}$}
We recall the structure of $H^{2}(S^{[3]},\mathbb{Z})$.
There is an identification
\[
H^{2}(S^{[3]},\mathbb{Z}) \;\cong\; H^{2}(S,\mathbb{Z}) \oplus \mathbb{Z}\delta,
\qquad (\delta,\delta)=-4.
\]
As in Section~\ref{sec:Kumlat}, we denote by
\[
N_x \in \Pic(S)\qquad (x\in A[2]),
\]
the classes of the exceptional curves of the resolution
$S \rightarrow S'=A/\{\pm1\}$ and, by abuse of notation, their images in
$H^{2}(S,\mathbb{Z})\subset \mathrm{Pic}(S^{[3]})$. 
The class $N_x\in \Pic(S^{[3]})$ is represented by a divisor
$E_{x}$, which parametrizes length--$3$ subschemes whose image in $S'$
contains the singular point defined by $x$, again denoted by
$x \in S'=K^s(A) \subset \mathbb{P}^{3}$.

We also denote by
\[
L \in \mathrm{Pic}(S)
\]
 the pullback of the class of a plane section of $S'$, and we use the same
notation for its image in $\mathrm{Pic}(S^{[3]})$.
Thus $\mathrm{Pic}(S^{[3]})$ is the saturation of the subgroup of rank $18$
\[
\langle L,\delta\rangle \,+\, \langle N_x:\;{x\in A[2]}\,\rangle
\]
inside $H^{2}(S^{[3]},\mathbb{Z})$.

\subsection{The divisors $D_p$ for $p\in\PP^3$}\label{sec:Dp}
To study the involution $\iota_{S^{[3]}}$, we introduce an invariant isotropic divisor class.
To every point $p \in \mathbb{P}^{3}$ we associate a divisor $D_{p}\subset S^{[3]}$ parametrizing length--$3$ subschemes $Z \subset S$ for which there exists a plane through $p$ in $\mathbb{P}^{3}$ whose inverse image under the morphism $S \rightarrow S' \subset \mathbb{P}^{3}$ contains $Z$.

By definition of $\iota_{S^{[3]}}$, the divisor $D_{p}$ is invariant, and hence
its cohomology class $[D_{p}]$ is invariant as well.
Moreover, the class $[D_{p}]$ does not depend on the choice of $p$.
Since no plane in $\mathbb{P}^{3}$ contains four general points
$p_{1}, p_{2}, p_{3}, p_{4}$, we have an empty set-theoretic intersection
\[
\bigcap_{i=1}^{4} D_{p_{i}}\,=\,\emptyset~.
\]
It follows that $[D_{p}]^{6}=0$, and by the Fujiki relation we have
\[
B([D_{p}], [D_{p}]) = 0~,
\]
that is, $[D_{p}]$ is isotropic with respect to the BBF form.
The following lemma determines the class of $D_p$ in $\Pic(S^{[3]})$.

\subsection{Lemma} The isotropic divisor class $[D_p]\in\Pic(S^{[3]})$ is
$$
[D_p]\, =\, L - \delta~.
$$

\ts
In order to determine the coefficients $a$, $b$, and $c_x$ in the expression
\begin{equation}\label{eq:Dp}
	[D_{p}] = aL + b\delta + \sum_{x\in A[2]} c_x N_x
\end{equation}
as a $\mathbb{Q}$--linear combination of the generators of the Picard group of $S^{[3]}$,
we make the following observation.

Recall that for $x\in A[2]$, the class $N_x$ is represented by a divisor
$E_{x}$, which parametrizes length--$3$ subschemes whose image in $S'$
contains the singular point defined by $x$ in $K^s(A)$.
Since $E_{x} \subset D_{x}$, an argument analogous to the one used above
shows that, if $p_1, p_2, p_3$ and $x$ are not coplanar in $\mathbb{P}^{3}$,
then there is an empty the intersection
\[
E_{x} \cap D_{p_1} \cap D_{p_2} \cap D_{p_3}\,=\,\emptyset~,
\qquad \mbox{hence}\quad B([D_p], E_x) \,=\, 0~,
\]
and consequently all coefficients $c_x$ vanish.

Since moreover $B([D_p],[D_p]) = 0$, we deduce that $b = \pm a$.
To determine $a$ and $b$, we evaluate both sides of \eqref{eq:Dp}
on the class $F$ of a one dimensional fiber of the Hilbert-Chow morphism.

For general points $p_1,p_2 \in S$, there exists a unique plane in
$\mathbb{P}^{3}$ containing $p$ whose inverse image in $S$ contains a subscheme
lying in the fiber $F_{2p_1 + p_2}$ of the Hilbert-Chow morphism over the cycle
$2p_1 + p_2$. Moreover, such a subscheme is unique, and therefore
\[
[D_p] \cdot F\, =\, 1~.
\]

On the other hand, since $L$ and the classes $E_i$ are pullbacks of classes
from the symmetric product $S^{(3)}$, we have
\[
L \cdot F = E_x \cdot F \,=\, 0~.
\]
Finally, $2\delta$ is represented by the exceptional divisor of the
Hilbert-Chow morphism, whose normal bundle has degree $-2$ on
$F_{2p_1+p_2}$; hence
\[
\delta \cdot F \,=\, -1~.
\]

It follows that $b = -1$ and $a = \pm 1$. Since $D_p$, $L$, and $\delta$ are
effective classes, we conclude that $a = 1$, and therefore
$[D_p]\, =\, L - \delta$. 
\qed

\

\noindent 
The following proposition describes  the action of $\iota_{S^{[3]}}$ on $H^{2}(S^{[3]},\mathbb{Z})$.  

\subsection{Proposition} \label{prop:antiinv-lattice}
Let $S$ be a very general Kummer surface and let $\iota_{S^{[3]}} : S^{[3]} \dashrightarrow S^{[3]}$ be 
the birational symplectic involution constructed in Section \ref{sub:birinvs3}. Then:
\begin{enumerate}[(a)]
    \item $\iota_{S^{[3]}}$ acts trivially on the transcendental lattice 
    $T\!\left(S^{[3]}\right)$.

    \item The anti-invariant sublattice of $\iota_{S^{[3]}}^{*}$ inside 
    $\Pic(S^{[3]})$ is the \emph{saturation} of the lattice generated by the 
    sixteen classes
    \[
        M_x\,:=\,L - \delta - 2N_x, \qquad x \in A[2]~,
    \]
    where $L$ is the pullback to $S$ of the class of plane section on the singular Kummer surface $S'$, $\delta$ is half the exceptional divisor of the Hilbert-Chow
    morphism, and the $N_x$, $x\in A[2]$,  are the exceptional curves of the Kummer
    surface $S$.
\end{enumerate}

\ts
(a)$\;$  
By Mukai’s construction, the holomorphic symplectic form on a moduli space of 
sheaves on a K3 surface is obtained from the Yoneda pairing
\[
\Ext^{1}(F,F)\otimes \Ext^{1}(F,F)\;\longrightarrow\;\Ext^{2}(F,F)
\]
followed by the trace morphism 
$\Ext^{2}(F,F)\rightarrow H^{2}(\mathcal{O}_{S})\simeq\mathbb{C}$.
Tensoring with a line bundle commutes with both operations.  
Since $\iota_{S^{[3]}}$ is defined by tensoring ideal sheaves of length--$3$ subschemes 
with the $2$--torsion line bundle $R$, the induced map on tangent spaces 
preserves the Mukai symplectic form.  
Thus $\iota_{S^{[3]}}$ is symplectic and acts trivially on $T(S^{[3]})$.

\

(b)$\;$
We can now determine the anti-invariant classes.
Let $x \in S'$ be one of the $16$ nodes, we already observed that $E_x\subset N_x$.
The divisor $D_{x}$ has exactly two reduced irreducible components:
\[
D_{x} \,=\, E_{x} \cup T_x~,
\]
where  $T_x$ parametrizes subschemes lying in the strict transforms of plane sections of $S'$ passing through $x$.

Since the line bundle $R$ has degree $1$ on the strict transform of a plane
section of $S'$ passing through $x$, and degree $-1$ on the nodal curve lying
over $x$, it follows from the definition of the involution
$\iota_{S^{[3]}}$ that
\[
\iota_{S^{[3]}}(T_x) \,=\, E_{x}~.
\]
As $\iota_{S^{[3]}}$ is an involution, we therefore deduce that
\[
\iota_{S^{[3]}}\bigl(E_{x}\bigr)\, =\, T_x~.
\] 
Hence $\iota_{S^{[3]}}$ exchanges the components $E_{x}$ and $T_x$ of $D_{x}$
and
\[
    \iota^{*}_{S^{[3]}}(T_{x} - E_{x})\,=\, -(T_{x} - E_{x})~.
\]
Since $T_x = D_{x} - E_{x}$ we have $[T_x- E_{x}]=L-\delta- 2N_x$.
We thus obtain sixteen linearly independent
$\iota_{S^{[3]}}$--anti-invariant classes, namely
\[M_x\,:=\,L - \delta - 2N_x,
\qquad x\,\in\,A[2]~.
\]

Since $\Pic(S^{[3]})$ has rank $18$ and there exist at least two linearly
independent $\iota_{S^{[3]}}$--invariant classes (namely $L-\delta$ and
$H+\iota_{S^{[3]}}^{*}H$ for any ample class $H$), the
$\iota_{S^{[3]}}$--invariant sublattice has rank at least $2$.
Consequently, the $\iota_{S^{[3]}}$--anti-invariant sublattice has rank at most
$16$.

Since we have exhibited $16$ linearly independent
$\iota_{S^{[3]}}$--anti-invariant classes, it follows that the
$\iota_{S^{[3]}}$--anti-invariant lattice has rank exactly $16$.
Hence, the $\iota_{S^{[3]}}$--anti-invariant lattice is the \emph{saturation}
in $H^{2}(S^{[3]},\mathbb{Z})$ of the lattice generated by the classes
\[
L - \delta - 2N_x, \qquad x\,\in\,A[2]~.  
\]
\qed

\

\subsection{Corollary}
The anti-invariant lattice of the rational involution
$\iota_{S^{[3]}} \colon S^{[3]} \dashrightarrow S^{[3]}$
is isomorphic to the Barnes-Wall lattice ${BW}$.

\ts
The sixteen classes $M_x=L - \delta - 2N_x$ in the anti-invariant lattice have the following
BBF-intersection products:
$$
B(M_x,M_x)\,=\,-8,\qquad B(M_x,M_y)\,=\,0\quad(x\neq y)~.
$$
There is thus an inclusion in the Kummer lattice
$$
\langle M_x:\;x\in A[2]\,\rangle\;\hookrightarrow\; Kum,\qquad
M_x\,\longmapsto\,2N_x~.
$$
The classes $(M_x-M_y)/2=N_x-N_y$ are also in the anti-invariant lattice.
Similarly, the $\mbox{$\frac{1}{4}$}(\sum_{x\in W}\pm M_x)=\half(\sum_{x\in W}\pm N_x)$,
with an even number of minus signs and with $W\subset A[2]$ a subgroup of order eight or its complement, 
are in the anti-invariant lattice. 
Extending the inclusion map $\QQ$--linearly,
we see that $BW\subset Kum$ is isometric to a sublattice of the anti-invariant lattice, 
cf.\ Section \ref{sec:propBW}.

To see that there are no other classes in the anti-invariant lattice one observes that
such a class must be in the dual lattice $BW^*$ of $BW$, that is, it must be a class $c\in BW \otimes\QQ$ such that $B(c,b)\in\ZZ$ for all $b\in BW$. The discriminant group $BW^*/BW\cong (\ZZ/2\ZZ)^8$ and one checks, with a computer, that there are no non-trivial classes with a representative in the anti-invariant lattice.
\qed

\subsection{A birational model of the maximal fixed component}
We now describe a  subvariety of $S^{[3]}$ birational to the Kummer 
fourfold $K(A\times A)$ which is a maximal dimensional component of the 
fixed locus of $\iota_{S^{[3]}}$.

\subsection{Proposition} \label{prop:vfix}
Let $S$ be a very general Kummer surface and let
\[
\iota_{S^{[3]}} \colon\, S^{[3]}\, \dashrightarrow\, S^{[3]}
\]
be the birational symplectic involution constructed in Section \ref{sub:birinvs3}. Then:
\begin{enumerate}[(a)]
    \item There exists a projective subvariety
    \[
        V \;\subset \;S^{[3]}
    \]
    which is birationally equivalent to the Kummer fourfold $K(A \times A)$.

    \item The involution $\iota_{S^{[3]}}$ is well defined on a Zariski open subset
    $U \subset S^{[3]}$ such that $U \cap V$ is dense in $V$ and is pointwise fixed by $\iota_{S^{[3]}}$.
    Moreover $V$ is a maximal irreducible component of the fixed locus of~$\iota_{S^{[3]}}$. 
\end{enumerate}

\ts
(a)$\;$  
For the first assertion it is enough to construct a projective  subvariety 
\[
V' \;\subset\; S'^{(3)}
\]
birational to the Kummer fourfold $K(A\times A)$ and not entirely 
contained in the singular locus of $S'^{(3)}$.  
We then define 
\[
V \;\subset\; S^{[3]}
\]
to be the strict transform of $V'$ under the resolution 
$S^{[3]} \rightarrow S'^{(3)}$.

Recall that $A$ is the Jacobian of a genus two curve $C_0$: the symmetric square $C_{0}^{(2)}$ 
can be realized as the blow-up of
$\Pic^{2}(C_{0})$ at the point corresponding to the canonical bundle $K_{C_{0}}$ of $C_{0}$.
Using the natural isomorphisms
\[
A \,\simeq \,\Pic^{0}(C_{0})\, \simeq \, \Pic^{2}(C_{0})
\]
induced by tensorization with $K_{C_{0}}$, the variety $A \setminus \{0\}$ is
identified with $C_{0}^{(2)} \setminus |K_{C_{0}}|$, while the involution $-1$ on $A$
is induced by the hyperelliptic involution on $C_{0}^{(2)}$.

Let $j \colon C_{0} \rightarrow A$ be a symmetric embedding.
Since the self-intersection of $j(C_{0})$ is $j(C_{0})^2=2$ and the equality
$j(C_{0})+a = j(C_{0})-a$ holds if and only if $a$ is a 2--torsion point,
for $a \in A \setminus A[2]$ the translates $j(C_{0})+a$ and $j(C_{0})-a$
intersect in a symmetric subscheme of length two.

Let $\pi \colon A \rightarrow S'$ be the quotient map.
Define
\[
\phi \colon \,\bigl(C_{0}^{(2)} \setminus |K_{C_{0}}|\bigr)
\times \bigl(A \setminus A[2]\bigr)\,
\longrightarrow\, S'^{(3)}
\]
by
\[
(p+q,a)\;\longmapsto\;
\Bigl\{
\pi\bigl(j(p)+a\bigr),\;
\pi\bigl(j(q)+a\bigr),\;
\pi\bigl((j(C_{0})+a)\cap (j(C_{0})-a)\bigr)
\Bigr\}~.
\]

The map $\phi$ is invariant under the involution acting diagonally as the 
hyperelliptic involution on $C^{(2)}$ and as $-1$ on $A\setminus A[2]$, hence 
descends to a rational map
\[
\Phi:\,K(A\times A)\,\dashrightarrow\, S'^{(3)}~.
\]

We now show that $\phi$ is generically injective.
For $p+q$ and $a$ general, the cycle $\phi(p+q,a)$ consists of three distinct
non-collinear points in the smooth locus of $S'$.
Moreover, by construction, these three points are contained in the singular
plane section $\pi\bigl(j(C_0)+a\bigr)$.

In particular, if
\[
\phi(p_1+q_1,a_1)\,=\,\phi(p_2+q_2,a_2)~,
\]
then necessarily
\[
\pi\bigl(j(C_0)+a_1\bigr)\,=\,\pi\bigl(j(C_0)+a_2\bigr)~,
\]
which implies that $a_1=\pm a_2$.
If $a_1=a_2$, it follows that $p_1+q_1=p_2+q_2$.
Otherwise, the hyperelliptic involution on $C_0^{(2)}$ sends $p_1+q_1$
to $p_2+q_2$.

Finally, for general $(p+q,a)$, the point $\phi(p+q,a)$ is a triple sum of smooth distinct points and 
lies in the smooth locus of $S'^{(3)}$.
Hence
\[
V' \,:=\, \overline{\phi\bigl(C_0^{(2)} \times (A \setminus A[2])\bigr)}
\]
is not contained in the singular locus of $S'^{(3)}$.
Its strict transform $V$ in $S^{[3]}$ is therefore birational to the
Kummer fourfold $K(A \times A)$.

(b)$\;$  A general point of $V'$ corresponds to a triple of distinct, non-collinear
points on a plane section $\pi\bigl(j(C_0)+a\bigr)$ of $S'$ lying entirely in
the smooth locus.
This plane section has exactly one ordinary node, namely
\[
\pi\bigl(\bigl(j(C_0)+a\bigr)\cap\bigl(j(C_0)-a\bigr)\bigr)~,
\]
and this node is one of the three points of the  triple.

Since the plane section lies entirely in the smooth locus of $S'$, its inverse
image $C \subset S$ is isomorphic to it.
We denote by $I$ the ideal sheaf of the corresponding length--$3$ subscheme on
$C$.

In order to check that the birational involution on $S'^{(3)}$ induced by
$\iota_{S^{[3]}}$ is well defined at the generic point of $V'$, and that this
generic point is fixed, it is enough to verify that
\[
I \otimes R|_{C} \,\simeq\, I~,
\]
where $R$ is the square root of the line bundle on $S$ associated with the
exceptional divisor of the resolution $S \rightarrow S'$.

Since the nodal point is one of the three points of the triple, the ideal
sheaf $I$ is a rank--$1$ torsion free sheaf that is not locally free on the
nodal curve $C$.
Hence $I$ is the pushforward of a line bundle $M$ from the normalization
\[
\nu \colon\, C_0\, \longrightarrow\, C~.
\]

By the projection formula, we have
\[
I \otimes R|_{C}
\,=\, \nu_{*}(M) \otimes R|_{C}
\,= \,\nu_{*}\bigl(M \otimes \nu^{*}(R|_{C})\bigr)~.
\]

The double cover $A \rightarrow S'$ sends the reducible curve
$\bigl(j(C_0)+a\bigr)\cup\bigl(j(C_0)-a\bigr)$ onto the image of $C$.
As a consequence, the restriction of $R$ to $C$ induces a double cover with
two irreducible components.
Such a double cover splits on the normalization
$\nu \colon C_0 \rightarrow C$.
Therefore, the pullback $\nu^{*}(R|_{C})$ is trivial, and we obtain
\[
\nu_{*}\bigl(M \otimes \nu^{*}(R|_{C})\bigr)\, =\, \nu_{*}(M) \,=\, I~.
\]
Thus, the generic point of $V$ is fixed by the involution $\iota_{S^{[3]}}$.

Hence $V \subset S^{[3]}$ contains a Zariski open subset on which
$\iota_{S^{[3]}}$ is well defined and acts trivially.
Since $S^{[3]}$ has dimension $6$ and the fixed locus of a non-trivial
symplectic involution has dimension at most $4$, the subvariety $V$
is a maximal irreducible component of the fixed locus of $\iota_{S^{[3]}}$.
\qed

\subsection{Remark} \label{rem:BM-system}
As anticipated in the introduction, another model of a hyperk\"ahler manifold
of K3$^{[3]}$ type, on which the birational symplectic involution admits a
natural geometric description, is the Beauville-Mukai system $M_{|L|,d}$.
It parametrizes rank--$1$ torsion free sheaves of degree $d$, stable with
respect to a generic polarization, on the curves of the complete linear
system $|L|$ on the Kummer surface $S$ associated with the Jacobian $A$.
Here $|L|$ is, as above, the linear system given by the preimages in $S$ of the plane
sections of the quartic surface $S' \subset \mathbb{P}^3$.

The variety $M_{|L|,d}$ carries a natural Lagrangian fibration
\[
\Psi_d :\, M_{|L|,d}\, \longrightarrow\, |L|~,
\]
which associates to any sheaf its Fitting support. If $C \in |L|$ is a smooth
curve, the fiber of $\Psi_d$ over $C$ identifies naturally with $\Pic^d(C)$.
The rational symplectic involution $\iota_{M_{|L|,d}}$ preserves the fibers
of $\Psi_d$ and acts on $\Pic^d(C)$ by tensoring with the restriction to $C$
of the line bundle $R$, which is a square root of the line bundle associated
with the sum of the exceptional curves of $S$.

In this geometrical model, it is easier to describe the subvariety $V_{|L|,d}$ 
birational to
the Kummer fourfold $K(A \times A)$. Let $a \in A$ be a general  point. 
The image of the translate of the symmetric theta divisor  $j(C_0)$ by $a$ in $S'$ does not
intersect the singular locus, and hence can be identified with a curve
$C'_a$ in $S$. Moreover, as $j(C_0)+a\cap j(C_0)-a$ consists of two points $x$ and $-x$,
the curve $C'_a$ has a single node. Let
\[
f_a : \,C_0\, \longrightarrow\, C'_a
\]
denote the corresponding map.

The component $V_{|L|,d}$ of the fixed locus of $\iota_{M_{|L|,d}}$ which is birational to
$K(A \times A)$ is then given by the locus parametrizing sheaves of the form
\[
f_{a*}(M), \qquad M \in \Pic^{d-1}(C_0)~.
\] 
In this case, $V_{|L|,d}$ can also be described as an irreducible
component of maximal dimension parametrizing sheaves which are not locally free on their support.

Finally, we observe that in the case relevant for our purposes, namely $d=-3$, assigning to a generic length--$3$ subscheme of $S$ its ideal sheaf, regarded as a subscheme of the unique curve in the linear system $|L|$ containing it, yields a birational map
\[
S^{[3]}\, \dashrightarrow\, M_{|L|,-3}~.
\]
In $M_{|L|,-3}$, the strict transform of the isotropic divisor $D_p$ on $S^{[3]}$ introduced in 
Section \ref{sec:Dp} is the pullback, via $\psi_{-3}$, 
of the plane $\pi_p\subset |L|$ parametrizing curves passing through the point $p\in \mathbb{P}^3$.

\

\noindent 
In the next proposition we seek a birational model of $S^{[3]}$ on which the 
indeterminacy of the birational symplectic involution $\iota_{S^{[3]}}$ can be resolved 
explicitly.

\subsection{Proposition}\label{prop:contrazione}
Let $S$ be a very general Kummer surface, then
there exists a hyperk\"ahler  manifold $M$ birational to $S^{[3]}$, together with a
contraction 
\[
c_M :\, M\, \longrightarrow \,\widehat{M}~,
\]
which contracts precisely the curves in $M$ whose classes are anti–invariant under 
the birational symplectic involution $\iota_{S^{[3]}}$ on $S^{[3]}$.  

\ts
To prove the existence of a birational model of $S^{[3]}$ with the desired
contraction, it is enough to find a birational model whose ample cone has a face
which is an open subset of the linear subspace orthogonal to the anti-invariant
lattice.
If this holds, a general divisor lying on that face is big and nef, and therefore
defines the required contraction (see also \cite[Corollary 5.9]{BLjems} for the
same result in the K\"ahler  case).

Recall that the Picard lattice of $S^{[3]}$ tensored with $\mathbb{Q}$ is generated by 
$L$, the $N_x$, and $\delta$, and  by Proposition~\ref{prop:antiinv-lattice}(2)
the $\iota_{S^{[3]}}$--anti-invariant $\mathbb{Q}$ vector subspace of
$\Pic(S^{[3]}) \otimes \mathbb{Q}$ is 
generated by the $16$ classes
\[
L- \delta - 2N_x~.
\]
A straightforward computation then shows that the orthogonal complement in $\Pic(S^{[3]}) \otimes \mathbb{Q}$ of the 
anti–invariant subspace (i.e.\ the $\iota_{S^{[3]}}$--invariant subspace) is the plane $\Pi$ generated by
\[\Pi\,:=\,\langle  L - \delta,\,
L - \sum_{x\in A[2]} N_x\rangle \;\subset\; \Pic(S^{[3]}) \otimes \mathbb{Q}~.
\]

Observe that the isotropic  class $L-\delta$, being represented by the divisors $D_p$
introduced in Section \ref{sec:Dp} is movable.
By~\cite[Theorem~1.2]{MZ}, it becomes nef on a suitable hyperk\"ahler  birational model $N$ of
$S^{[3]}$.

Moreover, in the case of varieties of K3$^{[3]}$ type, the SYZ conjecture has
been proved by Markman~\cite[Theorem~1.3]{Ma14} and independently by
Soldatenkov and Verbitsky~\cite[Theorem~3.8]{SoV}. 

For clarity, although we shall not make use of this fact, we note that the variety $N$ and the morphism $\Psi$ may be identified with the moduli space $M_{|L|,-3}$ and the morphism $\Psi_{-3}$ introduced in Remark~\ref{rem:BM-system}.

As a consequence, this birational model $N$ carries a Lagrangian fibration $\Psi: N\rightarrow \mathbb{P}^3$, and the
class $L-\delta$ is the pullback of the  ample generator of ${\rm Pic}(\mathbb{P}^3)$.

The geometry of this Lagrangian fibration allows us to prove the following lemma.

\subsection{Lemma}\label{lemmaEPD}
There exists a Euclidean open subset $U_{\Pi}$ of the plane $\Pi$, contained in
the positive cone $\mathcal{C}_N$ of $N$ such
that every class in $U_{\Pi}$ intersects positively all prime exceptional
divisors of $N$, i.e.\ all irreducible effective divisors of negative square.

\

We now complete the proof of the proposition, assuming the lemma.
We will prove the lemma afterwards.

By \cite[Theorem 6.17]{Msurv}, the fundamental uniruled chamber $\mathcal{FU}_N$ of $N$,
consisting of the classes in the positive cone that intersect positively every
prime exceptional divisor, coincides the fundamental exceptional chamber $\mathcal{FE}_N$, 
i.e.\ with the interior of the closure in the positive cone $\mathcal{C}_N$ of the
birational K\"ahler cone $\mathcal{BK}_N$.
It follows that $U_{\Pi}$ is contained in $\mathcal{FE}_N$.

Since, by Theorem \ref{thm:coni},
\[
\mathcal{BK}_N \,=\, \mathcal{FE}_N \setminus
\bigcup_{\alpha \in \MBM(N)\cap \mathrm{Pic}(N)}\alpha^{\perp}~,
\]
and the hyperplanes $\alpha^{\perp}$, as $\alpha$ varies among the algebraic  MBM classes
(or wall divisors), form a locally finite arrangement inside the positive cone
$\mathcal{C}_N$ (see \cite[Theorem 3.6]{AV21'}), 
up to shrinking the open subset $U_{\Pi}$ we may assume that
all its points lie on the boundary of the K\"ahler cone of a
K3$^{[3]}$ type manifold $M$ (possibly $M=N$ itself) that is birational to
$S^{[3]}$.

This means that a general rational element of $U_{\Pi}$ is big and nef and
therefore defines a contraction
\[
c_M \colon\, M\, \longrightarrow\, \widehat{M}
\]
which contracts precisely the curves whose classes lie in the
$\iota_{S^{[3]}}$--anti-invariant lattice.
This completes the proof of Proposition \ref{prop:contrazione}, assuming Lemma \ref{lemmaEPD}. \qed

\vspace{\baselineskip}\noindent{\bf Proof of Lemma \ref{lemmaEPD}.}$\;\;$
In order to prove the lemma it is useful to divide the prime exceptional divisors on $N$ into two types:

\noindent
(I)  those dominating the base $\mathbb{P}^3$ of the Lagrangian fibration $\Psi:N\rightarrow \mathbb{P}^3$;

\noindent
(II)  those mapped to a divisor on the base  of the Lagrangian 
fibration.

We will need several numerical properties of the prime exceptional divisors on $N$,
which we list and verify below.

\begin{enumerate}
\item Let $D$ be a prime exceptional divisor, and let $P$ be any prime divisor
such that $P \neq D$. Then
\[
B(D,P) \,\geq\, 0.
\]

Indeed, passing to a birational model on which $D$ can be contracted (see
\cite[Proposition~1.4]{Dr}), and denoting by $C$ the class of the curve
contracted by this contraction, there exists a constant $k>0$ such that, for
every prime divisor $P$, one has
\[
B(D, P)\, =\, k\, C \cdot P~.
\]
Since $P \ne D$ and the curves whose class is $C$ sweep out $D$, we have
\[
C \cdot P\, \ge\, 0~,
\]
and the assertion follows.

\

\item Let $D'$ be a prime exceptional divisor of  type~(I). Then
\[
B(D', L - \delta)\, >\, 0~.
\]

Indeed, by a theorem of Matsushita (see \cite[Corollary~1.3]{Mat}),
the Lagrangian fibration on $N$ deforms along the Hodge locus of
$L - \delta$.

By a result of Voisin (see \cite[Proposition 2.1]{Cam}), the fibers of any Lagrangian fibration are projective.
As a consequence, a curve $C'$ given by a complete intersection of ample
divisors on a general fiber of $N$ also deforms along the Hodge locus of
$L - \delta$.

This implies that, via the isomorphism
\[
H_2(N,\mathbb{Q})\, \simeq\, H^2(N,\mathbb{Q})
\]
induced by the BBF form, the class $L - \delta$ is the image
of the class of a positive multiple of $C'$. In particular, there exists $a > 0$ such that
\[
B(D', L - \delta)\, =\, a \,  D'\cdot C' > 0~.
\]

The last inequality holds since $D'$ intersects the general fiber properly
and $C'$ is a complete intersection of ample divisors on that fiber.

\

\item Let $D''$ be a prime exceptional divisor of type~(II). Then
\[
B(D'', L - \delta)\, =\, 0~.
\]

Indeed, since $D''$ is of  type~(II), one already has
\[
D'' \cdot (L - \delta)^3\, =\, 0~,
\]
hence
\[
D''^{3} \cdot (L - \delta)^3 \,=\, 0~.
\]
Since $B(L - \delta, L - \delta) = 0$, by the Fujiki relation it follows that
\[
B(D'', L - \delta) \,=\, 0~.
\]

\

\item Let $D''_1$ and $D''_2$ be prime exceptional divisors of type~(II).
If their images in $\mathbb{P}^3$ are distinct surfaces, then
\[
B(D''_1, D''_2)\, =\, 0~.
\]

Indeed, the image of the intersection $D''_1 \cap D''_2$ in $\mathbb{P}^3$ is a
curve, and therefore the intersection product
\[
[D''_1] \cdot [D''_2] \cdot (L - \delta)^2
\]
vanishes.
Consequently, for any ample class $H$, we have
\[
[D''_1] \cdot [D''_2] \cdot (L - \delta)^2 \cdot H^2 \,=\, 0~.
\]

On the other hand, since $D_{i}''^{3} \cdot (L - \delta)^3 = 0$, we deduce that
\[
B(D''_i, L - \delta) \,=\, 0~,
\qquad
B(L - \delta, L - \delta)\, =\, 0~.
\]
By the Fujiki relation, there exists a non-zero constant $k$ such that
\[
[D''_1] \cdot [D''_2] \cdot (L - \delta)^2 \cdot H^2
\,=\,
k \, B(D''_1, D''_2)\, B(L - \delta, H)^2~.
\]
Since $L - \delta$ is effective, we have $B(L - \delta, H) > 0$, and therefore
\[
B(D''_1, D''_2) \,=\, 0~.
\]
\end{enumerate}

\

We are now in a position to determine all prime exceptional divisors of type~(II) on $N$.

As shown in Proposition~\ref{prop:antiinv-lattice}, if $x \in S'$ is one of the
$16$ singular points, then
\[
D_{x} \,=\, E_{x} + T_{x}~,
\]
where both $E_{x}$ and $T_{x}$ are prime exceptional divisors.
Since $D_{x}$ represents the class $L - \delta$, it follows that
$E_{x}$ and $T_{x}$ are exceptional divisors of type~(II).

We prove that these are precisely the prime exceptional divisors of type~(II) on $N$.
Since
\[
D_{x} \,=\, E_{x} + T_x
\quad\text{and}\quad
[D_{x}]\, =\, L - \delta~,
\]
the divisor $D_{x}$ is the inverse image of a plane $H_x \subset \mathbb{P}^3$,
and no other prime exceptional divisor can have the same plane as its image.

On the other hand, by property~(2), every prime exceptional divisor of type~(II)
belongs to the orthogonal complement $(L - \delta)^{\perp}$ inside the Picard
lattice.  
This orthogonal complement has rank $17$ and is degenerate, since $L - \delta$
lies in its radical.

If $D''$ were a prime exceptional divisor of type~(II) whose image in
$\mathbb{P}^3$ is different from all the planes $H_x$, then, by property~(4),
the class of $D''$ would be orthogonal to all the classes $E_x$, with $x \in A[2]$.

However, the classes $E_x$ already generate a negative definite sublattice of
rank $16$. Therefore, there cannot exist any other class of negative square in
$(L - \delta)^{\perp}$ that is orthogonal to all the $E_x$.

We can now conclude the proof of the Lemma.

Since the class $L-\delta$ is nef on $N$, and the classes $L-\delta$ and
$L-\sum_{x\in A[2]}N_x$ generate the plane $\Pi$, it is enough to show that,
for $\epsilon>0$ sufficiently small, the class
\[
\Delta_{\epsilon}\,:=\,
(L-\delta) + \epsilon\left(L-\sum_{x\in A[2]}N_x\right)
\]
has positive square and intersects positively every prime exceptional
divisor.

A direct computation gives
\[
B(\Delta_{\epsilon},\Delta_{\epsilon})\,=\, 4\epsilon - 28\epsilon^{2}~,
\]
which is positive for $\epsilon$ sufficiently small.

Moreover, we show that for $0 < \epsilon < \frac{1}{16}$, if $D'$ is a prime exceptional divisor of type~(I), 
then
\[
B(\Delta_{\epsilon}, D')\, >\, 0~.
\]

Indeed, in our situation $T_x$ and $D'$ are distinct prime exceptional divisors, and we have
\[
L - \delta \,=\, N_x + [T_x]~.
\]
Hence, by property~(1),
\[
B(L - \delta, D') - B(N_x, D') \,=\, B(T_x, D')\, \ge\, 0~,
\]
and therefore
\[
16\,B(L - \delta, D') - B\!\left(\sum_{x \in A[2]} N_x, D'\right) \,\ge\, 0~.
\]

Moreover, since $L$ is nef on $S^{[3]}$ and the BBF form is a birational invariant, we have
\[
B(L, D') \,\ge \,0~.
\]
It follows that
\[
16\, B(L - \delta, D') + B\!\left(L - \sum_{x \in A[2]} N_x, D'\right) \,\ge\, 0~.
\]

Since $D'$ is a prime exceptional divisor of type~(I), we have
\[
B(L - \delta, D') \,>\, 0~.
\]
Therefore, for any $\alpha > 16$,
\[
\alpha\, B(L - \delta, D') + B\!\left(L - \sum_{x \in A[2]} N_x, D'\right)\, >\, 0,
\]
and hence
\[
B(\Delta_{\epsilon}, D')\, >\, 0
\]
for $0 < \epsilon < \frac{1}{16}$.

Finally, if $D''$ is a prime exceptional divisor of type~(II), then either
$[D'']=N_x$ or $[D'']=[T_x]$ for some $x\in A[2]$.
This concludes the proof of Lemma \ref{lemmaEPD}, since a direct computation shows that
\[
B(\Delta_{\epsilon}, D'')\, =\, 2\epsilon\, >\, 0.
\]
\qed

\

\subsection{} In the following corollaries, we show that $M$ has properties analogous to those
of a hyperk\"ahler manifold of K3$^{[3]}$ type whose Picard lattice is
isometric to the  lattice $BW$.

\subsection{Corollary}\label{256P3M} Let $M$ be as in Proposition \ref{prop:contrazione}.
The morphism $c_M: M \longrightarrow \widehat{M}$ contracts exactly  $256$ projective three–spaces.  The isomorphism
\[
j_M:H_2(M,\mathbb{Q})\longrightarrow H^2(M,\mathbb{Q}),
\]
induced by the BBF form of $M$, sends twice the classes of the lines contained in these projective spaces to the $256$ ppMBM classes of $M$.

\ts
By \cite[Proposition~5.8]{BLjems} the \emph{targets} of such contractions 
deform locally trivially on the Hodge locus of 
the anti-invariant lattice.  
By Proposition~\ref{prop:256P3II}, for a very general point of this Hodge locus, the contracted model has exactly $256$ singular 
points arising from the contraction of the $256$ Lagrangian $\mathbb{P}^3$’s 
corresponding to the $256$  ppMBM classes in the Picard group.  
Local triviality implies that the same holds for $M$.  
Thus $c_M  \colon M \longrightarrow \widehat{M}$ contracts exactly these $256$ projective three–spaces. 
\qed 

\subsection{Corollary} \label{cor:involution-properties} Let $M$ be as in Proposition \ref{prop:contrazione},
let $\iota_{M}: M \dashrightarrow M$ and 
$\iota_{\widehat{M}}: \widehat{M} \dashrightarrow \widehat{M}$ be the birational involutions induced by 
$\iota_{S^{[3]}}$. Then:

\begin{enumerate}
\item $\iota_{\widehat{M}}$ is biregular and fixes each of the $256$ singular points of $\widehat{M}$.

\item $\iota_{M}$ is regular away from the $256$ projective three-spaces in $M$
contracted by $c_M$.
\end{enumerate}

\ts
To prove that $\iota_{\widehat{M}}$ is biregular, we show that there exists an ample line 
bundle $H$ on $\widehat{M}$ such that its pullback $\iota_{\widehat{M}}^{*}H$,
which is defined as a line bundle in codimension one, extends to a genuine line
bundle on $\widehat{M}$, and that this extension is again ample.

The pullback $c_M^{*}H$ of $H$ to $M$ is trivial on all curves contracted by $c_M$,
and therefore it is orthogonal, with respect to the BBF form, to the $256$ primitive MBM
classes generating the anti-invariant lattice.  
Since these classes span the entire anti-invariant sublattice, it follows that
$c_M^{*}H$ lies in the invariant sublattice of $\Pic(M)$.

Consequently, $c_M^{*}H$ is invariant under $\iota_{M}$, and on the open subset where
$\iota_{\widehat{M}}$ is defined we have
\[
\iota_{\widehat{M}}^{*}H \,\cong\, H~.
\]
Thus $\iota_{\widehat{M}}^{*}H$ extends across the indeterminacy locus to a line
bundle on all of $\widehat{M}$, and this extension coincides with $H$.
Since $\iota_{\widehat{M}}^{*}H \cong H$, the birational involution
$\iota_{\widehat{M}}$ preserves an ample line bundle; hence it is biregular.

The second part of item~(1) and item~(2) follow by the same argument as in the proof of
(b) and (c) of Corollary~\ref{cor:indet}.  
Indeed, the ppMBM of $M$
are
anti-invariant with respect to the action induced by $\iota_{M}$ on
$H^{2}(M,\mathbb{Z})$.  
As a consequence, the involution $\iota_{\widehat{M}}$ fixes the singular points
of $\widehat{M}$, while $\iota_{M}$ cannot be regular along the components of
the exceptional locus of the contraction.
\qed

\


\section{A component of the fixed point locus}\label{s:compfix}

\subsection{} We are now in a position to identify the singular Kummer variety $K^s(A\times A)$
of the abelian fourfold  $A\times A$ with a component of the fixed locus of the 
regular involution $\iota_{\widehat{M}} : \widehat{M} \rightarrow \widehat{M}$.
Here $\widehat{M}$ is the contraction of the hyperk\"ahler variety $M$ which is birational to $S^{[3]}$, where $S=K(A)$ is a very general Kummer surface, from Proposition \ref{prop:contrazione}.


\subsection{Proposition}\label{prop:fixhatM}
The fixed locus of the involution $\iota_{\widehat{M}} : \widehat{M} \rightarrow \widehat{M}$ has a connected 
component of dimension $4$ which is isomorphic to the singular Kummer fourfold 
\[
K^s(A\times A)\,=\, (A\times A)/\{\pm 1\}\qquad \subset\, 
\Big(\widehat{M}\Big)^{\iota_{\widehat{M}}}~.
\]

\ts 
\textbf{Step 1.}
As a first step, we show that the fixed locus of the involution
$\iota_{\widehat{M}} \colon \widehat{M} \rightarrow \widehat{M}$
contains a subvariety $\widehat{Z}$ which is birational to the Kummer fourfold
$K(A \times A)$.

By part~(2) of Proposition~\ref{prop:vfix}, the fixed locus of the birational
involution $\iota_{S^{[3]}}$ on $S^{[3]}$ contains a subvariety $V$ which is birational to
the Kummer fourfold $K(A \times A)$.
Since the birational involution
$\iota_M \colon M \dashrightarrow M$ is induced through a birational
morphism $S^{[3]}\dashrightarrow M$, and since $\widehat{M}$ is obtained from $M$ by contracting $256$
projective three-spaces, it is enough to show that the indeterminacy locus of
the birational map between $S^{[3]}$ and $M$ does not contain the fourfold $V$.
Indeed, if this is the case, we may denote by $\widehat{Z}$ the strict transform of $V$ in
$\widehat{M}$; this then gives an irreducible component of the fixed locus of
$\iota_{\widehat{M}} \colon \widehat{M} \rightarrow \widehat{M}$.

To see that $V$ cannot be contained in the indeterminacy locus of the birational
map between $S^{[3]}$ and $M$, observe that, by a result of Kawamata (see
\cite{Ka}), this birational map is a composition of a finite sequence of flops.
If the four dimensional variety $V$ were contained in the indeterminacy locus,
it would be covered by rational curves.

However, the Kummer fourfold $K(A \times A)$ has Kodaira dimension zero, and 
therefore its birational model $V$ cannot be covered by rational curves.
This proves the first step.

\noindent\textbf{Step 2.}
We now show that the normalization $\overline{Z}$ of $\widehat{Z}$ is isomorphic to the singular
Kummer fourfold $K^s(A\times A)$.

As a preliminary observation, since $\widehat{Z}$ is contained in the fixed locus of the 
symplectic involution $\iota_{\widehat{M}}$, the singular locus $Sing(\widehat{Z})$ of $\widehat{Z}$ is contained in the 
singular locus $Sing(\widehat{M})$ of $\widehat{M}$ and therefore consists only of isolated points.  
In particular, the smooth locus of $\widehat{Z}$ is a symplectic manifold.  
Thus $\widehat{Z}$ and its normalization $\overline{Z}$ have singularities in codimension 
four, and consequently they have terminal singularities (see~\cite{Fle}).

From Step~1 we know that $\widehat{Z}$ is birational to $K^s(A\times A)$.  
Since both $\overline{Z}$ and $K^s(A\times A)$ have terminal singularities and 
trivial canonical bundle, they are isomorphic in codimension~$1$.  
Equivalently, there exists a smooth open subset 
\[
U_{\overline{Z}} \subset \overline{Z}
\]
whose complement has codimension at least two, and a smooth open subset
\[
U_{K} \subset K^s(A\times A)
\]
also with complement of codimension at least two, together with an isomorphism
\[
U_{\overline{Z}} \;\xrightarrow{\;\sim\;}\; U_{K}.
\]

The \'etale double cover of $U_{K}$ induced by the branched covering 
\[
A\times A \longrightarrow K^s(A\times A)
\]
therefore yields an \'etale double cover of $U_{\overline{Z}}$, which extends to a 
finite  cover of degree~$2$ 
\[
W \longrightarrow \overline{Z}
\]
branched only along the singular locus of $\overline{Z}$, by purity of the  ramification locus.  
Since $W$ is finite over $\overline{Z}$ and étale over its smooth locus, it also 
has terminal singularities.  
Moreover, by construction, $W$ contains a Zariski open subset isomorphic to a 
Zariski open subset of $A\times A$.  
Thus there exists a birational map 
\[
W \dashrightarrow A\times A.
\]

As $A\times A$ contains no rational curves, this birational map must be 
everywhere regular.  
Since $W$ is normal with trivial canonical bundle, we conclude that the morphism 
\[
W \longrightarrow A\times A
\]
is an isomorphism.  
Therefore
\[
\overline{Z} \;\simeq\; W /\{\pm 1\} \;\simeq\; K^s(A\times A).
\]
\[\]

By the arguments developed so far, in order to complete the proof of the
proposition it remains to show that the subvariety
\[
\widehat{Z} \subset \widehat{M}
\]
is normal.

Since $Sing(\widehat{Z})\subset Sing(\widehat{Z})$, 
it is sufficient to prove that $\widehat{Z}$ is normal at the points lying in  
$Sing(\widehat{Z})\cap Sing(\widehat{Z})$.
To this end, we need to analyze the action of the involution
$\iota_{\widehat{M}}$ in a neighbourhood of such points.

Before proceeding, we fix the local setting.
Recall that if $p$ is a singular point of $\widehat{M}$, then, since
$\widehat{M}$ is obtained by contracting a copy of $\mathbb{P}^3$ inside a
$6$--dimensional irreducible holomorphic symplectic manifold, it follows (as
observed in Corollary~\ref{cor:gsing}) that the analytic germ of the
singularity at $p$ is isomorphic to the affine cone over the incidence variety
\[
I
=
\{([x],[L]) \in \mathbb{P}^3 \times (\mathbb{P}^3)^{\vee} \mid L(x)=0\}
\subset
\mathbb{P}^3 \times (\mathbb{P}^3)^{\vee}
\subset
\mathbb{P}^{15}.
\]

Let $\widetilde{M}$ denote the blow-up of $\widehat{M}$ along its singular
locus.
We then obtain the  commutative diagram:
\[
\begin{tikzcd}[row sep=large, column sep=large]
I 
  \arrow[r, hook] 
  \arrow[d, "b_0"'] 
  \arrow[loop left, "\iota_{I,p}"] 
& 
\widetilde{M} 
  \arrow[d, "b"] 
  \arrow[loop right, "\iota_{\widetilde{M}}"] \\
\mathbb{P}^3 
  \arrow[r, hook] 
  \arrow[d, "c_0"'] 
& 
M 
  \arrow[d, "c_M"]
  \arrow[loop right, dashed, "\iota_M"]  \\
\{p\} 
  \arrow[r, hook] 
& 
\widehat{M}
  \arrow[loop right, "\iota_{\widehat{M}}"]
\end{tikzcd}
\]

where $b$ is the blow-up of $M$ along the $256$ projective three-spaces, and
$c$ is the contraction morphism.
The diagram is Cartesian, so the morphisms in the first column are the
restrictions of those in the second column.
In particular, the morphism
\[
b_0 \colon I\, \longrightarrow \,\mathbb{P}^3
\]
is the natural projection induced by the inclusion
$I \subset \mathbb{P}^3 \times (\mathbb{P}^3)^{\vee}$. Finally $\iota_{\widetilde{M}}$ and $\iota_{I,p}$ are the regular involutions induced by $\iota_{\widehat{M}}$.

\noindent\textbf{Step 3.}
Let $p$ be a singular point of $\widehat{M}$ contained in $\widehat{Z}$, we 
now prove that the involution
\[
\iota_{I,p} \colon \,I \longrightarrow \,I
\]
is given by
\[
([x],[L])\, \longmapsto\, \bigl(\varphi^{-1}([L]),\, \varphi([x])\bigr),
\]
where
\[
\varphi\, =\, \mathbb{P}(f) \colon\, \mathbb{P}^3 \longrightarrow \,(\mathbb{P}^3)^{\vee}
\]
is the projective isomorphism induced by a linear map
\[
f \colon\, \mathbb{C}^4\, \longrightarrow\, (\mathbb{C}^4)^{\vee}.
\]
Moreover the map $f$ is antisymmetric, namely it satisfies
\[
L(x) \,=\, -\, f(x)\bigl(f^{-1}(L)\bigr)
\]
for all $x \in \mathbb{C}^4$ and $L \in (\mathbb{C}^4)^{\vee}$.
Equivalently, with respect to dual bases, the matrix representing $f$
is skew-symmetric.

Observe that the inclusion
\[
I \,\hookrightarrow\, \mathbb{P}^3 \times (\mathbb{P}^3)^{\vee}
\]
induces an isomorphism on Picard groups.
In particular, the Picard group
$
\Pic(I)
$ of $I$
is generated by the restrictions of the line bundles defining the two
projections of $I$ onto $\mathbb{P}^3$ and $(\mathbb{P}^3)^{\vee}$.
It follows that every automorphism of $I$ is induced by an automorphism of
$\mathbb{P}^3 \times (\mathbb{P}^3)^{\vee}$.

Any automorphism of $\mathbb{P}^3 \times (\mathbb{P}^3)^{\vee}$ is either the
product of an automorphism of $\mathbb{P}^3$ with an automorphism of
$(\mathbb{P}^3)^{\vee}$, or it exchanges the two factors.
In our situation, the first possibility cannot occur: indeed, if
$\iota_{I,p}$ were induced by the product of two automorphisms, then it would
descend to an automorphism of $\mathbb{P}^3$, and consequently the involution
$\iota_{\widetilde{M}}$ would descend to an automorphism in a neighbourhood of
the corresponding $\mathbb{P}^3$.
This contradicts the fact that the involution $\iota_M$ is not defined on any
of the $256$ projective three-spaces contracted in the passage to
$\widehat{M}$ (see Corollary~\ref{cor:involution-properties}).

Therefore, there exist isomorphisms
\[
\varphi \colon \mathbb{P}^3 \,\longrightarrow\, (\mathbb{P}^3)^{\vee},
\qquad
\psi \colon (\mathbb{P}^3)^{\vee}\, \longrightarrow\, \mathbb{P}^3
\]
such that, for every $[x]\in \mathbb{P}^3$ and $[L]\in (\mathbb{P}^3)^{\vee}$,
the involution $\iota_{I,p}$ is given by
\[
\iota_{I,p}([x],[L]) \,=\, \bigl(\psi([L]),\, \varphi([x])\bigr)
\]
and, since $\iota_{I,p}\circ \iota_{I,p}=id$, it follows that $\psi = \varphi^{-1}$.

Moreover, since the involution $\iota_{I,p}$ preserves the incidence variety
$I$, letting
\[
f\colon \mathbb{C}^4\, \longrightarrow\, (\mathbb{C}^4)^{\vee}
\]
be a linear isomorphism inducing $\varphi$, we must have
\[
L(x)=0 \quad \text{if and only if} \quad
f(x)\bigl(f^{-1}(L)\bigr)\,=\,0~.
\]
Equivalently, the two bilinear forms
\[
(x,L)\,\longmapsto\, L(x),
\qquad
(x,L)\,\longmapsto\, f(x)\bigl(f^{-1}(L)\bigr)
\]
vanish on the same pairs of vectors.
Therefore, they are proportional.
Since $\iota_{I,p}$ is an involution, the proportionality constant can only be
equal to $1$ or $-1$.

In the first case we say that the isomorphism $f$ is \emph{symmetric}, whereas 
in the second case we say that $f$ is \emph{antisymmetric}.  
In either case, the fixed locus $\Fix(\iota_{I,p})$ of the involution 
$\iota_{I,p}$ consists of the pairs of the form
\[
([x],[f(x)])
\]
which belong to the incidence variety $I$.  

If $f$ is antisymmetric, then $f(x)(x)=0$ for every $x\in \mathbb{P}^{3}$, so 
that $([x],[f(x)])\in \Fix(\iota_{I,p})$ for all $x$, and hence
\[
\Fix(\iota_{I,p})\,\simeq \,\mathbb{P}^{3}.
\]
If $f$ is symmetric, there exists $x\in \mathbb{P}^{3}$ such that 
$f(x)(x)\neq 0$, and therefore the dimension of $\mathrm{Fix}(\iota_I)$ is 
strictly smaller than $3$.

We now show that, under our assumption that $p \in \widehat{Z}$, the linear map $f$ cannot
be symmetric.
Indeed, the subvariety $\widehat{Z}$ has dimension $4$, and its strict transform in
$\widetilde{M}$, which is also of dimension $4$, must intersect the exceptional
divisor over $p$ in a subvariety of dimension at least $3$.
Moreover, this intersection is contained in the fixed locus of the involution
$\iota_{I,p}$.
It follows that the fixed locus of $\iota_{I,p}$ inside the incidence variety
$I$ has dimension at least $3$.

This excludes the possibility that $f$ is symmetric.
Therefore, the linear map $f$ must be antisymmetric.

\noindent\textbf{Step~4.}
Let $p$ be a singular point of $\widehat{M}$ contained in $\widehat{Z}$.
Then the normal cone of $\widehat{Z}$ at $p$ is normal.

We first observe that the fixed locus $\Fix(\iota_{I,p})$
of the involution $\iota_{I,p}$ is projectively normal in
$\mathbb{P}(\mathbb{C}^4 \otimes (\mathbb{C}^4)^{\vee}) \simeq \mathbb{P}^{15}$.
Indeed, $\Fix(\iota_{I,p})$ is isomorphic to the diagonal
\[
\Delta \,:=\, \{([x],[x])\} \;\subset\; \mathbb{P}^3 \times \mathbb{P}^3
\;\subset\; \mathbb{P}^{15}~,
\]
via an isomorphism induced by a linear isomorphism 
$\mathbb{P}(\mathbb{C}^4 \otimes (\mathbb{C}^4)^{\vee})\simeq \mathbb{P}(\mathbb{C}^4 \otimes \mathbb{C}^4)$.
Since $\Delta\simeq \mathbb{P}^3$ is linearly  normal in $\mathbb{P}^{15}$, it is also projectively normal
and the same holds for
$\Fix(\iota_{I,p})$.

As schemes, the fixed locus $\Fix(\iota_{I,p})$ coincides with the exceptional
divisor over $p$ of the blow-up $\widetilde{Z} \rightarrow \widehat{Z}$ along the singular locus
of $\widehat{M}$.
Indeed, since $\widehat{Z}$ is the fixed locus of $\iota_{\widehat{M}}$, the involution
\[
\iota_{\widetilde{M}} \colon\, \widetilde{M} \rightarrow\, \widetilde{M}
\]
acts trivially on the strict transform $\widetilde{Z}$, and in particular on the
scheme $\widetilde{Z} \cap I$, which is the exceptional divisor over $p$ of
$\widetilde{Z} \rightarrow \widehat{Z}$.
Therefore,
\[
\widetilde{Z} \cap I\; \subset\; \Fix(\iota_{I,p})~.
\]
Since the exceptional divisor $\widetilde{Z} \cap I$ has dimension $3$, and
$\Fix(\iota_{I,p})$ is a smooth irreducible variety of dimension $3$, we conclude
that
\[
\widetilde{Z} \cap I \,=\, \Fix(\iota_{I,p})~.
\]

Finally, the normal cone to $\widehat{Z}$ at $p$ is a normal variety, because the
projectivization of this cone is the exceptional divisor over $p$ of
$\widetilde{Z} \rightarrow \widehat{Z}$, namely $\Fix(\iota_{I,p}) \subset \mathbb{P}^{15}$, which
is projectively normal.

\noindent\textbf{Step~5.} 
Let $p$ be a singular point of $\widehat{M}$ contained in $\widehat{Z}$.
Then the singularity of $\widehat{Z}$ at $p$ is normal.

We first observe that the singularity of the normal cone to $\widehat{Z}$ at $p$ is
isomorphic to the quotient singularity $\mathbb{C}^4/\{\pm 1\}$.
Indeed, let $C$ be the affine cone over the diagonal
\[
\Delta\; \subset\; \mathbb{P}^3 \times \mathbb{P}^3 \;\subset \;\mathbb{P}^{15}.
\]
The map
\[
\mathbb{C}^4 \,\longrightarrow\, C,\qquad
(x_0,x_1,x_2,x_3)\longmapsto
(x_i x_j)_{0\le i\le j\le 3}~,
\]
is finite of degree $2$ and factors through the quotient
\[
\mathbb{C}^4 \,\longrightarrow\, \mathbb{C}^4/\{\pm 1\}~.
\]
Since both $\mathbb{C}^4/\{\pm 1\}$ and $C$ are normal varieties, this map induces
an isomorphism between them.

We now observe that the singularity of $\widehat{Z}$ at $p$ is also isomorphic to the
quotient singularity $\mathbb{C}^4/\{\pm 1\}$.
Indeed, by deformation to the normal cone, there exists a flat family over
$\mathbb{A}^1$ whose special fibre over $0$ is the normal cone to $\widehat{Z}$ at $p$,
namely $\mathbb{C}^4/\{\pm 1\}$, and whose restriction over
$\mathbb{A}^1 \setminus \{0\}$ is trivial with fibre $\widehat{Z}$, that is
$\widehat{Z} \times (\mathbb{A}^1 \setminus \{0\})$.
Since the singularity $\mathbb{C}^4/\{\pm 1\}$ is rigid (see~\cite{Sch71}), it
follows that the singularity of $\widehat{Z}$ at $p$ is isomorphic to that of
$\mathbb{C}^4/\{\pm 1\}$ which is normal.

As already observed, since $\widehat{Z}$ is an irreducible component of the fixed locus of
the involution $\iota_{\widehat{M}}$, its singular points are contained among the
$256$ singular points of $\widehat{M}$.
By Step~5, the singularity of $\widehat{Z}$ at each such point is normal.
Therefore, the variety $\widehat{Z}$ is normal, and the normalization morphism
\[
K^{\mathrm{s}}(A \times A) \,\longrightarrow\, \widehat{Z}
\]
is an isomorphism.

Finally, we show that $\widehat{Z}$ is a connected component of the fixed locus of
$\iota_{\widehat{M}}$.

Indeed, since both $\widehat{M}$ and $K^{\mathrm{s}}(A \times A)$ have $256$
singular points, it follows that $\widehat{Z}$ contains all the singular points of
$\widehat{M}$.
As seen in Step~4, for each singular point $p$ of $\widehat{M}$, the fixed locus
of $\iota_{I,p}$ is isomorphic to $\mathbb{P}^3$ and is contained in the strict
transform $\widetilde{Z}$.

Since the fixed locus of the involution $\iota_{\widetilde{M}}$ on the smooth
variety $\widetilde{M}$ is itself smooth, it follows that $\widetilde{Z}$ is the
unique connected component of the fixed locus of $\iota_{\widetilde{M}}$
intersecting the exceptional locus of the morphism
$c_M\circ b:\widetilde{M} \rightarrow \widehat{M}$.
Consequently, $\widehat{Z}$ is the unique irreducible component of the fixed locus of
$\iota_{\widehat{M}}$ intersecting the singular locus of $\widehat{M}$, and
therefore it is a connected component of the fixed locus of
$\iota_{\widehat{M}}$.
\qed

\subsection{Corollary} \label{cor:Z}  
Let $M$ be a hyperk\"ahler variety which is birational to $S^{[3]}$ where $S=K(A)$ is a very
general Kummer surface.
Let $\widehat{Z}$ be the irreducible component of the fixed locus of the involution
\[
\iota_{\widehat{M}} :\, \widehat{M}\, \longrightarrow \,\widehat{M}
\]
as in the Proposition \ref{prop:fixhatM}, and let $\widetilde{M}$ be the blow-up of $\widehat{M}$
along its singular locus.  
Let $Z$ and $\widetilde{Z}$ denote the strict transforms of $\widehat{Z}$ in $M$
and in $\widetilde{M}$, respectively.

Then the varieties $Z$ and $\widetilde{Z}$ are smooth and are both
isomorphic to the Kummer fourfold \[K(A\times A)\] of $A\times A$, i.e.\ the blow-up
of the singular  Kummer fourfold $K^s(A\times A)$ along its singular locus.

\ts
Since the singular locus of $\widehat{Z}$ coincides with the singular locus of $\widehat{M}$, the 
strict transform $\widetilde{Z}$ of $\widehat{Z}$ inside $\widetilde{M}$ is isomorphic to 
the blow-up $\widetilde{K}(A\times A)$ of the Kummer fourfold $K(A\times A)$ at 
its $256$ singular points.

Finally, the variety $Z$ is isomorphic to $\widetilde{Z}$.  
Indeed, the restriction to $\widetilde{Z}$ of the morphism
\[
\widetilde{M} \,\longrightarrow\, M,
\]
which identifies with the blow-up of $M$ along the $256$ projective three-spaces 
contracted by $c$, is injective and has injective differential.  
This concludes the proof.
\qed

\


\section{Deformation theory and the proof of Theorem \ref{thm:WeilK33}}\label{s:defth}

\subsection{} Our proof of Theorem \ref{thm:WeilK33} is based on deformation theory. More precisely, we construct a
locally complete deformation of the pair $(Z,M)$ over a smooth $5$ dimensional
base, where $M$ is the six dimensional hyperk\"ahler manifold introduced in
Proposition~\ref{prop:contrazione} and $Z$ is the smooth four dimensional subvariety, isomorphic to
the Kummer variety of $A\times A$, introduced in Corollary~\ref{cor:Z}. 

Subsequently we show that along these deformations, 
the subvariety $Z$ remains  of a Kummer fourfold.  Finally we observe that if $(Z_t,M_t)$ is 
a deformation of $(Z,M)$ where $M$ is algebraic, then $Z_t$ is of Weil type, with an imaginary quadratic  field $K$ and with trivial discriminant, moreover any such field $K$ occurs. We find complete four dimensional families of abelian varieties in this way and this proves Theorem \ref{thm:WeilK33}.

The following proposition constructs the deformations of $(Z,M)$.

\subsection{Proposition}\label{prop:defo}
Let $D$ be the Kuranishi space of the holomorphic symplectic manifold $M$, and let
\[
D_{BW}\; \subset\; D
\]
be the Hodge locus associated with the anti-invariant lattice of the birational
involution $\iota_M$.  
Then $D_{BW}$ is a smooth complex manifold of dimension $5$, and there exists a
Euclidean open neighbourhood $U\subset D_{BW}$ of the point $0$ corresponding to $M$
such that the smooth pair $(Z,M)$ deforms over $U$.

\ts Since the anti-invariant lattice is non-degenerate of rank $16$,
the Hodge locus $D_{BW}$ is smooth of dimension  $5= h^{1,1}(M) -16$.
 
By 
Corollary \ref{256P3M} the contraction
\[
c_M : M \,\longrightarrow\, \widehat{M} 
\]
contracts precisely the curves whose classes belong to the anti-invariant
lattice.  
By \cite[Theorem~1.1]{BLjems}, an open neighbourhood $U$ of $0$ in the Hodge locus
$D_{BW}$ of the anti-invariant lattice is also the base of the locally trivial
universal deformation $\widehat{\mathcal{M}}$ of $ \widehat{M}$, and the contraction $c_M : M \rightarrow  \widehat{M}$
deforms over $U$, i.e.\ there exists a deformation $\mathcal{M}$ over $U$ and a commutative diagram  

\[
\begin{tikzcd}
M \arrow[r, hook] \arrow[d,"c_M"]
& \mathcal{M}\arrow[d, "c_{\mathcal{M}}"] \\
 \widehat{M} \arrow[r, hook] \arrow[d] 
  & \widehat{\mathcal{M}} \arrow[d, "\varphi'"] \\
\{0\} \arrow[r, hook] 
  & U
\end{tikzcd}.
\] 

To show that the pair $(Z,M)$ deforms over $U$, we first prove that the
birational involution $\iota_M : M \dashrightarrow M$ extends to a
bimeromorphic involution over $U$. We then deduce that the component $Z$
of the fixed locus of $\iota_M$ also deforms over $U$.

Deforming the birational involution $\iota_M$ over $U$ is equivalent to
deforming the regular involution
\[
\iota_{\widehat{M}} :\, \widehat{M} \,\longrightarrow\, \widehat{M}~.
\]
To achieve this, we argue as follows.

By the universal property of the family $\widehat{\mathcal{M}}$, the regular
involution $\iota_{\widehat{M}}$ of $\widehat{M}$ induces, after possibly
shrinking $U$, an involution
\[
\iota_{\widehat{\mathcal{M}}} \colon \,\widehat{\mathcal{M}}\, \longrightarrow\,
\widehat{\mathcal{M}}
\]
which descends to an involution  $\iota_U$ of the base $U$ of $\varphi'$.

If we show that $\iota_U$ is trivial, that is,
it coincides with the identity on $U$, then the involution
$\iota_{\widehat{\mathcal{M}}}$ induces an involution on each fiber and
provides the desired deformation of $\iota_{\widehat{M}}$.

Since $\iota_U$ has finite order, it is sufficient to show that it acts
trivially on the tangent space at the point $0 \in U$ corresponding to
$\widehat{M}$ in order to conclude that it is the identity.

The tangent space $T_{0}U$ to the base of the locally trivial universal
deformation is naturally identified with
\[
T_{0}U\,\cong\,H^{1}(T_{\widehat{M}})~.
\]
By \cite[Proposition~3.6]{BLjems}, this space is canonically isomorphic to
\[
H^{1}\!\left(c_{M*}\Omega^{1}_{M}\right)~,
\]
where $\Omega^{1}_{M}$ denotes the cotangent bundle of $M$.

Moreover, by \cite[Corollary~2.3]{BLjems}, we have an identification
\[
H^{1}\!\left(c_{M*}\Omega^{1}_{M}\right)\, \simeq\, H^{1,1}(\widehat{M})~,
\]
and this space embeds naturally into $H^{1,1}(M)$.
Under these identifications, the action of the involution on $T_{0}U$
coincides with the restriction of the action of $\iota_{M}$ on $H^{1,1}(M)$
to the subspace $H^{1,1}(\widehat{M})$.

Finally, this restriction is trivial because, by construction, the second
cohomology group $H^{2}(\widehat{M})$ is the orthogonal complement of the
anti-invariant part of $H^{2}(M)$.
Therefore, it coincides with the invariant part of $H^{2}(M)$ under the action
of $\iota_{M}$.

In order to prove that also the pair $(Z,M)$ deforms we recall that the blow-up $\widetilde{M}$ of $\widehat{M}$ along its singular locus is
smooth, and since $\widehat{\mathcal{M}}$ is a locally trivial deformation of
$\widehat{M}$, it follows that the blow-up
\[
\widetilde{\mathcal{M}}\, \longrightarrow\, \widehat{\mathcal{M}}
\]
along the singular locus of $\widehat{\mathcal{M}}$ is also smooth.

Moreover, the singular locus of $\widehat{\mathcal{M}}$ is invariant under
$\iota_{\widehat{\mathcal{M}}}$; hence the involution lifts to a regular
involution
\[
\iota_{\widetilde{\mathcal{M}}} \colon \,\widetilde{\mathcal{M}}
\,\longrightarrow\, \widetilde{\mathcal{M}}~.
\]

We thus obtain a family of involutions on a family of smooth varieties. After
possibly shrinking $U$, there exists a connected component
$\widetilde{\mathcal{Z}}$ of the fixed locus of
$\iota_{\widetilde{\mathcal{M}}}$ which intersects the central fibre
$\widetilde{M}$ along $\widetilde{Z}$. Furthermore,
$\widetilde{\mathcal{Z}}$ is smooth over $U$, and the pair
$(\widetilde{\mathcal{Z}}, \widetilde{\mathcal{M}})$ provides a smooth
deformation of the pair $(\widetilde{Z}, \widetilde{M})$ over $U$.

On the other hand, $\widetilde{M}$ is the blow-up of $M$ along the locus
contracted by
\[
c_M \colon\, M\, \longrightarrow\, \widehat{M},
\]
and the family $\widehat{\mathcal{M}}$ is a locally trivial deformation of
$\widehat{M}$. It follows that $\widetilde{\mathcal{M}}$ is the blow-up of
$\mathcal{M}$ along the locus contracted by
\[
c_{\mathcal{M}} \colon\, \mathcal{M}\, \longrightarrow\, \widehat{\mathcal{M}}~.
\]

Therefore, we obtain a morphism over $U$
\[
b_{\widetilde{\mathcal{M}}} \colon \,\widetilde{\mathcal{M}} \,\longrightarrow\, \mathcal{M}~.
\]

We can now define
\[
\mathcal{Z} \,:=\, b_{\widetilde{\mathcal{M}}}(\widetilde{\mathcal{Z}})
\]
and we claim that the pair $(\mathcal{Z}, \mathcal{M})$ provides the desired
deformation of the pair $(Z, M)$ over $U$.  
Since $b(\widetilde{Z}) = Z$, it is enough to show that $\mathcal{Z}$ is smooth maximal dimensional
over $U$.

For this purpose, it suffices to prove that the restriction of $b_{\widetilde{\mathcal{M}}}$ to
$\widetilde{\mathcal{Z}}$ induces an isomorphism onto its image, that is, that
this restriction is injective and has injective differential.  
After possibly shrinking $U$, it is enough to check this property on the
restriction of $b_{\widetilde{M}}$ to the central fibre $\widetilde{Z}$ of
$\widetilde{\mathcal{Z}}$.  
This has already been established in Corollary \ref{cor:Z}.
\qed

\

\noindent The following lemma shows that the deformations of the Kummer variety $Z\subset M$ 
over an open subset of $U$ are Kummer varieties.

\subsection{Lemma}\label{a}
	Let $A$ be a complex torus of dimension $n\ge 4$ and $K(A)$ be the associated Kummer manifold, i.e.\ the blow up of $A/\{\pm 1\}$ at its singular points.
	\begin{enumerate}
		\item{Small deformations of $K(A)$ are Kummer manifolds.}
		\item{Let $f: \mathcal{K}\rightarrow W$ be a smooth algebraic family of smooth projective varieties parametrized by an irreducible variety $W$ and let $U\subset W$ be a non-empty Euclidean open subset of $W$ such that for all $t\in U$ the fiber $K_t:=f^{-1}(t)$ is a Kummer manifold. Then there exists a Zariski open non-empty subset $W^0\subset W$ such that $K_t$ is a Kummer manifold for all $t\in W^0$. }
	\end{enumerate}   

\ts
1) Since complex tori of dimension $n$ vary in a family of dimension $n^2$, the
	Kuranishi space of a Kummer manifold $K:=K(A)$ has dimension at least $n^2$. Hence,
	it suffices to show that its tangent space, namely the first cohomology group
	$H^1(T_K)$ of the tangent bundle $T_K$ of $K$, has dimension at most $n^2$.
	
	Since $K$ is obtained as the quotient of the blow-up $\widetilde{A}$ of a
	complex torus  $A$ along its $2$--torsion points $A[2]$, by the involution
	induced by $-1$ on $A$, we have a short exact sequence of sheaves
	\[
	0 \longrightarrow T' \longrightarrow T_K \longrightarrow
	\bigoplus_{x \in A[2]} \mathcal{O}_{\mathbb{P}^{n-1}_x}(2)
	\longrightarrow 0,
	\]
	where $T'$ is the $(-1)$--invariant part of the pushforward of the tangent
	bundle of $\widetilde{A}$, and $\mathcal{O}_{\mathbb{P}^{n-1}_x}(2)$ arises
	from the normal bundle of the component of the exceptional divisor of the
	desingularization $K \rightarrow A/\{\pm1\}$ corresponding to the $2$--torsion point
	$x$.
	
	The induced long exact sequence
	\[
	H^1(T') \longrightarrow H^1(T_K) \longrightarrow 
	\bigoplus_{x \in A[2]} H^1\bigl(\mathcal{O}_{\mathbb{P}^{n-1}_x}(2)\bigr)
	\]
	allows us to prove the desired inequality. Indeed, $H^1(T')$ is isomorphic to
	the invariant subspace $H^1(T_{\widetilde{A}})^{(-1)}$ of the first
	cohomology group of the tangent bundle $T_{\widetilde{A}}$ of
	$\widetilde{A}$. Moreover,
	\[
	H^1(T_{\widetilde{A}})^{(-1)} \,\simeq\, H^1(T_A)\, \simeq \,\mathbb{C}^{n^2}~,
	\]
	while
	\[
	H^1\bigl(\mathcal{O}_{\mathbb{P}^{n-1}_x}(2)\bigr)\,=\,0~.
	\]
	Therefore, $\dim H^1(T_K) \leq n^2$.

2) For $t \in U$, the Kummer manifold $K_t$ contains exactly $2^{2n}$ pairwise
	disjoint projective spaces. The class of their union $E$ is divisible by $2$
	in the Picard group of $K_t$, and the associated abelian variety can be
	recovered by considering the double covering defined by a square root of the
	class of $E$, and then contracting the inverse image of $E$.
	
	We now show that this construction can be carried out for every $t$ in a
	Zariski open subset of $W$. Since $E$ deforms over $U$, the algebraicity of
	the relative Hilbert scheme implies that it also deforms over $W$. By openness
	of smoothness, there exists a Zariski open subset $W^0 \subset W$ and a
	Cartesian diagram
	\[
	\begin{tikzcd}
		\mathcal{E} \arrow[r, hook] \arrow[d, "f_{\mathcal{E}}"]
		& \mathcal{K} \arrow[d, "f"] \\
		W^0 \arrow[r, hook]
		& W
	\end{tikzcd}
	\]
	where $f_{\mathcal{E}}$ is smooth and proper, and
	$E_t := f_{\mathcal{E}}^{-1}(t) = E$ for all $t \in U$.
	
	Since smooth deformations of projective spaces are again projective spaces,
	we have $E_t \simeq E$ for every $t \in W^0$. Moreover, the class of $E_t$ in
	$\Pic(K_t)$ remains divisible by $2$, and the inverse image
	$\widetilde{E}_t$ of $E_t$ in the double cover $\widetilde{A}_t$ of $K_t$
	ramified along $E_t$ is isomorphic to $E$.
	
	Since the normal bundle of each component of $\widetilde{E}_t$ is isomorphic
	to $\mathcal{O}(-1)$, Nakano's contraction theorem (see \cite{Na}) yields a
	smooth complex space $A_t$ and a contraction morphism
	\[
	g : \, \widetilde{A}_t \,\longrightarrow \, A_t
	\]
	contracting exactly the components of $\widetilde{E}_t$.
	
	As this construction can be performed in families, we obtain a smooth
	projective morphism
	\[
	h : \mathcal{A}\, \longrightarrow \, W^0
	\]
	such that $A_t := h^{-1}(t)$ is the abelian variety whose associated Kummer
	manifold is $K_t$ for all $t \in U$. Since smooth deformations of complex tori
	are complex tori (see \cite[Theorem 4.6]{Ca}), it follows that $A_t$ is an abelian variety
	for every $t \in W^0$, and hence $K_t$ is the associated Kummer manifold.
\qed

\

\noindent We now use Ratner theory, as introduced by Verbitsky in the study of
hyperk\"ahler manifolds, to obtain global results on Kummer fourfolds in K3$^{[3]}$ type manifolds.

\subsection{The proof of Theorem \ref{thm:WeilK33}}\label{prf:ratner}
We use the notation from Proposition \ref{prop:defo}. 
The invariant sublattice of the rational
involution $\iota_{M}$ on $H^2(M,\ZZ)\cong v^\perp$ is (see Section~\ref{sub:invret})
\[
\Lambda\, :=\, U(2)^{\oplus 3} \oplus \langle -4 \rangle~.
\]
By \cite[Lemma~3.5]{BLjems}, this lattice is isomorphic to the
BBF lattice of $\widehat{M}$.
Therefore, the local Torelli theorem for locally trivial deformations of the
singular symplectic variety $\widehat{M}$ (see \cite[Proposition~4.7]{BLjems})
identifies $U$ with an Euclidean open subset of a five dimensional period domain
\[
U\;\hookrightarrow\;
\Omega_{\Lambda} :=
\left\{ [\sigma] \in \mathbb{P}(\Lambda \otimes_{\mathbb{Z}} \mathbb{C})
\;\middle|\;
\sigma \cdot \sigma = 0,\;
\sigma \cdot \overline{\sigma} > 0
\right\}.
\]
 
For any positive $d$, there is a primitive $\lambda\in U(2)\subset \Lambda$ with $(\lambda,\lambda)=4d$, and then
$$
\lambda^\perp\,\cong\,U(2)^{\oplus 2}\,\oplus\,\langle -4\rangle\,\oplus\,\langle -4d\rangle~,
\qquad \det(\lambda^\perp)\,=\,2^8d~.
$$ 
Hence, up to products with squares of integers, $\det(\lambda^\perp)$ can be any positive integer.

\

By Ratner's theorem (see \cite[Proposition~3.11]{BLjems} or 
\cite[Theorem~2.5]{Ve17}), the action of the orthogonal group $O(\Lambda)$ on the period
domain $\Omega_{\Lambda}$ has many dense orbits.
More precisely, the orbits are dense in the Euclidean topology for all points
corresponding to Hodge structures whose $(2,0)\oplus(0,2)$ part contains no
rational vectors.
This implies that every K3 type Hodge structure on $\Lambda$ without rational classes in the
$(2, 0) \oplus (0, 2)$ part is represented in every open subset of $\Omega_{\Lambda}$, 
and in particular in $U$. 

Therefore any rank six sublattice $\lambda^\perp\subset \Lambda$ as above occurs as the 
transcendental lattice of a deformation $\widehat{M}_t$ for some $t\in U$. 
Notice that $\widehat{M}_t$ is projective with $\Pic(\widehat{M}_t)=\ZZ\lambda\oplus BW$.
The induced Hodge structure on $\lambda^\perp$ has Hodge numbers $(1,4,1)$. 

The corresponding deformation of $Z$ is a Kummer variety $Z_t=K(B)\subset M_t$ for some abelian
fourfold $B$ and, as in the proof of Theorem \ref{thm:c2TX} given in Section \ref{sec:prf03}, 
the pull back to $B$ of $\lambda^\perp$ defines a sub Hodge structure $T$ of rank 6, 
with the same Hodge numbers, in $H^2(B,\ZZ)$. 
Therefore $B$ is of Weil type by Lemma \ref{lem:noKS}, for some field $K$, and has trivial discriminant. 
As in the proof of that lemma, the polarization on $\lambda^\perp$ given by the restriction of the BBF form to 
is unique up to scalar multiple. The same is true for $T$ and thus these quadratic forms on $\lambda^\perp$ and $T$ are proportional. As the ranks of these lattices is six, this implies that their determinants are the same up to squares, thus $K= \QQ(\sqrt{-d})$ again by Lemma \ref{lem:noKS}, and any $d>0$ occurs for a suitable choice of $\lambda$.

From the four dimensional family of Hodge structures on $\lambda^\perp$ obtained by varying $t\in\PP(\lambda^\perp)_\CC\cap \Omega_\Lambda$  we thus obtain a complete family of Weil type
abelian fourfolds,
with field $\QQ(\sqrt{-d})$ and discriminant one such that the Kummer variety of the general member
is a submanifold of a deformation of $M$. \qed

\

\subsection{An alternative proof of Thm.\ \ref{thm:WeilK33}}
We now present an alternative proof, still based on Proposition~\ref{prop:defo}, that avoids the use of Ratner's ergodic-theoretic results. To this end, we first prove the following lemma and its corollary, which provide the necessary density results.

\

\subsection{Lemma} \label{lemma:dens}

Let $\Gamma'$ be a lattice, and set
\[
\Gamma := U^{\oplus 3} \oplus \Gamma'.
\]
Let $d$ and $e$ be positive integers such that $\Gamma$ contains an element
of square $2d$ and divisibility $e$.

Let $\Omega_{\Gamma}$ denote the associated period domain, i.e.\
\[
\Omega_{\Gamma} :=
\left\{ [\sigma] \in \mathbb{P}(\Gamma \otimes \mathbb{C})
\;\middle|\; \sigma^2 = 0 \text{ and } \sigma \cdot \overline{\sigma} > 0 \right\}.
\]

Then, for every non-empty Euclidean open subset $U \subset \Omega_{\Gamma}$,
there exist $[\sigma_0] \in U$ and $\gamma \in \Gamma$ such that
\[
\sigma_0 \cdot \gamma = 0,\quad \gamma^2 = 2d,
\quad \text{and} \quad {\rm div}(\gamma) = e.
\]

\ts
Let $\sigma \in \Omega_{\Gamma}$ and denote by
\[
\sigma^{\perp} \subset \Gamma \otimes \mathbb{C}
\]
the orthogonal hyperplane, and by
\[
\mathbb{P}(\sigma^{\perp}) \subset \mathbb{P}(\Gamma \otimes \mathbb{C})
\]
its projectivization, which is the projective tangent space to $\Omega_\Gamma$.

Since $U$ is open, the union
\[
V := \bigcup_{\sigma \in U} \mathbb{P}(\sigma^{\perp})
\subset \mathbb{P}(\Gamma \otimes \mathbb{C})
\]
is a non-empty Euclidean open subset.

Indeed, let $Q\subset \mathbb{P}(\Gamma\otimes \mathbb{C})$ be the Zariski closure of $\Omega_{\Gamma}$, i.e. the complex quadric defined by the quadratic form on $\Gamma$. Let
\[
J\subset Q\times \mathbb{P}(\Gamma\otimes \mathbb{C})
\]
be the incidence variety of tangent hyperplanes to $Q$, and let
\[
p_1:J\,\longrightarrow\, Q,\qquad p_2:J\,\longrightarrow\, \mathbb{P}(\Gamma\otimes \mathbb{C})
\]
be the natural projections. Then
$V=p_2\bigl(p_1^{-1}(U)\bigr)$
and, since $p_2$ is flat, hence open, it follows that $V$ is an open subset of $\mathbb{P}(\Gamma\otimes \mathbb{C})$.

For every $\sigma \in \Omega_{\Gamma}$, the intersection
\[
\sigma^{\perp} \cap \overline{\sigma}^{\perp}
\]
gives a non-degenerate, real, vector subspace of $\Gamma \otimes \mathbb{R}$ on which 
the quadratic form is not definite. In particular, it contains isotropic real vectors.
It follows that $V$ contains a non-empty open subset of the real quadric
\[
Q_{\mathbb{R}} := \{[\alpha] \in \mathbb{P}(\Gamma \otimes \mathbb{R})
\mid \alpha^2 = 0\}.
\]

Since $\Gamma $ contains isotropic vectors, the quadric
$Q_{\mathbb{R}}$ admits rational points, and the set of rational points
$Q_{\mathbb{Q}}$ is dense in $Q_{\mathbb{R}}$. 

More precisely, we show that rational
points represented by primitive vectors of divisibility $1$ are also dense.
Let $\{e_i,f_i\}_{i=1}^3$ be a standard basis of $\Gamma$ such that
$e_i \cdot e_j = f_i \cdot f_j = 0$ and $e_i \cdot f_j = \delta_{ij}$.
Let $\alpha$ be a primitive integral representative of a point of
$Q_{\mathbb{Q}}$. By the Eichler criterion (see \cite[Section~10]{Ei} or
\cite[Lemma~3.5]{GHS}), since $U^{\oplus 3}$ is unimodular and contains two
orthogonal isotropic planes, the automorphism group of $U^{\oplus 3}$ acts
transitively on primitive isotropic vectors of a given divisibility.
Hence we may assume that the component of $\alpha$ in $U^{\oplus 3}$ is of
the form $a e_3 + b f_3$.

It follows that, for every integer $n$, the vector $e_1 + n\alpha$ is
primitive, isotropic, and of divisibility $1$, and its projective class
approaches $[\alpha]$. Therefore, primitive isotropic classes of
divisibility $1$ are dense in $Q_{\mathbb{R}}$.

In view of this density, It thus suffices to prove that, for every primitive $\alpha \in \Gamma$
with $\alpha^2 = 0$ and ${\rm div}(\alpha)=1$, the point $[\alpha]$ is a limit
of classes of primitive elements in $\Gamma$ of square $2d$ and
divisibility $e$. By the Eichler criterion, applied to $\Gamma$, there
exists an isometry of $\Gamma$ sending $\alpha$ to $e_1$, so we may assume
$\alpha = e_1$.

Let $\gamma \in \langle e_3, f_3 \rangle \oplus \Gamma'$ be an element of
square $2d$ and divisibility $e$, and set
\[
\alpha_n := n e \;e_1 + \gamma.
\]
Then $\alpha_n^2 = 2d$, ${\rm div}(\alpha_n) = {\rm div}(\gamma) = e$, and
\[
\lim_{n \rightarrow \infty} [\alpha_n] = [\alpha].
\]
\qed

\subsection{Corollary}

\begin{enumerate}
\item For $X$ a K3 surface or a complex torus of dimension $2$, the locus in
its Kuranishi space parametrizing projective deformations with a primitive
polarization of degree $2d$ is dense in the Euclidean topology.

\item For a hyperk\"ahler manifold $X$ of the known deformation types, i.e.\
K3$^{[n]}$ type, Kummer type, or one
of O'Grady's examples, let $d,e>0$ be such that $X$ admits a class of square $2d$
and divisibility $e$ in $H^2(X,\ZZ)$.
Then the locus in the Kuranishi space of $X$ parametrizing deformations with a primitive
polarization of degree $2d$ and divisibility $e$ is dense in the Euclidean
topology.
\end{enumerate}
\ts
The 
second cohomology group of $X$ satisfies the assumptions on the
lattice $\Gamma$ in Lemma~\ref{lemma:dens}. By the local Torelli theorem,
the Kuranishi space of $X$ can be identified with a Euclidean open subset
$U$ of the period domain $\Omega_{\Gamma}$.

By Lemma~\ref{lemma:dens}, the locus in the Kuranishi space of $X$
parametrizing deformations admitting an integral $(1,1)$--class of square
$2d$ and divisibility $e$ is dense. This property still holds if one further
imposes that the Néron Severi group be cyclic.

By the projectivity criterion, such deformations are projective, and an ample generator of their Néron Severi group defines the desired polarization.
\qed

\vspace{\baselineskip}\noindent{\bf Proof of Theorem 0.2.}$\;\;$ We use the notation from Proposition \ref{prop:defo}.
As shown in Section~\ref{sub:invret}, the invariant sublattice of the rational
involution $\iota_{M}$ on $H^2(M,\ZZ)\cong v^\perp$ is
\[
\Lambda\, :=\, U(2)^{\oplus 3} \oplus \langle -4 \rangle~.
\]

The local Torelli theorem for hyperk\"ahler manifolds of $K3^{[3]}$--type 
identifies $U$ with a Euclidean open subset of the five dimensional period domain
\[
U\;\hookrightarrow\;
\Omega_{\Lambda} :=
\left\{ [\sigma] \in \mathbb{P}(\Lambda \otimes_{\mathbb{Z}} \mathbb{C})
\;\middle|\;
\sigma \cdot \sigma = 0,\;
\sigma \cdot \overline{\sigma} > 0
\right\}.
\] for the Hodge locus of $BW$.
 
For any positive $d$, there is a primitive $\lambda\in U(2)\subset \Lambda$ with $(\lambda,\lambda)=4d$, and then
$$
\lambda^\perp\,\cong\,U(2)^{\oplus 2}\,\oplus\,\langle -4\rangle\,\oplus\,\langle -4d\rangle~,
\qquad \det(\lambda^\perp)\,=\,2^8d~.
$$ 
Hence, up to products with squares of integers, $\det(\lambda^\perp)$ can be any positive integer.

By  Lemma \ref{lemma:dens}, $\lambda^{\perp}$ is isomorphic  to the  transcendental lattice of some hyperk\"ahler manifold $M_t$ of $K3^{[3]}$--type representing a point $t$ of $U$, moreover $(\lambda,\lambda)>0$ and by Huybrechts criterion $M_t$ is projective and, shrinking $U$ if necessary, we  may assume that  the Hodge locus $U_{\lambda\oplus BW}\subset U$ of $\lambda\oplus BW$ is a four dimensional complex manifold  parametrizing hyperk\"ahler manifolds belonging to a fixed component $Y$ of the moduli space of polarized $K3^{[3]}$ manifolds.

This component $Y$ is algebraic (see \cite[Theorem~1.5]{GHS}), and the same
holds for the Hodge locus 
\[
Y_{\lambda \oplus BW} \subset Y
\]
of the sublattice $\lambda \oplus BW$: indeed  this is precisely the locus where the
$256$ projective three--spaces corresponding to the primitive MBM classes of
$BW$, together with the effective divisorial class $\lambda$, deform.

The modular map sends $U$ to  $Y_{\lambda \oplus BW}$: hence 
there exists a nonempty Euclidean open subset of $Y_{\lambda \oplus BW}$ parametrizing
hyperk\"ahler manifolds of $K3^{[3]}$--type containing a Kummer fourfold as a
smooth subvariety.

By the algebraicity of the relative Hilbert scheme together with
Lemma~\ref{a}(2), one obtains a Zariski open subset
\[
Y_{\lambda \oplus BW}^0 \subset Y_{\lambda \oplus BW}
\]
having the same property.

Since the image of the period map for the moduli space of polarized
hyperk\"ahler manifolds is the complement of a finite union of divisors in the
period domain (see \cite[Theorem~B.1]{De}), it follows that a  general
K3 type Hodge structure is isomorphic to the Hodge structure on the
transcendental lattice of a $K3^{[3]}$ type manifold in $Y_{\lambda \oplus BW}^0$ containing a
Kummer fourfold as a smooth subvariety.
Finally, by Lemma \ref{lem:noKS}, the associated abelian variety is a Weil abelian variety  for a quadratic imaginary field   $K$ and trivial discriminant and, by uniqueness of the polarization for a generic K3 type Hodge structure, $K=\sqrt{-d}$. Conversely, since there exists a single four dimensional family parametrizing, up to isogeny, all abelian Weil fourfolds with fixed quadratic imaginary field and trivial discriminant, it follows that, up to isogeny,  the Kummer fourfold of a general  abelian fourfold of Weil type with trivial discrimant appears as a smooth subvariety of  a $K3^{[3]}$ type manifold.
\qed

\

\subsection{Remark}
The five dimensional family of Kummer fourfolds over $\Omega_\Lambda$
(notation as in Section \ref{prf:ratner})
is very similar to the family considered by Markman in \cite{Markman23}.
We briefly recall that family and, for simplicity, we follow the exposition in \cite{vG23}.

Let $S^+:=U^{\oplus 4}$, a lattice of rank $8$ and signature $(4,4)$. Let $Q\subset \PP(S^+_\CC)$ be the six dimensional quadric defined by this quadratic form in the projectivization of the complex vector space
$S^+_\CC:=S^+\otimes_\ZZ\CC$. The quadric $Q$ parametrizes certain $4$--dimensional complex
subspaces in $V_\CC\cong \CC^8$ that are isotropic in the complexification of a lattice $V\cong\ZZ^8$.

Let $\Omega\subset Q$ be the (Euclidean) open subset of $Q$ whose points are the $\ell\in Q$ with
$(\ell,\overline{\ell})>0$. Each $\ell\in\Omega$ defines a $4$--dimensional complex torus $\cT_\ell$,
it is $V_\CC/(Z_\ell+V)$, where $Z_\ell\subset V_\CC$ is the subspace defined by $\ell$,
\cite[\S 1.12]{vG23}.
There is a natural isomorphism $H^4(\cT_\ell,\ZZ)=\wedge^4V$.

Fix a non-zero $s\in S^+$, then $s$ defines a class $c_s\in H^4(\cT_\ell,\ZZ)$, called a Cayley class.
For any $\ell\in \Omega$ with $(s,\ell)=0$,  the Cayley class is a Hodge class, i.e.\
$c_s\in H^{2,2}(\cT_\ell,\ZZ)$, \cite[Proposition 2.7]{vG23}.
Thus the five dimensional domain $\Omega\cap s^\perp$ parametrizes
complex tori with a Hodge class in codimension two.

For $h\in S^+$ such that $(h,s)=0$ and $\langle h,s\rangle\subset S^+$ is a positive definite rank two sublattice, define $d:=(h,h)(s,s)$. Then for $\ell\in \langle h,s\rangle^\perp\cap \Omega$,
a four dimensional domain, the complex torus $\cT_\ell$ is an abelian fourfold of Weil type, with trivial discriminant and field $K=\QQ(\sqrt{-d})$. Moreover, for general such $\ell$,
the Cayley class $c_s$ is not an intersection of divisor classes, \cite[Theorem 4.6]{vG23}.

The family of the $\cT_\ell$ over $\Omega\cap s^\perp$ has properties that are very similar
to those of the five dimensional family of  complex tori $B_t$ over $\Omega_\Lambda$ that we constructed.
The role of the Cayley class $c_s$ is now played by the pull-back to $B_t$ of the second Chern class
of the tangent bundle $c_2(\cT_M)$ of the hyperk\"ahler sixfold $M=M_t$ for general $t\in \Omega_\Lambda$. The choice of $h$  corresponds to the choice $\lambda\in \Lambda$, with $(\lambda,\lambda)=4d$; both determine complete 4--dimensional families of abelian fourfolds of Weil type with field $K$.

\


\begin{thebibliography}{PQ}

\bibitem[AV1]{AVimrn} E.\ Amerik, M.\ Verbitsky,
{\it Rational Curves on Hyperk\"ahler Manifolds},
International Mathematics Research Notices, {\bf 23} (2015), 1300--1345.

\bibitem[AV2]{AVsel} E.\ Amerik, M.\ Verbitsky,
{\it Contraction centers in families of hyperk\"ahler manifolds},
Selecta Mathematica {\bf 27} (2021), no 60, pp 26.


\bibitem[AV3]{AV21'} E.\ Amerik, M.\ Verbitsky, {\it Rational Curves and MBM Classes on Hyperkähler Manifolds: A Survey}. In:
Rationality of Varieties (2021), 75--96, Progr. Math.  342, Birkh\''auser/Springer, Cham 2021. 


\bibitem[AV4]{Amerik_V} E.\ Amerik, M.\ Verbitsky,
{\it MBM classes and contraction loci on low-dimensional
hyperk\"ahler manifolds of K3$^{[n]}$ type},
Algebraic Geometry {\bf 9} (2022), 252--265.

\bibitem[Be]{Be} A.\ Beauville,
{\it Symplectic singularities},
Invent.\ Math.\ {\bf139} (2000), 541--549.

  
\bibitem[BaL]{BLjems} B.\ Bakker, C.\ Lehn,
{\it A global Torelli theorem for singular symplectic varieties},
J.\ Eur.\ Math.\ Soc.\ {\bf 23} (2021), 949--994.


\bibitem[BiL]{LB} C.\ Birkenhake, H.\ Lange,
{Complex abelian varieties},
Second edition, Springer Verlag 2004.

\bibitem[Bo]{Bo} S.\ Boucksom,
{\it Divisorial Zariski decompositions on compact complex manifolds},
Ann.\ Sci.\ \'Ecole Norm.\ Sup.\ (4) {\bf 37} (2004), 45--76.

\bibitem[C]{Ca} F.\ Catanese, {\it Deformation types of real and complex manifolds}, in: Contemporary Trends in Algebraic Geometry and Algebraic Topology (Proceedings of the Chern-Chow Memorial Conference), Nankai Tracts in Mathematics 5, (2002) 193--236.

\bibitem[Ca]{Cam} F.\ Campana,
{\it  Isotrivialit\'e de certaines familles k\"ahl\'eriennes de vari\'et\'s non projectives},
Mathematische Zeitschrift,  {\bf 252} (2006), 147--156.

\bibitem[CF]{Cacciatori_F} S.\ Cacciatori, S.\ Filippini,
{\it The $E^3/\ZZ_3$ orbifold, mirror symmetry, and Hodge structures of Calabi-Yau type},
Journal of Geometry and Physics {\bf 138} (2019), 70--89.

\bibitem[Ch]{Char} F. Charles,
{\it Two results on the Hodge structure of complex tori},
Math.\ Z. {\bf 300} (2022),  3623--3643.

\bibitem[CS]{Conway_S} J.H.\ Conway, N.J.A.\ Sloane,
Sphere Packings, Lattices and Groups,
Grundlehren der Mathematischen Wissenschaften, vol. 290 (2nd ed.), Berlin, New York: Springer-Verlag 1988.

\bibitem[De]{De} O.\ Debarre,
{\it Hyperk\"ahler manifolds}.
Milan J. Math.  {\bf 90} (2022), 305--387.

\bibitem[Do]{Dolgachev} I.\ Dolgachev,
Algebraic Geometry. Cambridge University Press 2012.

\bibitem[Dr]{Dr} S.\ Druel,
{\it Quelques remarques sur la d\'ecomposition de Zariski divisorielle sur les vari\'et\'es dont la premi\`ere classe de Chern est nulle},
Math.\ Z.\ {\bf 267} (2011), 413--423.

\bibitem[E]{Ei}  M.\ Eichler, {\it Quadratische Formen und orthogonale Gruppen},
Grundlehren der mathematischen Wissenschaften 63. SpringerVerlag, Berlin-New York, 1974.

\bibitem[Fl]{Fle} H.\ Flenner, {\it Extendability of differential forms on non-isolated singularities},
Invent.\ Math.\  {\bf 94} (1988), 317--326.

\bibitem[FF25]{Floccari_F} S.\ Floccari, L.\ Fu,
{\it The Hodge conjecture for Weil fourfolds with discriminant 1 via singular OG6-varieties},
J.\ Math.\ Pures Appl.\  {\bf 210} (2026), Paper No. 103876, 17 pp.


\bibitem[FF26]{Floccari_F26} S.\ Floccari, L.\ Fu,
{\it The hyper-Kummer construction},
arXiv:2607.07528.

\bibitem[Fu]{Fu} O.\ Fujino, {\it Notes on rational chain connectedness}
arXiv:2602.19415 2026.


\bibitem[FH]{Fulton_H} W.\ Fulton, J.\ Harris,
Representation Theory. GTM 129, Springer-Verlag 1991. 

\bibitem[vG1]{vG_hcav} B.\ van Geemen,
{\it An introduction to the Hodge Conjecture for abelian varieties},
in: Algebraic cycles and Hodge theory (Torino, 1993), 233--252, Lecture Notes in Math., 1594, Springer, Berlin, 1994.



\bibitem[vG2]{vG23} B.\ van Geemen,
{\it Fourfolds of Weil type and the spinor map},
Expositiones Mathematicae {\bf 41} (2023), 418--447.


\bibitem[GHS]{GHS} V.\ Gritsenko, K.\ Hulek,  G.\ K.\ Sankaran,  {\it Moduli spaces of irreducible symplectic manifolds}, Compositio Mathematica, {\bf 146} (2010),  404--434.


\bibitem[HT]{Hassett_T_intnum} B.\ Hassett, Y.\ Tschinkel,
{\it Intersection numbers of extremal rays on holomorphic symplectic varieties},
Asian J.\ Math.\ {\bf 14} (2010), 303--322.

\bibitem[HHT]{Harvey_HT} D.\ Harvey, B.\ Hassett, Y.\ Tschinkel,
{\it Characterizing Projective Spaces on Deformations
of Hilbert Schemes of K3 Surfaces},
Communications on Pure and Applied Mathematics, Vol. LXV, 0264–0286 (2012).


\bibitem[K]{Ka} Y.\ Kawamata,
{\it Flops Connect Minimal Models},
Publ.\ Res.\ Inst.\ Math.\ Sci.\ {\bf 44} (2008),  419–423.

\bibitem[L]{Lombardo} G.\ Lombardo,
{\it Abelian varieties of Weil type and Kuga-Satake varieties}, 
Tohoku Math.\ J.\ {\bf 53} (2001), 453--466.

\bibitem[Ma1]{Msurv} E.\ Markman,
{\it A survey of Torelli and monodromy results for holomorphic-symplectic varieties}
in: Complex and differential geometry, 257--322,
Springer Proc. Math., 8, Springer, Heidelberg, 2011.

\bibitem[Ma2]{MaSPE} E.\ Markman,
{\it Prime exceptional divisors on holomorphic symplectic varieties and monodromy-reflections},
Kyoto Journal of Mathematics, {\bf 53} (2013), 345-403.

\bibitem[Ma3]{Ma14} E.\ Markman.
{\it Lagrangian fibrations of holomorphic-symplectic varieties of K3$^{[n]}$-type},
in Anne Fr\"uhbis-Kr\"uger, Remke Nanne Kloosterman, and Matthias Sch\"utt, editors,
Algebraic and Complex Geometry, 241--283. Springer International Publishing (2014).

\bibitem[Ma4]{Markman23} E.\ Markman,
{\it The monodromy of generalized Kummer varieties and algebraic cycles on their intermediate Jacobians},
J.\ Eur.\ Math.\ Soc.\ {\bf 25} (2023), 231--321.

\bibitem[Ma5]{Markman25} E.\ Markman,
{\it Cycles on abelian 2n-folds of Weil type from secant sheaves on abelian n-folds}.
arXiv:2502.03415.

\bibitem[Mat]{Mat} D.\ Matsushita,
{\it On Deformations of Lagrangian Fibrations},
in:  K3 Surfaces and Their Moduli, Vol. 315, 2016,  237--243. Birkhäuser, Cham.


\bibitem[MZ]{MZ} D.\ Matsushita, D.-Q.\ Zhang,
{\it Zariski F-decomposition and Lagrangian fibration on hyperk\"ahler manifolds},
Math.\ Res.\ Lett.\ {\bf 20} (2013), no. 5, 951--959.
 
\bibitem[Mo1]{Mongardi15} G.\ Mongardi,
{\it A note on the K\"ahler and Mori cones of hyperk\"ahler manifolds},
Asian J.\ Math.\ {\bf 19} (2015), 583--591.

\bibitem[Mo2]{Mongardi} G.\ Mongardi,
{\it Towards a classification of symplectic automorphisms
on manifolds of K3$^{[n]}$ type},
Math.\ Z.\ {\bf 282} (2016), 651--662.

\bibitem[MRS18]{Mongardi_RS18} G.\ Mongardi, A.\ Rapagnetta, G. Sacc\`a,
{\it The Hodge diamond of O’Grady’s six-dimensional example}, 
Compositio Math.\ {\bf 154} (2018),  984--1013.

\bibitem[Na]{Na} S.\ Nakano,
{\it On the inverse of monoidal transformations},
Publications of the Research Institute for Mathematical Sciences, {\bf 10}(3) (1974), 837--851.

\bibitem[N1]{Nikulin} V.V.\ Nikulin,
{\it On Kummer surfaces},
Izv.\ Akad.\ Nauk SSSR Ser.\ Mat.\ {\bf 39} (1975) 278--293,
Math.\ USSR Izv.\ {\bf 9} (1975), 261--275 (1976).

\bibitem[N2]{NikulinL} V.V.\ Nikulin,
{\it Integer symmetric bilinear forms and some of their geometric applications},
Izv.\ Akad.\ Nauk SSSR Ser.\ Mat.\ {\bf 43} (1979), 111--177,
Math.\ USSR Izv.\ {\bf 14} (1979), 103--167 (1980).

\bibitem[Po]{Po} G.\ Pourcin,
{\it Th\'eor\`eme de Douady au-dessus de S},
Annali della Scuola Normale Superiore di Pisa, Classe di Scienze 3
série, tome 23 (1969),
no 3, 451--459.


\bibitem[O'G]{O'Grady} K.\ O'Grady,
{\it Irreducible symplectic 4-folds numerically equivalent to (K3)$^{[2]}$},
Commun.\ Contemp.\ Math.\ {\bf 10} (2008), 553--608.

\bibitem[R]{Rizzo} F.\ Rizzo,
{\it On the fixed locus of the antisymplectic involution of an EPW cube},
arXiv:2505.15717.


\bibitem[Sc]{Sch71} M.\ Schlessinger,
{\it Rigidity of quotient singularities},
Invent.\  Math.\ {\bf 14}  (1971), 17--26.


\bibitem[SoV] {SoV}A.\ Soldatenkov,  M.\ Verbitsky,
{\it The abundance and SYZ conjectures in families of hyperk\"ahler manifolds}
epiga:14370 - \'Epijournal de G\'eométrie Alg\'ebrique, 9 décembre 2025, Volume 9.

\bibitem[Vo]{Vlagr} C.\ Voisin,
{\it Sur la stabilit\'e des sous-vari\'et\'es lagrangiennes des vari\'et\'es
symplectiques holomorphes},
in: Complex projective geometry, 294--303. Cambridge Univ. Press, 1992.

\bibitem[Ve]{Ve17} M.\ Verbitsky,
{\it Ergodic complex structures on hyperk\"ahler manifolds: an erratum}, arXiv:1708.05802.

\bibitem[W]{Weil} A.\ Weil, {\it Abelian varieties and the Hodge ring},
in: Collected Papers, Vol. III, 421--429. Springer Verlag (1980).

\bibitem[WW]{WW} J.\ Wierzba, J.\ A.\ Wi\'sniewski,
{\it Small contractions and symplectic 4-folds},
Duke Mathematical Journal, {\bf 120} (2003), 65--95.

\bibitem[Wl]{Wl} J.\ Włodarczyk, 
{\it Resolution of singularities of analytic spaces}, in Proceedings of G\"okova Geometry-Topology Conference 2008, Int. Press, Somerville, MA, 2009,  13--47.

\

\end{thebibliography}
\end{document}